\documentclass[11pt,a4paper,reqno]{amsart}

\usepackage[margin=1.0in]{geometry}
\usepackage{eucal,color,graphicx,subfigure,mathrsfs}
\usepackage{amsmath}
\usepackage{amssymb}
\usepackage{stmaryrd}
\usepackage{hyperref}
\usepackage{multicol}
\usepackage{multirow}
\usepackage{booktabs}
\usepackage{subfigure}
\usepackage{cite}
\allowdisplaybreaks

\numberwithin{equation}{section}

\newtheorem{theorem}{Theorem}[section]
\newtheorem{lemma}[theorem]{Lemma}

\newtheorem{assumption}{Assumption}
\newtheorem{assumptionalt}{Assumption}[assumption]

\newenvironment{assumptionp}[1]{
	\renewcommand\theassumptionalt{#1}
	\assumptionalt
}{\endassumptionalt}

\theoremstyle{definition}
\newtheorem{definition}[theorem]{Definition}
\newtheorem{example}[theorem]{Example}

\theoremstyle{remark}
\newtheorem{remark}[theorem]{Remark}

\numberwithin{equation}{section}

\newcommand{\dt}{\frac{{\rm d}}{{\rm d}t}} 

\newcommand{\BF}{{\boldsymbol{F}}}

\newcommand{\BPhi}{{\boldsymbol{\Phi}}}
\newcommand{\BPsi}{{\boldsymbol{\Psi}}}

\newcommand{\CP}{{\mathcal P}}
\newcommand{\CQ}{{\mathcal Q}}

\newcommand{\CT}{{\mathcal T}}

\newcommand{\CHu}{{\check{u}}}
\newcommand{\CHv}{{\check{v}}}
\newcommand{\CHchi}{{\check{\chi}}}

\newcommand{\CHpsi}{{\check{\psi}}}
\newcommand{\CHCP}{{\check{\mathcal{P}}}}
\newcommand{\CHCT}{{\check{\mathcal{T}}}}
\newcommand{\CHE}{{\check{E}}}
\newcommand{\CHI}{{\check{I}}}
\newcommand{\CHK}{{\check{K}}}

\newcommand{\CHS}{{\check{S}}}

\newcommand{\vertiii}[1]{{\left\vert\kern-0.25ex\left\vert\kern-0.25ex\left\vert #1 
    \right\vert\kern-0.25ex\right\vert\kern-0.25ex\right\vert}}

\newcommand{\R}{\mathbb{R}}

\def\d{{\mathrm d}}

\begin{document}

\title[]{Convergence of arbitrary Lagrangian-Eulerian\\ second-order projection method for\\ the Stokes equations on an evolving domain}

\author[]{\,\,Qiqi Rao}
\address{Department of Mathematics, City University of Hong Kong,
83 Tat Chee Avenue, Kowloon, Hong Kong, China}
\email{qi-qi.rao@connect.polyu.hk}

\author[]{Jilu Wang}
\address{
Department of Mathematics, School of Sciences, Harbin Institute of Technology Shenzhen.} 
\email{wangjilu@hit.edu.cn}

\author[]{\,\,Yupei Xie}
 \address{Department of Applied Mathematics, The Hong Kong Polytechnic University, Hung Hom, Hong Kong, China} 
\email{yu-pei.xie@connect.polyu.hk}

\begin{abstract} 
The numerical solution of the Stokes equations on an evolving domain with a moving boundary is studied based on the arbitrary Lagrangian-Eulerian finite element method and a second-order projection method along the trajectories of the evolving mesh for decoupling the unknown solutions of velocity and pressure. 
The error of the semidiscrete arbitrary Lagrangian-Eulerian method is shown to be $O(h^{r+1})$ for the Taylor--Hood finite elements of degree $r\ge 2$, using Nitsche's duality argument adapted to an evolving mesh, by proving that the material derivative and the Stokes--Ritz projection commute up to terms which have optimal-order convergence in the $L^2$ norm. Additionally, the error of the fully discrete finite element method, with a second-order projection method along the trajectories of the evolving mesh, is shown to be $O(\ln(1/\tau+1)\tau^{2}+\ln(1/h+1)h^{r+1})$ in the discrete $L^\infty(0,T; L^2)$ norm using newly developed energy techniques and backward parabolic duality arguments that are applicable to the Stokes equations with an evolving mesh. To maintain consistency between the notations of the numerical scheme in a moving domain and those in a fixed domain, we introduce the equivalence class of finite element spaces across time levels. Numerical examples are provided to support the theoretical analysis and to illustrate the performance of the method in simulating Navier--Stokes flow in a domain with a rotating propeller. 

\end{abstract}



\setlength\abovedisplayskip{3.5pt}
\setlength\belowdisplayskip{3.5pt}

\maketitle
\section{Introduction}
The Stokes and Navier--Stokes equations are partial differential equations (PDEs) widely used to describe the motion of viscous fluids such as water and air. Solving the Stokes and Navier--Stokes equations is a critical area of research in fluid dynamics, particularly when the domain is not fixed, such as in moving boundary/interface or fluid-structure interaction problems. The inclusion of such a dynamic domain introduces an additional layer of intricacy to the problem. 

This article concerns the numerical solution of the Stokes equations in a time-dependent domain $ \Omega (t) \subset \R^d$ with $d\in \{2,3\}$, i.e., 
\begin{subequations}\label{main-equation-stokes}
	\begin{align}
		\partial_tu-\Delta u+\nabla p =f \quad & \text{in } \quad \bigcup_{t\in (0,T]}\Omega(t)\times \{t\},\label{main-equation-stokes-a}\\
		\nabla\cdot u=0 \quad &  \text{in }\quad \bigcup_{t\in (0,T]}\Omega(t)\times \{t\},\\
		 u=0 \quad & \text{on } \quad \bigcup_{t\in (0,T]}\partial\Omega(t)\times \{t\},\\
		 u=u_0 \quad & \text{on} \quad\quad \Omega^0 = \Omega(0) .
	\end{align}
\end{subequations}
where the domain $ \Omega (t)$ has a smooth boundary $\Gamma(t) = \partial \Omega(t)$ which moves under a given smooth vector field $w(\cdot,t)$. For simplicity, we assume that the vector field $w$ has a smooth extension (which we do not need to know explicitly) to the entire space $\R^d$ and generates a smooth flow map $\Phi(\cdot,t)$. The equation also includes a source term $f$, a given smooth function that depends on both space and time variables. To ensure uniqueness of the solutions, we assume that $p(\cdot, t)\in L^2_0(\Omega(t))$, which is the space of functions $p$ in $L^2(\Omega(t))$ such that $\int_{\Omega(t)}p\, \d x = 0$.

The arbitrary Lagrangian-Eulerian (ALE) formulation is widely used to handle the complexities arising from domain evolution by allowing the mesh to move according to an ALE mapping, such as the interpolation $\Phi_h$ of $\Phi$, to fit the evolving domain. To employ the ALE formulation, one can define the material derivative of the solution $u$ with respect to the velocity field $w$ as
\begin{align}\label{u-material-derivative}
	D_t u(x,t) := \frac{\d}{\d t}u(\Phi(\xi,t),t) = \partial_t u + w \cdot \nabla u \,\,\text{at}\,\, x = \Phi(\xi,t)\in \Omega(t)\;\text{for}\;\xi \in \Omega^{0}.
\end{align}
Using this definition of material derivative, the first two equations in \eqref{main-equation-stokes} can be rewritten as 
\begin{subequations}\label{ALE-formulation}
\begin{align}
	D_t u - w\cdot \nabla u - \Delta u + \nabla p &= f  \\
	\nabla\cdot u &=0 ,
\end{align}
\end{subequations}
and the ALE method can be employed to discretize the material derivative $ D_{t}u $ along the characteristic lines of the evolving mesh. 

In an early investigation of ALE methods, Formaggia \& Nobile \cite{Formg1999} provided stability results for two different ALE finite element schemes. Subsequently, Gastaldi \cite{Gas2001} established a priori error estimates of ALE finite element methods (FEMs) for parabolic equations, illustrating that a piecewise linear element can yield an $ L^{2} $ error of order $O(h)$ when the mesh size $ h $ is sufficiently small. In a related study \cite{Nobie2001}, Nobile obtained an error estimate of $O(h^k)$ in the $L^2$ norm for spatially semidiscrete ALE finite element schemes, with $ k $ denoting the degree of the piecewise polynomials utilized. The stability of time-stepping schemes in the context of ALE formulations, such as implicit Euler, Crank--Nicolson, and backward differentiation formulae (BDF), were proved in \cite{Boffi2004} and \cite{ForNo2004}. Under specific generalized compatibility conditions and step size restrictions, these investigations yielded $L^2$ error estimates of $O(\tau^s + h^k)$, where $s = 1,2$ corresponds to the order of the time schemes and $ k $ denotes the degree of the finite element space employed. Moreover, Badia and Codina \cite{Badia2006} obtained $ L^{2} $ error bounds of $O(\tau^s + \tau^{-1/2}h^{k+1})$ for $s = 1, 2$ for fully discretized ALE methods that employ BDF in time and FEM in space. These sub-optimal error bounds were obtained when the mesh dependent stabilization parameter appearing in fully discrete scheme is as small as the time step size. 

Optimal convergence of $O(h^{r+1})$ in the $L^{\infty}(0,T;L^2)$ norm of ALE semidiscrete FEM for diffusion equations in a bulk domain with a moving boundary was established by Gawlik \& Lew in \cite{GawE2015} for finite element schemes of degree $r \ge 1$. We also refer to \cite{ElliottUFEM2021} and \cite{Edelmann-2022} for a unified framework of ALE evolving FEMs and an ALE method with harmonically evolving mesh, respectively. Optimal-order $H^1$ convergence of the ALE FEM for PDEs coupling boundary evolution arising from shape optimziation problems was proved in \cite{GongLiRao-2023}. These results were established for high-order curved evolving mesh. 
Optimal convergence of $O(h^{r+1})$ in the $L^{\infty}(0,T;L^2)$ norm, with flat evolving simplices in the interior and curved simplices exclusively on the boundary, was prove in \cite{buyang22} for the ALE semidiscrete FEM utilizing a standard iso-parametric element of degree $r$ in \cite{Lenoir1986}. 

In addition to the ALE spatial discretizations mentioned above, the stability and error estimates of discontinuous Galerkin (dG) semi-discretizations in time for diffusion equations in a moving domain using ALE formulations were established in \cite{Bonito2013} and \cite{BonKyerror2013}, respectively. The ALE methods for PDEs in bulk domains \cite{GongLiRao-2023} are also closely related to the evolving FEMs for PDEs on evolving surfaces. Optimal-order convergence in the $L^2$ and $H^1$ norms of evolving FEMs for linear and nonlinear PDEs on evolving surfaces were shown in \cite{Dziuk-Elliott-2007,EllioESFEM2015,KLLP-2017}. 


The above-mentioned research efforts have focused on diffusion equations with and without advection terms. The analysis of ALE methods for the Stokes and Navier--Stokes equations has also yielded noteworthy results but remained suboptimal, as discussed below. In \cite{Legendre2008}, Legendre \& Takahashi introduced a novel approach that combines the method of characteristics with finite element approximation to the ALE formulation of the Navier--Stokes equations in two dimensions, and established an $ L^{2} $ error estimate of $O(\tau + h^{1/2})$ for the $P_{1b}$--$P_1$ elements under certain restrictions on the time step size. In a related work \cite{SanMart2009}, an error estimate of $O(h^2|\log h|)$ was obtained for the ALE semidiscrete FEM with the Taylor--Hood $P_2$--$P_1$ elements for the Stokes equations in a time-dependent domain. Moreover, for a fully discrete ALE method with the implicit Euler scheme in time, convergence of $O(\tau + h^2 + h^2/\tau)$ was proved in \cite{SanMart2009}. 
The errors of ALE finite element solutions to the Stokes equations on a time-varying domain, with BDF-$k$ in time (for $1\le k \le 5$) and the Taylor--Hood $P_r$--$P_{r-1}$ elements in space (with degree $r\ge 2$), was shown to be $O(\tau^k + h^r)$ in the $ L^{2} $ norm in \cite{Liu2013}. 

As far as we know, optimal-order convergence of ALE semidiscrete and fully discrete FEMs were not established for the Stokes and Navier--Stokes equations in an evolving domain. As shown in \cite{GawE2015,buyang22}, the optimal-order convergence of ALE semidiscrete FEM requires proving the following optimal-order convergence for the material derivative of the Ritz projection: 
\begin{equation}
    \| D_{t,h}R_{h} u - R_{h}D_{t,h} u \|_{L^{2}} \le C h^{r+1} . 
\end{equation} 
This can be proved using the techniques in \cite{GawE2015,buyang22} by additionally establishing and utilizing the $H^{-1}$ error estimate for the Ritz projection of pressure, i.e., \eqref{estimate-Rhq-hat}, which is used in Lemma \ref{L^2-duality-2} and proved in Appendix \ref{appendix-A}. This leads to optimal-order convergence of the ALE semidiscrete FEM, as a minor result of this article (see Theorem \ref{THM-SD}). 

The main result of this article is the formulation and convergence analysis for a fully discrete second-order projection method along the trajectories of the evolving mesh for decoupling the unknown solutions of velocity and pressure. Although second-order convergence in the discrete $L^2(0,T;L^2)$ norm for the projection method has been proved in \cite{ShenJ1996} in a fixed domain, the convergence of projection methods along trajectories of an evolving mesh for Stokes equations on an evolving domain remains an open question. We fill in this gap by providing comprehensive analysis of the aforementioned error bound, with second-order convergence (up to a logarithmic factor) in the discrete $L^\infty(0,T;L^2)$ norm; see Theorem \ref{main-theorem}. A crucial aspect of the analysis involves proving a discrete $L^1(0,T;L^2)$ estimate for the discretized Stokes equations (in a space-time duality argument), as detailed in Appendix \ref{ell1L2-estimate-section}. 

For simplicity, we focus on the analysis of ALE semidiscrete and fully discrete projection methods for the Stokes equations. However, the numerical scheme and analysis presented in this article can be readily extended to the Navier--Stokes equations. The methodologies employed can be effectively utilized to tackle the nonlinear terms as well. 

The rest of this article is organized as follows. In Section \ref{Preliminary} we present preliminary results for the evolving mesh, ALE finite element spaces, and boundary-skin estimates. In Section \ref{section-semidiscrete}, we present the formula for the semidiscrete finite element approximation and introduce an optimal-order error estimate for it. The second-order fully discrete projection method scheme is then presented in Section \ref{second-fully-section}, where we also prove the optimal-order error estimate of the fully discrete scheme. To reinforce the theoretical analysis, Section \ref{sec:Numerical experiments} includes numerical results for the Stokes equation, Navier--Stokes equation, and Navier--Stokes equation with propeller rotation. These results serve as empirical evidence supporting our theoretical findings.


\section{Preliminary}\label{Preliminary} 

\subsection{Evolving mesh and ALE finite element spaces}
Suppose that the initial smooth domain $\Omega^0$ is divided into a set $\mathcal{T}_h^0$ of shape-regular and quasi-uniform curved simplices with maximal mesh size $h$. Each curved simplex $K$ is associated with a unique polynomial $F_K$ of degree $r$, referred to as the parametrization of $K$ (as described in \cite{ElliottUFEM2021}). This parametrization maps the reference simplex $\hat{K}$ onto the curved simplex $K$.
Additionally, each boundary simplex $K$ (with one face or edge attached to the boundary) may contain a curved face or edge that needs to interpolate the boundary $\Gamma^0 = \partial\Omega^0$. To achieve this interpolation, we employ iso-parametric finite elements of Lenoir's type (see \cite{Lenoir1986} for further details) at time $t=0$ based on the parametrization of the boundary which is denoted by $\varUpsilon:\partial \tilde{D}\to \Gamma^0$. Here, $\partial \tilde{D}$ represents the flat boundary face of the triangulated flat domain, which has the same vertices as the curved triangulated domain $\Omega_h^0=\bigcup_{K\in\mathcal{T}_h^0}K$. In practical implementations, the parametrization $\varUpsilon$ can be chosen as the normal projection onto $\Gamma^0$. In other words, it computes the unique point $\varUpsilon(x)\in\Gamma^0$ satisfying the equation:
\begin{align*}
	x = \varUpsilon(x) + \text{sign}(x,\Omega^{0})|x-\varUpsilon(x)|\textbf{n}(\varUpsilon(x)),
\end{align*}
where $\textbf{n}(\varUpsilon(x))$ is the unit outward normal vector at point $\varUpsilon(x)$ and 
\begin{align*}
	\text{sign}(x,\Omega^{0}) = 
	\left\{\begin{array}{lll}
		1 & \text{for} & x\in \R^d\backslash\overline{\Omega^{0}},\\
		-1 & \text{for} & x\in \Omega^{0}.
	\end{array}\right.
\end{align*}

Let us denote the nodes of the triangulation $\mathcal{T}_h^0$ as $\xi_j\in\R^d$, where $j = 1,\cdots,N$. Each node $\xi_j$ undergoes a time evolution with velocity $w$, resulting in the movement of the node to a point $x_j(t)\in\R^d$ at time $t$. This evolution is governed by an ordinary differential equation (ODE):
\begin{align*}
	\frac{\d}{\d t}x_j(t) = w(x_j(t),t)\quad\text{and}\quad x_j(0) = \xi_j.
\end{align*}
Consequently, the points $x_j(t)$, where $j = 1,\cdots,N$, constitute the nodes of a time-dependent triangulation denoted as $\mathcal{T}_h(t)$. The relations among these points mirror those among the original nodes $\xi_j$, namely, a set of nodes $x_j(t)$ form the vertices of a simplex in $\mathcal{T}_h(t)$ if and only if the corresponding nodes $\xi_j$ form the vertices of a simplex in $\mathcal{T}_h^0$. Hence, the evolving domain $\Omega_h(t) = \bigcup_{K\in\mathcal{T}_h(t)}K$ serves as an approximation of the exact domain $\Omega(t)$. This approximation is achieved by employing piecewise polynomial interpolation of degree $r$ on the reference simplex, with an associated interpolation error of $O(h^{r+1})$. 
Note that the approximation to $\Omega(t)$ by $\Omega_h(t)$ may not be Lenoir's type for $t>0$.

In a manner similar to the initial triangulation $\mathcal{T}_h^0$, each simplex $K\in\mathcal{T}_h(t)$ is associated with a unique polynomial of degree $r$, denoted as $F_K^t:\hat K \to K$, which serves as a parametrization of $K$ over time. Therefore, the finite element space defined on the evolving discrete domain $\Omega_h(t)$ is given by: 
\begin{align*}
	V_h^r(\Omega_h(t)):=\{v_h\in C(\Omega_h(t)): v_h\circ F_K^t\in P^r(\hat K) \,\,\text{for}\,\,\text{all}\,\,K\in\mathcal{T}_h(t)\},
\end{align*}
where $P^r(\hat K)$ represents the set of polynomials on $\hat K$ with degree less than or equal to $r$. The finite element basis functions of $V_h^r(\Omega_h(t))$ are denoted as $\phi_j^t$, where $j = 1,\cdots,N$. These basis functions satisfy the property: 
\begin{align*}
	\phi_j^t(x_i(t)) = \delta_{ij},\,\, i,j = 1,\cdots,N.
\end{align*}
In terms of these basis functions, the approximated flow map $\Phi_h(\cdot,t)\in V_h^k(\Omega_h^0)^d$ can be expressed as
\begin{align*}
	\Phi_h(\xi,t) = \sum_{j = 1}^N x_j(t)\phi_j^0(\xi)\,\,\text{for}\,\, \xi\in\Omega_h^0.
\end{align*}
The flow map $\Phi_h(\cdot,t)$ establishes a one to one correspondence between $\Omega_h^0$ and $\Omega_h(t)$ at time $t$, with a velocity field $w_h\in V_h^r(\Omega_h(t))^d$ satisfying:
\begin{align}\label{w-h-represent}
	w_h(\Phi_h(\xi,t),t) = \frac{\d}{\d t}\Phi_h(\xi,t) = \sum_{j = 1}^N w(x_j(t),t)\phi_j^0(\xi)   \,\,\text{for}\,\,\xi\in\Omega_h^0.
\end{align}
This representation corresponds to the unique Lagrange interpolation of the exact velocity $w(\Phi(\cdot,t),t)$. Analogous to definition \eqref{u-material-derivative}, we can define the material derivative of any finite element function $v_h\in V_h^r(\Omega_h(t))$ with respect to the discrete velocity field $w_h$ as follows: 
\begin{align*}
	D_{t,h} v_h(x,t) := \frac{\d}{\d t}v_h(\Phi_h(\xi,t),t) = \partial_t v_h+ \nabla v_h\cdot w_h\,\,\text{at}\,\, x= \Phi_h(\xi,t)\in\Omega_h(t)\;\text{for}\;\xi \in \Omega_{h}^{0}.
\end{align*}
The pullback of the finite element basis function $\phi_h^t$ from the domain $\Omega_h(t)$ to $\Omega_h(s)$, i.e., $\phi_h^t\circ\Phi_h(\cdot,t)\circ\Phi_h(\cdot,s)^{-1}$, gives rise to a finite element function defined on $ \Omega_{h}(s) $. Remarkably, the nodal values of this function coincide with those of $\phi_h^s$. As a result, we establish the equality $\phi_h^t\circ\Phi_h(\cdot,t)\circ\Phi_h(\cdot,s)^{-1} = \phi_h^s$. Exploiting this relationship, we can derive the well-known transport property of the basis function $\phi_h^t$, which states:  
\begin{align}
	D_{t,h}\phi_h^t(x) = \frac{\d}{\d t}\phi_h^0(\xi) = 0 \,\,\text{at}\,\, x = \Phi_h(\xi,t).
\end{align}

The analysis of integrals over dynamically evolving domains necessitates the application of the Transport Theorem, as established in \cite[Lemma 5.7]{Walker2015}. This pivotal theorem provides a concise and indispensable description of the intrinsic relationship between the time derivative of an integral over a domain that evolves with time and the derivatives of the integrated function and domain velocity. 

\begin{lemma}[Transport Theorem]\label{transport-theorem-lemma}
	If the domain $\Omega$ undergoes motion with a velocity field $w\in W^{1,\infty}(\Omega)$, we have
	\begin{align}
		\frac{\d}{\d t}\int_{\Omega} f \,\d x = \int_\Omega D_t f + f \nabla \cdot w\,\d x,
	\end{align}
	where $D_t f$ is the material derivative of $f$ with respect to the velocity $ w $.
\end{lemma}
The interaction between the operators $ D_{t} $ and $ \nabla $ plays an essential role in the error analysis. Consequently, we establish the following lemma as a direct consequence of \eqref{u-material-derivative}:
\begin{lemma}\label{commute-lemma}
    For any scalar or vector-valued function $ f $, the material derivative of $ \nabla f $ and $ \nabla \cdot f $ with respect to the velocity field $ w $ can be expressed as follows:
    \begin{align}
	    D_{t}\nabla f & = \nabla D_{t}f - \nabla f \nabla w,\\
	    D_{t}\nabla \cdot f & = \nabla \cdot D_{t}f- (\nabla f):(\nabla w)^{\top}.
    \end{align} 
\end{lemma} 
In order to facilitate the implementation of projection methods, we adopt the Taylor--Hood type finite element spaces on the evolving domain $\Omega(t)$, which allow for a continuous approximation of the pressure. Specifically, we define the following spaces:
\begin{align*}
	\mathring{V}^r_h(\Omega_h(t)):= & \{ u\in V_h^r(\Omega_h(t)): u|_{\partial\Omega_h(t)} = 0\}, \\
	Q_h^{r-1}(\Omega_h(t)) := & \{ p\in V_h^{r-1}(\Omega_h(t)): \int_{\Omega_h(t)} p\,\d x = 0\}.
\end{align*}

By employing Verf\"urth's trick and utilizing the macros-element criterion, as described in \cite[Section 8.5 and Section 8.8]{boffi13}, we establish the inf-sup condition for the Taylor--Hood type isoparametric elements.
\begin{lemma}[Inf-sup condtion]
	There exists a constant $\kappa>0$, independent of $h$ and $t\in [0,T]$ for $r\ge 2$, such that 
	\begin{align}\label{inf-sup-condition}
		\sup_{0\neq v_h\in \mathring V^r_h(\Omega_{h}(t))}\frac{({\rm div}v_h, p_h)_{\Omega_h(t)}}{\|\nabla v_h\|_{\Omega_h(t)}}\geq \kappa \|p_h\|_{\Omega_h(t)} \quad \forall p_h\in Q^{r-1}_h(\Omega_{h}(t)), 
	\end{align}
	where $(\cdot,\cdot)_{\Omega_h(t)}$ and $\|\cdot\|_{\Omega_h(t)}$ denote the $L^2$ innerproduct and norm over $\Omega_h(t)$, respectively.
\end{lemma}

\subsection{Mesh velocity and its material derivative}
For any function $u$ defined on $\bigcup_{0\leq t\leq T} \Omega(t)\times \{t\}$, there is an extension function $\widetilde{u}$ defined on $\mathbb{R}^d \times [0,T]$ such that
\begin{align}\label{stein-extension} 
	\widetilde{u}(\cdot, t):=E(u(\cdot, t)\circ \Phi(\cdot,t))\circ \widetilde{\Phi}(\cdot,t)^{-1},
\end{align}
where the operator $E:L^1(\Omega(0))\to L^1(\mathbb{R}^d)$ refers to Stein's extension operator, and $ \widetilde{\Phi} $ is the extended flow map generated by the extended velocity $ E (w (\Phi (\cdot,t),t)) $. It holds that
\begin{align}
	\|\widetilde{u}(\cdot,t)\|_{W^{k,p}(\mathbb{R}^d)}\leq C\|u(\cdot,t)\|_{W^{k,p}(\Omega(t))}.
\end{align}
We denote the interpolation operators as $ I_{h}^{\mathring V}(t):H^1_0(\Omega_h(t))\to \mathring{V}^r_h(\Omega_{h}(t)) $ and $I_h^Q(t):H^1(\Omega_h(t))$ $\to Q^{r-1}_h(\Omega_{h}(t))$. Throughout this discussion, the explicit time dependency $ t $ is often omitted, and we will use $ I_{h}^{\mathring V} $ and $ I_{h}^{Q} $ instead.

In certain cases, we come across vector-valued spaces such as $\mathring V_h^r(\Omega_h(t))^d$ and vector-valued interpolation operators such as $I_h^{\mathring V^d}$. To streamline the notation, we will use $\mathring V_h^r(\Omega_h(t))$ and $I_h^{\mathring V}$ when referring to vector-valued objects, provided there is no ambiguity within the context.

By \eqref{w-h-represent}, the interpolation $w_h$ serves as an approximation of $w$. Consequently, we can establish an error estimate in the $W^{1,\infty}$ norm for $w_h$ as follows: 
\begin{align}\label{w1infty-estimate-wh}
	\|w_h(\cdot,t) - w(\cdot,t)\|_{W^{1,\infty}(\Omega_h(t))} \le C h^r\|w(\cdot,t)\|_{W^{r+1,\infty}(\Omega_h(t))}.
\end{align}
Similarly, we can obtain the $ W^{1,\infty} $ error estimates for $ D_{t,h}w_{h} $ and $ D_{t,h}^{2}w_{h} $. By employing Gronwall's inequality and exploiting the fact that $ w $ corresponds to the velocity of the flow map, we can derive the error between $ \Phi (t) $ and $ \Phi_{h} (t) $ as follows: 
\begin{align}\label{flow-map-err}
    \|\Phi (\cdot,t) - \Phi_{h} (\cdot,t)\|_{W^{1,\infty}(\Omega_{h}^{0})} \le C h^{r}.
\end{align}
The estimates in \eqref{w1infty-estimate-wh} and \eqref{flow-map-err}, and the similar estimates for $ D_{t,h}w_{h} $ and $ D_{t,h}^{2}w_{h} $ as discussed above, lead to the following result (via using the inverse inequality of finite element functions):
\begin{subequations}\label{assumption-w-phi}
\begin{align}
	\| D^{2}_{t,h}w_{h} \|_{W^{1,\infty}(\Omega_{h}(t))}+\| \Phi_h (\cdot,t) \|_{W^{1,\infty}(\Omega_{h}^{0})}+\| \Phi_{h}^{-1}(\cdot,t) \|_{W^{1,\infty}(\Omega_{h}(t))} & \le C, \\
	\| w_{h}(t) - \widetilde{w} (t) \|_{W^{1,\infty}(\Omega_{h}(t))} + \| D_{t,h} (w_{h}(t) - \widetilde{w}(t)) \|_{W^{1,\infty}(\Omega_{h}(t))} & \le C h, 
\end{align} 
\end{subequations}
where $ C $ is a constant independent of the mesh size $ h $ and time $ t $. This serves as a basic condition on the mesh velocity in the subsequent analysis.

\subsection{Discrepancy between $\Omega(t)$ and $\Omega_{h}(t)$}

To address the discrepancy between $\Omega(t)$ and its finite element approximation, $\Omega_{h}(t)$, we utilize the boundary-skin estimate. This estimate is essential for effectively managing errors that arise from the finite element approximation of the domain. The rigorous proof of this estimate can be found in \cite[Section 3.4]{buyang22}.
\begin{lemma}\label{boundary-skin}
	For any finite element function $v_h\in \mathring V_h^r(\Omega_h(t))$, the following inequalities hold:
	$$
	\|v_h\|_{L^2(\Omega_h(t)\setminus \Omega(t))}\leq Ch^{3(r+1)/2-d/2}\|\nabla v_h\|_{L^2(\Omega_{h}(t))}\leq Ch^{3(r+1)/2-d/2-1}\|v_h\|_{L^2(\Omega_h(t))}.
	$$
\end{lemma}

Due to the inherent discrepancy between the finite element domain $ \Omega_{h}(t) $ and the exact domain $ \Omega (t) $, the exact solution $ u $ does not vanish on $ \partial \Omega_{h}(t) $. To handle this situation, we rely on the following lemma to derive an estimate for the integral over the boundary $ \partial\Omega_h(t)$. A proof of this lemma can be found in \cite[(3.32)]{buyang22}.

\begin{lemma}\label{boundary-integral}
	Let $g\in W^{1,1}(\mathbb{R}^d)$. Then the following inequality holds: 
	\begin{equation}
		\|g\|_{L^1(\partial\Omega_h(t))}\leq C\|g\|_{L^1(\partial\Omega(t))}+C\|\nabla g\|_{L^1(\mathbb{R}^d)},
	\end{equation}
	where $ C $ is a constant independent of the mesh size $ h $ and time $ t $.
\end{lemma}

The significance of the ensuing lemma lies in its pivotal role in acquiring optimal $H^{-1}$-norm estimates for pressure through implementation of a duality argument. A rigorous proof of this lemma can be found in \cite[Corollary 1.5]{farwig94}.
\begin{lemma}\label{div-regular-inverse}
	Let $\lambda\in H^1(\Omega(t))\cap L_{0}^{2}(\Omega (t))$. Then, there exists a function $\chi\in H^2(\Omega(t))\cap H^1_0(\Omega(t))$ such that ${\rm div}\chi=\lambda$, and the following inequality holds:
	$$
	\|\chi\|_{H^2(\Omega(t))}\leq C\|\lambda\|_{H^1(\Omega(t))},
	$$
	where the constant $C$ is independent of $t\in [0,T]$.
\end{lemma}




\section{The semidiscrete finite element approximation}\label{section-semidiscrete}
In this section, we present a semidiscrete FEM for solving problem \eqref{main-equation-stokes} using the ALE formulation \eqref{u-material-derivative}. We also provide an error estimate for this method. Let $u_h\in \mathring V^r_h(\Omega_h(t))$ be a finite element function. It is worth noting that $u_h\circ \Phi_h(\cdot,t)\circ \Phi_h^{-1}(\cdot,s)\in \mathring V^r_h(\Omega_h(s))$ and shares the same nodal values as $u_h$. This observation allows us to define an equivalence relation $\sim$ on the disjoint union $\bigcup_{0\leq t\leq T} \mathring V^r_h(\Omega_h(t))$ as follows:
$$
u_h \sim v_h \text{ iff } v_h\circ \Phi_h(\cdot,s)\circ \Phi_h^{-1}(\cdot,t) = u_h,\,\, \forall \,u_h\in \mathring V^r_h(\Omega_h(t)), v_h\in \mathring V^r_h(\Omega_h(s)).
$$
Let $ \{u_{h}\} $ denote the equivalence class generated by $u_h$ in $\bigcup_{0\leq t\leq T} \mathring V^r_h(\Omega_h(t))$, and let $[\mathring V_h^r]$ represent the set of equivalence classes. We can define the following calculus on $[\mathring V_h^r]$:
\begin{enumerate}
	\item Given $\{u_h\}\in [\mathring V_h^r]$, we define $\{u_h\}[s]$ as the unique function in $ \mathring V^r_h(\Omega_h(s))$ that belongs to the equivalence class $ \{u_{h}\} $. In other words, $ \{u_{h}\}[s] $ serves as the representative of the equivalence class $\{u_h\}$ at time $s$.
	
	\item Consider $\{u_h\}, \{v_h\}\in [\mathring V_h^r]$. For any $0\leq s,t\leq T$, we observe that $\{u_h\}[s]+\{v_h\}[s]\sim \{u_h\}[t]+\{v_h\}[t]$. This equivalence relation allows us to define addition between equivalence classes as follows: 
	$$
	\{u_h\}+\{v_h\}:= \left\{\{u_h\}[t]+\{v_h\}[t]\right\}\in [\mathring V_h^r].
	$$
	This definition implies that $[\mathring V_h^r]$ forms a finite dimensional vector space.
	\item We can establish the integral calculus on equivalence classes as follows:
	\begin{align*}
	\big(\{u_h\},\{v_h\}\big)_{\Omega_h(t)} & :=\big(\{u_h\}[t],\{v_h\}[t]\big)_{\Omega_h(t)},\\
	\big(\nabla \{u_h\}, \nabla \{v_h\}\big)_{\Omega_h(t)} & :=\big(\nabla (\{u_h\}[t]), \nabla (\{v_h\}[t])\big)_{\Omega_h(t)}.
	\end{align*}
	Moreover, we define the $L^2$ norm of $\{u_h\}$ on $\Omega_h(t)$ as 
	\begin{align*}
		\|\{u_h\}\|_{\Omega_h(t)} := \sqrt{(\{u_h\},\{u_h\})_{\Omega_h(t)}}.
	\end{align*}
\end{enumerate} 
These established definitions allow us to treat equivalence classes on an equal footing with individual functions, enabling us to define inner products and norms for these classes. 

Similarly, we extend the concept of equivalence classes to $[Q_h^{r-1}]$, which represents the set of equivalence classes in $\bigcup_{0\leq t\leq T}Q^{r-1}_h(\Omega_h(t))$. We define the following pairing between $[\mathring V_h^r]$ and $[Q^{r-1}_h]$:
$$
\big({\rm div} \{u_h\}, \{q_h\}\big)_{\Omega_h(t)}:=\big({\rm div}(\{u_h\}[t]), \{q_h\}[t]\big)_{\Omega_h(t)}.
$$
In the subsequent analysis, for the sake of simplicity, we adopt a simplified notation where $u_h$ directly represents the equivalence class $\{u_h\}$. Similarly, we use $\mathring V_{h}$ to denote the set $[\mathring V_{h}^{r}]$ and $Q_{h}$ to represent $[Q_{h}^{r-1}]$, provided there is no ambiguity. 

The semidiscrete finite element problem can be formulated as follows: Seek solutions $u_h\in \mathring V_h$ and $p_h\in Q_h$ that satisfy the following equations for all test functions $v_h \in \mathring V_h$ and $q_h \in Q_h$:
\begin{subequations}\label{main-semidiscrete}
	\begin{align}
	(D_{t,h}u_h-w_h\cdot \nabla u_h, v_h)_{\Omega_h(t)}+(\nabla u_h, \nabla v_h)_{\Omega_h(t)}-({\rm div} v_h, p_h)_{\Omega_h(t)}=&\,(f,v_h)_{\Omega_h(t)} ,\label{main-semidiscrete-a}\\
	({\rm div} u_h, q_h)_{\Omega_h(t)}=&\, 0,\label{main-semidiscrete-b}\\
	u_h(0)=&\,u_{h,0} \in \mathring V_{h}.\label{main-semidiscrete-c}
	\end{align}
\end{subequations}

The main result of this section is the following theorem.

\begin{theorem}[Error estimate of the semidiscrete FEM]\label{THM-SD}
	Consider the semidiscrete finite element solution $u_h$ given by \eqref{main-semidiscrete}, and denote the error by $e:=\widetilde{u}-u_h$. Assuming that the exact solution to problem \eqref{main-equation-stokes} is sufficiently smooth, the following estimate holds under condition \eqref{assumption-w-phi}{\rm:}
\begin{align}
	\sup_{t\in [0,T]}\|e(t)\|_{\Omega_h(t)}\leq Ch^{r+1} + C \| e (0) \|_{\Omega_{h}(0)}.
\end{align}
where $ C $ is a constant independent of the mesh size $ h $ and time $ t $.
\end{theorem}

\subsection{The Stokes--Ritz projection}

Analogous to the Stokes--Ritz projection in a fixed domain, we introduce the concept of the Stokes--Ritz projection for a pair $(v (\cdot,t),q (\cdot,t))\in H^1(\Omega_h(t))\times L^2(\Omega_h(t))$ for $ t \in [0,T] $ over a time-dependent domain, denoted as $(R_hv, R_hq)\in \mathring V_{h}\times Q_h$. The Stokes--Ritz projection satisfies the following equations for all test functions $\chi_h\in \mathring V_h$ and $\lambda_h\in Q_h$: 
\begin{subequations}\label{stokes-ritz}
	\begin{align}
		(\nabla R_hv, \nabla \chi_h)_{\Omega_h(t)}- ({\rm div}\chi_h, R_hq)_{\Omega_h(t)}=&\,(\nabla v, \nabla \chi_h)_{\Omega_h(t)}- ({\rm div}\chi_h, q)_{\Omega_h(t)},\\
		({\rm div}R_h v, \lambda_h)_{\Omega_h(t)}=&\,({\rm div} v, \lambda_h)_{\Omega_h(t)} ,
	\end{align}
\end{subequations}
Additionally, we define the norm $\|\cdot\|^{'}$ over any domain $D\subset\R^d$ as follows:
\begin{align}
	\|f\|^{'}_{D} := \|f - \bar f\|_{L^2(D)},
\end{align}
where $\bar f$ denotes the average of $f$ over $ D $, given by $\bar f :=\frac{1}{|D|}\int_D f\,\d x$.

To assess the errors associated with the Stokes--Ritz projection, it is essential to rely on the following stability estimate of the discrete Stokes operator. This estimate can be rigorously proven by leveraging the inf-sup condition \eqref{inf-sup-condition}.
\begin{lemma}\cite[Chapter 1, Theorem 4.1]{Girault2012}\label{stokes-ritz-estimate-0}
Let $(u_h,p_h)\in \mathring V_h\times Q_h$ satisfy
\begin{subequations}
	\begin{align*}
		(\nabla u_h,\nabla v_h)_{\Omega_h(t)}-({\rm div}v_h, p_h)_{\Omega_h(t)}= &\,\ell^{t}(v_h) \quad\forall v_h\in \mathring V_h,\\
		({\rm div}u_h, q_h)_{\Omega_h(t)}=&\,\phi^{t} (q_h) \quad \forall q_h\in Q_h,
	\end{align*}
\end{subequations}
where $\ell^{t}$ is a linear functional on $\mathring V_h$ and $\phi^{t}$ is a linear functional on $Q_h$, both defined for $ t \in [0,T] $. We may use $ \ell $ and $ \phi $ without the superscript $ t $ when there is no ambiguity. Then the following estimate holds:
$$
\|\nabla u_h\|_{\Omega_h(t)}+\|p_h\|_{\Omega_h(t)}\leq C\left(\|\ell\|_{-1,h}+\|\phi\|_h \right),
$$
where
$$
\|\ell\|_{-1,h}:=\sup_{0\neq v_h\in \mathring V_h}\frac{|\ell(v_h)|}{\|\nabla v_h\|_{\Omega_h(t)}}, \quad \|\phi\|_h:=\sup_{0\neq q_h\in Q_h}\frac{|\phi(q_h)|}{\|q_h\|_{\Omega_{h}(t)}},
$$
and $C$ is a constant independent of the mesh size $h$ and time $t$.
\end{lemma}
As a direct consequence of Lemma \ref{stokes-ritz-estimate-0}, the Stokes--Ritz projection exhibits quasi-optimal error estimates, as stated in the following lemma:
\begin{lemma}\cite[Chapter 2, Theorem 1.1]{Girault2012}\label{H^1-error-1}
	Let $(R_h v, R_h q)$ be the Stokes--Ritz projection of $(v,q)\in H^1(\Omega_h(t))\times L^2(\Omega_h(t))$ as defined in \eqref{stokes-ritz}. Then the following estimate holds:
	$$
	\|\nabla (v-R_h v)\|_{\Omega_h(t)}+\|q-R_hq\|^{'}_{\Omega_h(t)}\leq C(\inf_{v_h\in \mathring V_h} \|\nabla(v-v_h)\|_{\Omega_h(t)}+\inf_{q_h\in Q_h} \|q-q_h\|_{\Omega_h(t)}).
	$$
	In particular, when $v(t)\in H^{r+1}_h(\Omega_h(t))$ and $p(t)\in H^r_h(\Omega_h(t))$ for all $ t \in [0,T] $, the following estimate holds: 
	$$
	\|\nabla (v-R_h v)\|_{\Omega_h(t)}+\|q-R_hq\|^{'}_{\Omega_h(t)}\leq Ch^r\left(\|v\|_{H^{r+1}_h(\Omega_h(t))}+\|q\|_{H^r_h(\Omega_h(t))}\right).
	$$
\end{lemma}

In order to prove the optimal-order estimate, we need to facilitate the estimation of errors such as $\|D_{t,h}v-D_{t,h}R_h v\|_{\Omega_h(t)}$ and $\|D_{t,h}q-D_{t,h}R_h q\|^{'}_{\Omega_h(t)}$. It is convenient to introduce the operator $E_{t,h}$ defined as:
\begin{align}
    E_{t,h}:= D_{t,h}R_{h} - R_{h}D_{t,h}.
\end{align}

For the $ H^{1} $-estimate of $ E_{t,h} $, by taking discrete material derivative of the equation \eqref{stokes-ritz} and applying Lemma \ref{transport-theorem-lemma}, Lemma \ref{commute-lemma} and Lemma \ref{stokes-ritz-estimate-0} under the assumption \eqref{assumption-w-phi}, we can establish the following lemma using techniques similar to those employed in \cite[Lemma 3.5]{buyang22}.
\begin{lemma}\label{H^1-error-2}
Let $(R_h v, R_h q)$ and $(R_hD_{t,h}v,R_hD_{t,h}q)$ denote the Stokes--Ritz projections of $(v,q)$ and $(D_{t,h}v,D_{t,h}q)$, respectively, where $ (v,q) \in H^{1}(\Omega_{h}(t))\times L^{2}(\Omega_{h}(t)) $. Then the following estimate holds: 
\begin{align}
	\|\nabla E_{t,h}v\|_{\Omega_h(t)}+\|E_{t,h}q\|^{'}_{\Omega_h(t)} \leq \,C\left(\|\nabla(R_hv-v)\|_{\Omega_h(t)}+\|R_hq-q\|^{'}_{\Omega_h(t)}\right).
\end{align}
\end{lemma}
For the $ H^{1} $-estimate of $ D_{t,h}E_{t,h} $, by taking the second order discrete material derivative of the equation in \eqref{stokes-ritz} and employing  similar methods as those used in the proof of Lemma \ref{H^1-error-2}, we can derive the following lemma.
\begin{lemma}\label{H^1-error-3}
	Let $ (v,q) \in H^{1}(\Omega_{h}(t))\times L^{2}(\Omega_{h}(t))  $ and $ (D_{t,h}v,D_{t,h}q)\in H^{1}(\Omega_{h}(t))\times L^{2}(\Omega_{h}(t)) $. Then there exists a constant $ C > 0 $ independent of $ h $ and $ t $ such that
	\begin{align}
		\|\nabla D_{t,h}E_{t,h}v\|_{\Omega_h(t)}+\|D_{t,h}E_{t,h}q\|^{'}_{\Omega_h(t)} \leq &\, C\|\nabla(R_hv-v)\|_{\Omega_h(t)}+C\|R_hq-q\|^{'}_{\Omega_h(t)}\nonumber\\
		+ C\|\nabla (R_hD_{t,h}v- & D_{t,h}v)\|_{\Omega_h(t)}+C\|D_{t,h}q-R_hD_{t,h}q\|^{'}_{\Omega_h(t)}.
	\end{align}	
\end{lemma}
\subsection{The Nitche's trick and duality argument}
In order to obtain an optimal order error estimate of $\|R_h v-v\|_{\Omega_h(t)}$, we will apply Nitche's trick. Let $g_h$ be a function in $ \mathring V_h$ that we can extend outside of $\Omega_h(t)$ by setting it to zero. We solve the following equations in $\Omega(t)$:
\begin{subequations}\label{Nitche_trick}
	\begin{align}
		-\Delta \psi +\nabla \varphi = &\, g_h \quad \text{in }\Omega(t), \\
		{\rm div} \psi  = &\, 0 \quad\:\: \text{in }\Omega(t),\quad \psi|_{\partial\Omega(t)}=0.
	\end{align}
\end{subequations}
By applying regularity estimates for the Stokes equation in $\Omega(t)$, we obtain the following result:
\begin{equation}\label{regularity-stokes}
	\|\psi\|_{H^2(\Omega(t))}+\|\nabla \varphi\|_{L^2(\Omega(t))}\leq C\|g_h\|_{\Omega_h(t)}.
\end{equation} 
To extend the functions $\psi$ and $\varphi$ to $\widetilde{\psi}$ and $\widetilde{\varphi}$, respectively, we employ the Stein extension operator as in \eqref{stein-extension}. By applying this operator, we can define $\widetilde{\eta}$ as $\widetilde{\eta}:= -\Delta \widetilde{\psi}+\nabla\widetilde{\varphi}-g_h$ and arrive at the following expression:
\begin{align}\label{nitche-1}
	\|g_h\|^2_{\Omega_h(t)}
	=(\nabla\widetilde{\psi}, \nabla g_h )_{\Omega_h(t)}-({\rm div}g_h, \widetilde{\varphi})_{\Omega_h(t)}-(g_h, \widetilde{\eta})_{\Omega_h(t)}.
\end{align} 
Notably, since $\widetilde{\eta}$ vanishes in $\Omega(t)$ and we have $r \geq 2$, we can utilize Lemma \ref{boundary-skin} along with the regularity estimate \eqref{regularity-stokes} to obtain the following inequality:
\begin{align*}
	|(g_h, \widetilde{\eta})_{\Omega_h(t)}|&=|(g_h, \widetilde{\eta})_{\Omega_h(t)\setminus\Omega(t)}|\leq \|\widetilde{\eta}\|_{L^2(\mathbb{R}^d)}\|g_h\|_{L^2(\Omega_h(t)\setminus\Omega(t))}\leq Ch^{2}\|g_h\|^2_{\Omega_h(t)}.
\end{align*}
Consequently, when $h>0$ is sufficiently small, we can absorb $|(g_h, \widetilde{\eta})_{\Omega_h(t)}|$ on the right-hand side of \eqref{nitche-1} by the left-hand side. This yields the following estimate:
\begin{align}\label{estimate-g-h}
	\|g_h\|^2_{\Omega_h(t)}
	\leq& C\left|(\nabla(\widetilde\psi_h-\widetilde{\psi}), \nabla g_h )_{\Omega_h(t)}\right|+C\left|({\rm div}g_h, \widetilde{\varphi}-\widetilde\varphi_h)_{\Omega_h(t)}\right| \notag\\
	    & + C\left|(\nabla\widetilde\psi_h, \nabla g_h )_{\Omega_h(t)}-({\rm div}g_h, \widetilde\varphi_h)_{\Omega_h(t)}\right|\nonumber\\
	\leq & Ch\|g_h\|_{\Omega_h(t)}\|\nabla g_h\|_{\Omega_h(t)}+C\left|(\nabla\widetilde\psi_h, \nabla g_h )_{\Omega_h(t)}-({\rm div}g_h, \widetilde\varphi_h)_{\Omega_h(t)}\right|,
\end{align}
where $\widetilde\psi_h:=I_h^{\mathring V}\widetilde{\psi}$, $\widetilde\varphi_h:=I_h^Q\widetilde{\varphi}$, and the last inequality follows from the interpolation error estimates.

Again, let $(R_{h}v, R_{h}q)$ denote the Stokes--Ritz projection of $(v,q)$. We choose $g_h=R_h v-\chi_h$ for $\chi_h\in \mathring V_h$ in \eqref{Nitche_trick}. By appropriately estimating the last terms on the right-hand side of \eqref{estimate-g-h} and applying Young's inequality, Lemma \ref{L^2-duality-1} (see Appendix \ref{appendix-A}) establishes an estimate for $ \| R_{h}v - \chi_{h} \|_{\Omega_{h}(t)} $ as follows:
\begin{align}\label{L^2-duality-1-inequality}
	\|R_{h}v - \chi_{h}\|_{\Omega_h(t)}\leq &Ch\|\nabla (R_{h}v - \chi_{h})\|_{\Omega_h(t)}+Ch\left(\|R_hq-q\|^{'}_{\Omega_h(t)}+\|\nabla (v-\chi_h)\|_{\Omega_h(t)}\right)\nonumber\\
	&+C\|v-\chi_h\|_{\Omega_h(t)}+C\|v\|_{L^\infty(\partial\Omega_h(t))}.
\end{align}

To obtain the optimal order estimate of $E_{t,h}v$, we rely on the negative norm estimate of $R_hq-q$. Lemma \ref{negative-norm-pressure-error} (see Appendix \ref{appendix-A}) establishes the following estimate:
\begin{equation}\label{estimate-Rhq-hat}
	\frac{|(\widehat{R_hq},\lambda)_{\Omega_h(t)}|}{\|\lambda\|_{H^1(\mathbb{R}^d)}}\leq C\left[h\left(\|\widehat{R_hq}\|_{\Omega_h(t)}+\|\nabla(R_hv-v)\|_{\Omega_h(t)}\right)+\|v-R_hv\|_{\Omega_h(t)}+\|v\|_{L^\infty(\partial\Omega_h(t))}\right].
\end{equation}

Having completed the necessary preparations, we are now poised to establish the $ L^{2} $-estimate of $ E_{t,h} $. To achieve this, we adopt a proof technique akin to that used in Lemma \ref{L^2-duality-1}, leveraging the insights gained from \eqref{estimate-g-h} and \eqref{estimate-Rhq-hat}. 
\begin{lemma}\label{L^2-duality-2}
	For $ (v,q) \in H^{1}(\Omega_{h}(t))\times L^{2}(\Omega_{h}(t)) $ and $ (D_{t,h}v,D_{t,h}q) \in H^{1}(\Omega_{h}(t))\times L^{2}(\Omega_{h}(t)) $, if we choose $g_h= E_{t,h}v $ in \eqref{Nitche_trick}, then the following estimate holds:
\begin{align}
	\|E_{t,h}v\|_{\Omega_h(t)}\leq &Ch\|\nabla E_{t,h}v\|_{\Omega_h(t)}+Ch\left(\|R_hq-q\|^{'}_{\Omega_h(t)}+\|\nabla (R_hv-v)\|_{\Omega_h(t)}+\|E_{t,h}q\|^{'}_{\Omega_h(t)}\right)\nonumber\\
	&+C\|R_hv-v\|_{\Omega_h(t)}+C\|v\|_{L^\infty(\partial\Omega_h(t))}.
\end{align}
\end{lemma}

\subsection{Error estimate for the semidiscrete FEM problem}\label{semi-discrette-3}
Consider the exact solutions $(u,p)$ of equation \eqref{main-equation-stokes}, and let $(\widetilde{u},\widetilde{p})$ be the extension of $(u,p)$ using the method described in \eqref{stein-extension}. We define the auxiliary function $\widetilde{\eta}$ as follows:
\begin{equation}
	\widetilde{\eta}:= \partial_t\widetilde{u}-\Delta\widetilde{u}+\nabla \widetilde{p}-f.
\end{equation}
By integrating the above equation against $v_h\in \mathring V_h$ over $\Omega_h(t)$, we obtain:
\begin{align*}
		(D_{t,h}\widetilde{u}-w_h\cdot \nabla \widetilde{u}, v_h)_{\Omega_h(t)}+(\nabla \widetilde{u}, \nabla v_h)_{\Omega_h(t)}-({\rm div}v_h, \widetilde{p})_{\Omega_h(t)}=(f,v_h)_{\Omega_h(t)}+\langle\mathscr{G}(t), v_h\rangle, 
\end{align*}
where the remainder $$\langle\mathscr{G}(t), v_h\rangle:= (\widetilde{\eta}, v_h)_{\Omega_h(t)}$$ is caused by the perturbation of the domain from $\Omega(t)$ to $\Omega_h(t)$. 
Applying H\"older's inequality, Lemma \ref{boundary-skin} and the fact $r\geq 2$, we can derive the following estimate:
\begin{align}\label{boundary-caused-error-1}
	 |(\widetilde{\eta}, v_h)_{\Omega_h(t)}|=|(\widetilde{\eta}, v_h)_{\Omega_h(t)\setminus \Omega(t)}|
	 \leq  Ch^{r+1}\|\nabla v_h\|_{\Omega_h(t)}\|\widetilde{\eta}\|_{L^\infty(\mathbb{R}^d)},
\end{align}
Let us define $(\hat{u}_h, \hat{p}_h):=(R_h\widetilde{u}, R_h\widetilde{p})$. It follows that $(\hat{u}_h, \hat{p}_h)$ satisfies the following equation:
\begin{subequations}\label{equation-u-hat}
	\begin{align}
		(D_{t,h}\hat{u}_h - w_h\cdot \nabla \hat{u}_h, v_h)_{\Omega_h(t)} + &\, (\nabla \hat{u}_h, \nabla v_h)_{\Omega_h(t)}-({\rm div} v_h, \hat{p}_h)_{\Omega_h(t)}	\nonumber\\
		=&\,(f,v_h)_{\Omega_h(t)}+\langle \mathscr{G}(t), v_h\rangle+\langle \mathscr{F}(t), v_h\rangle && \forall v_h\in \mathring V_h, \label{equation-u-hat-a}\\
		({\rm div} \hat{u}_h, q_h)_{\Omega_h(t)}=&\, 0 && \forall q_h\in Q_h,\label{equation-u-hat-b}
	\end{align}
\end{subequations}
where the remainder 
\begin{align}
\langle \mathscr{F}(t), v_h\rangle:=(D_{t,h}(\hat{u}_h-\widetilde{u}), v_h)_{\Omega_h(t)}-(w_h\cdot\nabla(\hat{u}_h-\widetilde{u}), v_h)_{\Omega_h(t)} 
\end{align}
represents the consistency error of the spatial discretization. Via integration by parts, we can estimate the second term on the right-hand side above as follows: 
\begin{align}
	|(w_h\cdot\nabla(\hat{u}_h-\widetilde{u}), v_h)_{\Omega_h(t)}|&=|-(v_h{\rm div}w_h, \hat{u}_h-\widetilde{u})_{\Omega_h(t)}-(\hat{u}_h-\widetilde{u}, w_h\cdot\nabla v_h)_{\Omega_h(t)}|\nonumber\\
	&\leq C\|\hat{u}_h-\widetilde{u}\|_{\Omega_h(t)}\|\nabla v_h\|_{\Omega_h(t)}, 
\end{align}
where we have used the $W^{1,\infty}$ boundedness of the mesh velocity, i.e., $\|w_h(t)\|_{W^{1,\infty}(\Omega_h(t))}\leq C$, which follows from \eqref{w1infty-estimate-wh} and the triangle inequality. Thus, we have
\begin{align}\label{FEM-consistency-error-1}
	|\langle \mathscr{F}(t), v_h\rangle|\leq C\|D_{t,h}(\hat{u}_h-\widetilde{u})\|_{\Omega_h(t)}\|v_h\|_{\Omega_h(t)}+ C\|\hat{u}_h-\widetilde{u}\|_{\Omega_h(t)}\|\nabla v_h\|_{\Omega_h(t)}.
\end{align}
The right-hand side of \eqref{FEM-consistency-error-1} is estimated below. 

Firstly, based on the results and Lemmas presented in the preceding subsection, the following $ H^{1} $ error estimate for the Stokes--Ritz projection follows in a standard way.
\begin{lemma}\label{error-Ritz-H1-optimal}
	Let $ (\hat{u}_{h} , \hat{p}_{h} )$ denote the Stokes--Ritz projection $ (R_{h}\tilde{u} ,R_{h} \tilde{p}) $ of $ (\tilde{u}, \tilde{p}) $. Then the following estimates holds: 
\begin{align}\label{inequality-uh-ph-dtuh}
	\|\nabla(\hat{u}_h-\widetilde{u})\|_{\Omega_h(t)}+\|\widetilde{p}-\hat{p}_h\|_{\Omega_h(t)}+\|\nabla D_{t,h}(\hat{u}_h-\widetilde{u})\|_{\Omega_h(t)}+\|D_{t,h}(\hat{p}_h-\widetilde{p})\|^{'}_{\Omega_h(t)}\leq&\, Ch^{r},
\end{align}
where $ C $ is a constant independent of the mesh size $ h $ and time $ t $.
\end{lemma}
\begin{proof}
	By taking $(v,q)=(\widetilde{u}, \widetilde{p})$ in Lemma \ref{H^1-error-1}, and noting that $I_h^{\mathring V}\widetilde{u}|{\partial\Omega_h(t)}=0$, we obtain the following inequality:
\begin{align}\label{lemma-3.10-temp0}
	\|\nabla(\hat{u}_h-\widetilde{u})\|_{\Omega_h(t)}+\|\widetilde{p}-\hat{p}_h-\overline{\widetilde{p}}\|_{\Omega_h(t)}\leq Ch^{r}\left(\|\widetilde{u}\|_{H^{r+1}(\mathbb{R}^d)}+\|\widetilde{p}\|_{H^r(\mathbb{R}^d)}\right).	
\end{align} 
To estimate $\overline{\widetilde{p}}$, we have
\begin{align*}
	|\overline{\widetilde{p}}|=\frac{1}{|\Omega_h(t)|}\Big|\int_{\Omega_h(t)} \widetilde{p}\,\d x\Big|&\leq  \frac{1}{|\Omega_h(t)|} \Big(\Big|\int_{\Omega(t)\setminus\Omega_h(t)}\widetilde{p}\,\d x\Big|+\Big|\int_{\Omega_h(t)\setminus\Omega(t)}\widetilde{p}\,\d x\Big| \Big)\leq Ch^{r+1}\|\widetilde{p}\|_{L^\infty(\mathbb{R}^d)},
\end{align*}
which allows us to estimate the first two terms on the left-hand side of \eqref{inequality-uh-ph-dtuh} using the triangle inequality. The estimates about the last two terms on the left-hand side of \eqref{inequality-uh-ph-dtuh} follow directly from Lemma \ref{H^1-error-1} and Lemma \ref{H^1-error-2}.
\end{proof}

Secondly, to establish the optimal-order convergence, it is imperative to derive the $L^2$ error estimate between the Stokes--Ritz projection and the exact solutions. Notably, we can achieve higher-order convergence for both the difference $\hat{u}_h-\widetilde{u}$ and its material derivative $D_{t,h}(\hat{u}_h-\widetilde{u})$, as demonstrated in the subsequent lemma.
\begin{lemma}\label{error-stokes-ritz}
Let $ (\hat{u}_{h} , \hat{p}_{h} )$ be the Stokes--Ritz projection $ (R_{h}\tilde{u} ,R_{h} \tilde{p}) $ of $ (\tilde{u}, \tilde{p}) $ as in Lemma \ref{error-Ritz-H1-optimal}. Then the following estimate holds: 
	\begin{align}
		\|\hat{u}_h-\widetilde{u}\|_{\Omega_h(t)}+\|D_{t,h}(\hat{u}_h-\widetilde{u})\|_{\Omega_h(t)}\leq &\,Ch^{r+1},
	\end{align}
	where $ C $ is a constant independent of the mesh size $ h $ and time $ t $.
\end{lemma}
\begin{proof}
By choosing $(v,q)=(\widetilde{u}, \widetilde{p})$, $\chi_h=I_h^{\mathring V}\widetilde{u}$ in \eqref{L^2-duality-1-inequality}, and considering the boundary condition $ \tilde{u}|_{\partial \Omega (t)}  = 0 $, we have
\begin{equation}\label{estimate-u-linfty-partial}
	\|\tilde{u}\|_{L^{\infty}( \partial\Omega_{h}(t) )}\le \sup_{x \in \partial \Omega_{h}(t)} \text{dist}( x,\partial \Omega_{h}(t) )\|\tilde{u}\|_{W^{1,\infty}( \R^{d} )} \le C h^{r+1}\|\tilde{u}\|_{W^{1,\infty}(\R^d)}.
\end{equation}
Using this result together with Lemma \ref{error-Ritz-H1-optimal}, we obtain the desired estimate for $ \hat{u}_{h} - \widetilde{u} $.

To estimate $ D_{t,h} (\hat{u}_{h} - \tilde{u}) $, we split it into two terms as follows: 
\begin{align}\label{middle-splitting-identity}
D_{t,h}(\hat{u}_h-\widetilde{u})=\left(D_{t,h}R_h \widetilde{u}-R_hD_{t,h}\widetilde{u}\right)+\left(R_hD_{t,h}\widetilde{u}-D_{t,h}\widetilde{u}\right).
\end{align}
By taking $(v,q)=(\widetilde{u},\widetilde{p})$ in Lemma \ref{L^2-duality-2} and applying Lemma \ref{error-Ritz-H1-optimal} and the estimate for $ \hat{u}_{h} - \widetilde{u} $, we obtain the estimate for the first term on the right-hand side of \eqref{middle-splitting-identity}.

Next, we take $(v,q)=(D_{t,h}\widetilde{u}, D_{t,h}\widetilde{p})$ in Lemma \ref{H^1-error-1} and \eqref{L^2-duality-1-inequality}. Noting that $I_h^{\mathring V}(D_{t,h}\widetilde{u})$ vanishes on the boundary of $ \Omega_{h}(t) $, we can use the same technique as employed in the estimation of $\|R_h\widetilde{u}-\widetilde{u}\|_{\Omega_h(t)}$ to obtain the estimate for the second term on the right-hand side of \eqref{middle-splitting-identity}. This completes the proof.
\end{proof}

By applying Lemma \ref{error-stokes-ritz} and \eqref{boundary-caused-error-1}, we obtain the estimate for the remainders in \eqref{equation-u-hat}: 
\begin{align}\label{FEM-consistency-error-2}
	|\langle\mathscr{F}(t), v_h\rangle|+|\langle\mathscr{G}(t), v_h \rangle|\leq Ch^{r+1}\|\nabla v_h\|_{\Omega_h(t)} ,
\end{align}
which can be used to derive the error estimate of the semidiscrete FEM in Theorem \ref{THM-SD}.
\begin{proof}[Proof of Theorem \ref{THM-SD}]
	Let us define $e_h:= \hat{u}_h-u_h$ and $\lambda_h:= \hat{p}_h-p_h$. By subtracting equation \eqref{main-semidiscrete-a} from equation \eqref{equation-u-hat-a} and testing the resulting equation with $v_{h} = e_{h}$, we obtain:
\begin{align*}
	\frac{1}{2} \dt\|e_h(t)\|^2_{\Omega_h(t)}+\|\nabla e_h(t)\|^2_{\Omega_h(t)}&=\langle \mathscr{G}(t), e_h\rangle+\langle \mathscr{F}(t), e_h\rangle\nonumber\leq Ch^{r+1}\|\nabla e_h(t)\|_{\Omega_h(t)},
\end{align*}
where the last inequality follows from \eqref{FEM-consistency-error-2}.

By applying Young's inequality and absorbing $\|\nabla e_h\|_{\Omega_h(t)}^{2}$ on the right-hand side, we can integrate the inequality from $t=0$ to $t=s$ to obtain:
\begin{align*}
	\|e_h\|^2_{\Omega_h(s)}+\int_{0}^{s}\|\nabla e_h\|^2_{\Omega_h(t)}dt\leq \|e_h(0)\|^2_{\Omega_h(0)}+Ch^{2(r+1)}.
\end{align*}
Combining this result with Lemma \ref{error-stokes-ritz}, we complete the proof.
\end{proof}

Before we end this section, we furthermore need to present the estimates for $D_{t,h}^2(\hat{u}_h-\widetilde{u})$ and $D_{t,h}^2(\hat{p}_h-\widetilde{p})$, which will be used in the estimations of fully discrete scheme. According to Lemma \ref{H^1-error-1} and Lemma \ref{H^1-error-2}, we deduce the following lemma. 
\begin{lemma}\label{error-stokes-ritz-H^1-3}
	Let $ (v,q) = (E_{t,h}\widetilde{u},E_{t,h}\widetilde{p}) $ in Lemma \ref{H^1-error-3}, the errors $D_{t,h}^2(\hat{u}_h-\widetilde{u})$ and $D_{t,h}^2(\hat{p}_h-\widetilde{p})$ satisfy the following inequality:
	\begin{align}
		\|\nabla(D^2_{t,h}\hat{u}_h-D_{t,h}^2\widetilde{u})\|_{\Omega_h(t)}+\|D_{t,h}^2\hat{p}_h-D_{t,h}^2\widetilde{p}\|^{'}_{\Omega_h(t)} \leq\, Ch^r,
	\end{align}
where $ C $ is a constant independent of the mesh size $ h $ and time $ t $.
\end{lemma}
As a corollary, we obtain the following stability result.
\begin{align}\label{H^1-stability}
	\|\nabla \hat{u}_h\|_{\Omega_h(t)}+\|\nabla D_{t,h}\hat{u}_h\|_{\Omega_h(t)}+\|\nabla D^2_{t,h}\hat{u}_h\|_{\Omega_h(t)}+\|\nabla \hat{p}_h\|_{\Omega_h(t)}+\|\nabla D_{t,h}^2\hat{p}_h\|_{\Omega_h(t)}\leq \, C.
\end{align}

\section{ALE second-order projection method on evolving mesh}\label{second-fully-section}
In this section, we present the ALE second-order projection method on an evolving mesh, and then present an optimal-order error estimate for the method. Throughout this section, we adopt the following notations:
$$
w^n_h:= w_h(t_{n}), \quad \Omega_h^n:= \Omega_h(t_n), \quad \Phi_h^n:=\Phi_h(t_n), \quad \Omega^n:=\Omega(t_n), \quad \Phi^n:=\Phi(t_n).
$$ 

Let $ M = [T/\tau] $, the largest integer that is no larger than $ T/\tau $. For the fully discrete second-order projection method scheme of problem \eqref{main-equation-stokes} in the ALE formulation, we employ a non-conservative second-order time discretization. To simplify the notation, we introduce the discrete versions of $[\mathring V_h^r]$ and $[Q_h^{r-1}]$, which are defined in Section \ref{section-semidiscrete}. Specifically, we define
	\begin{align*}
	    u_{h}\sim v_{h} \;\text{iff}\; u_{h}\circ \Phi_{h}^{n}\circ (\Phi_{h}^{m})^{-1} = v_{h}, \;\forall u_{h} \in \mathring V_{h}^{r}(\Omega_{h}^{n}),\;v_{h} \in \mathring V_{h}^{r}(\Omega_{h}^{m}).
	\end{align*} 
	We use $ \mathring V_{h} $ to represent $ [\mathring V_{h}^{r}] $ and $ Q_{h} $ to represent $ [Q_{h}^{r-1}] $, respectively.

	On discrete time levels, we introduce the following notation for calculations on $\Omega_h^{n+1/2}$, where $u_h$ and $v_h$ are arbitrary functions in $\mathring V_h$:
\begin{align*}
(u_h, v_h)_{\Omega^{n+1/2}_h}&:=\frac{1}{2}\left((u_h, v_h)_{\Omega^n_h}+(u_h, v_h)_{\Omega^{n+1}_h}\right),\\
\left(\nabla u_h, \nabla v_h\right)_{\Omega^{n+1/2}_h}&:=\frac{1}{2}\left(\left(\nabla u_h, \nabla v_h\right)_{\Omega^{n+1}_h}+\left(\nabla u_h, \nabla v_h\right)_{\Omega^{n}_h}\right)
\end{align*}

Consider a sequence $\{u^n_h: u^n_h\in \mathring V_h\}_{n=0}^{T/\tau}$, where $ u_{h}^{n} $ represents the equivalence class in $ \mathring V_{h} $ generated by the numerical solution at time $ t_{n} $. We introduce the following notations:
\begin{align*}
u^{n+1}_h-u^n_h&:= \{u^{n+1}_h\}-\{u^n_h\}=\{u^{n+1}_h-u^n_h\circ \Phi_h^n\circ(\Phi_{h}^{n+1})^{-1}\}\in \mathring V_h,\\
u^{n+1/2}_h&:= \frac{\{u^{n+1}_h\}+\{u^n_h\}}{2}=\frac{1}{2}\left\{u^{n+1}_h+u^n_h\circ \Phi_h^n\circ(\Phi_{h}^{n+1})^{-1}\right\}\in \mathring V_h,\\
d_t u^{n+1}_h&:= \frac{\{u^{n+1}_h\}-\{u^n_h\}}{\tau}=\frac{1}{\tau}\left\{u^{n+1}_h-u^n_h\circ \Phi_h^n\circ(\Phi_{h}^{n+1})^{-1}\right\} \in \mathring V_h.
\end{align*}


We assume that the initial values of the finite element solutions are sufficiently accurate, satisfying the following estimate. 
\begin{assumptionp}{A}\label{assump-initial-data} $\,$\vspace{-10pt} 
	\begin{align}
		\|e_h^0\|_{\Omega_h^0}+\tau\|\nabla\eta_h^0\|_{\Omega_h^0}\leq C\tau^2+Ch^{r+1}.
	\end{align}
\end{assumptionp}
For example, if $u_{h}^{0}$ is chosen as the interpolation or Stokes--Ritz projection of the exact initial value $u_{0}$, and $p_{h}^{0}$ is determined using a general coupled second-order scheme for the first step, then the above-mentioned error estimate holds. 

\begin{theorem}\label{main-theorem}
Let $u$ and $p$ be sufficiently smooth solutions to the , and let Assumption \ref{assump-initial-data} hold. Then there exists a constant $\tau_0$ which is independent of the exact solutions and the mesh size $h$, such that for $\tau\le \tau_0$ the error $ e^{n} = u (t_{n}) - u_{h}^{n} $ satisfies the following estimate: 
\begin{align}
		\max_{0\leq n\leq M}\|e^{n}\|_{\Omega_h^{n}} \le \max_{0\leq n\leq M}\|u (t_{n}) - \hat{u}_{h}^{n}\|_{\Omega_h^{n}} + \max_{0\leq n\leq M}\|e_h^{n}\|_{\Omega_h^{n}} \leq C\min\{\ell_\tau,\ell_h\}(h^{r+1}+\tau^2),
\end{align}
where we have used the notations $\ell_\tau:=\ln\left(\frac{1}{\tau}+1\right)$ and $\ell_h:=\ln\left(\frac{1}{h}+1\right)$.	
\end{theorem}

With these notations, the second-order fully discrete projection method is formulated as follows: Find $u_h^{n+1}\in \mathring V_h$ and $p_h^{n+1}\in Q_h$ at step $n+1$ such that
\begin{subequations}\label{fully-scheme}
\begin{align}
	& \left(d_t u^{n+1}_h, v_h\right)_{\Omega^{n+1/2}_h}-(w^{n+1/2}_h \cdot \nabla u_h^{n+1/2}, v_h)_{\Omega^{n+1/2}_h}+(\nabla u^{n+1/2}_h, \nabla v_h)_{\Omega^{n+1/2}_h}-\left({\rm div} v_h, p_h^{n}\right)_{\Omega^{n+1/2}_h}\nonumber\\
	&=\frac{1}{2}\{(f(t_{n+1}), v_h)_{\Omega^{n+1}_h} + (f(t_{n}), v_h)_{\Omega^{n}_h}\}\quad \forall v_h\in \mathring V_h,\label{main-scheme}\\
	& ({\rm div}u^{n+1}_h, q_h)_{\Omega^{n+1}_h}+\beta\tau \left(\nabla(p^{n+1}_h-p^n_h), \nabla q_h\right)_{\Omega^{n+1}_h}=0 \quad \forall q_h\in Q_h,\label{p-scheme}
\end{align}
\end{subequations}
where $ \beta > 1 $ is a constant. The solution $ u^{n+1}_{h} $ is obtained by solving equation \eqref{main-scheme}, and subsequently, $ p_{h}^{n+1} $ is computed using equation \eqref{p-scheme} and $ u_{h}^{n+1} $.


\subsection{Analysis of the consistency error}
By denoting $ (\hat{u}_{h}(t_{n}), \hat{p}_{h}(t_{n})) $ as $ (\hat{u}_{h}^{n},\hat{p}_{h}^{n}) $ and $ (\hat{u}_{h}^{n}+\hat{u}_{h}^{n+1})/2 $ as $ \hat{u}_{h}^{n+1/2} $, we can substitute $ (\hat{u}_{h}^{n}, \hat{p}_{h}^{n}) $ into scheme \eqref{main-scheme} to replace $ (u_{h}^{n}, p_{h}^{n}) $. As a result, the pair $ (\hat{u}_{h}^{n},\hat{p}_{h}^{n}) $ satisfies the following equations, accompanied by remainder terms $ \mathscr{E}_{1}^{n+1} $ and $ \mathscr E_{2}^{n+1} $.
\begin{subequations}\label{Ritz-consistency-equation}
\begin{align}
	& \left(d_t \hat{u}^{n+1}_h, v_h\right)_{\Omega^{n+1/2}_h}-(w^{n+1/2}_h \cdot \nabla \hat{u}_h^{n+1/2}, v_h)_{\Omega^{n+1/2}_h}+(\nabla \hat{u}^{n+1/2}_h, \nabla v_h)_{\Omega^{n+1/2}_h}-\left({\rm div} v_h, \hat{p}_h^{n}\right)_{\Omega^{n+1/2}_h}\nonumber\\
	&=\frac{1}{2}\{(f(t_{n+1}), v_h)_{\Omega^{n+1}_h} + (f(t_{n}), v_h)_{\Omega^{n}_h}\} + \langle \mathscr{E}_1^{n+1}, v_h\rangle\quad \forall v_{h} \in \mathring V_{h},\label{Ritz-consistency-equation-a} \\
	& ({\rm div}\hat{u}^{n+1}_h, q_h)_{\Omega^{n+1}_h}+\beta\tau \left(\nabla(\hat{p}^{n+1}_h-\hat{p}^n_h), \nabla q_h\right)_{\Omega^{n+1}_h}= \langle \mathscr E_{2}^{n+1},q_{h}\rangle \quad \forall q_h\in Q_h.\label{Ritz-consistency-equation-b}
\end{align}
\end{subequations}
Since $ (\hat{u}_{h},\hat{p}_{h}) $ satisfies equation \eqref{equation-u-hat}, integrating \eqref{equation-u-hat-a} from $t=t_n$ to $t=t_{n+1}$ yields
\begin{align}\label{equation-u-hat-integrated}
	&\int_{t_n}^{t_{n+1}}(D_{t,h}\hat{u}_h, v_h)_{\Omega_h(t)}-(w_h\cdot \nabla \hat{u}_h, v_h)_{\Omega_h (t)}+(\nabla \hat{u}_h, \nabla v_h)_{\Omega_h (t)}-({\rm div} v_h, \hat{p}_h)_{\Omega_h (t)}\d t\notag\\
	=&\int_{t_n}^{t_{n+1}}(f,v_h)_{\Omega_h(t)}+\langle \mathscr{G}(t), v_h\rangle +\langle \mathscr{F}(t), v_h\rangle \; \d t. 
\end{align}
Subtracting \eqref{equation-u-hat-integrated} from \eqref{Ritz-consistency-equation-a}, we obtain the expression for the remainder $ \mathscr E_{1}^{n+1} $ as a linear operator on $ \mathring V_{h} $,
\begin{align}\label{expression-E-1-n}
    \mathscr E_{1}^{n+1} = Q_{0}^{n+1} + Q_{1}^{n+1} + Q_{2}^{n+1} + Q_{3}^{n+1} + Q_{4}^{n+1} + \mathscr F^{n+1} + \mathscr G^{n+1},
\end{align} 
where these error terms are defined as follows:
\begin{subequations}\label{remainder-error-E-1}
\begin{align}
	\tau\langle Q^{n+1}_0, v_h\rangle:= & (\hat{u}_{h}(t_{n+1})-\hat{u}_{h}(t_{n}),v_{h})_{\Omega_{h}^{n+1/2}} - \int_{t_n}^{t_{n+1}}\left(D_{t,h}\hat{u}_h(t), v_h\right)_{\Omega_h(t)}\d t,\label{Q_0-def}\\
	\tau\langle Q^{n+1}_1,v_h\rangle := & - \tau (w_h^{n+1/2}\cdot \nabla \hat{u}_h^{n+1/2}, v_h)_{\Omega_h^{n+1/2}}+\int_{t_n}^{t_{n+1}}\left(w_h\cdot \nabla \hat{u}_h, v_h \right)_{\Omega_h(t)}\d t,\\
	\tau\langle Q^{n+1}_2, v_h\rangle := &\tau (\nabla \hat{u}_h^{n+1/2}, \nabla v_h)_{\Omega_h^{n+1/2}}-\int_{t_n}^{t_{n+1}}\left(\nabla \hat{u}_h, \nabla v_h \right)_{\Omega_h(t)}\d t,\\
	\tau\langle Q^{n+1}_3, v_h\rangle :=&-\frac{\tau}{2}\left(\left(f(t_{n+1}), v_h\right)_{\Omega^{n+1}_h} + \left(f(t_{n}), v_h\right)_{\Omega^{n}_h}\right)+\int_{t_n}^{t_{n+1}}\left(f(t), v_h \right)_{\Omega_h(t)}\d t,\\
	\tau\langle Q^{n+1}_4, v_h\rangle:=&-\tau\left({\rm div}v_h, \hat{p}_h^n\right)_{\Omega_h^{n+1/2}}+\int_{t_n}^{t_{n+1}}\left({\rm div} v_h, \hat{p}_h\right)_{\Omega_h(t)}\d t, \label{def-E_2}\\
	\tau\langle \mathscr{F}^{n+1}, v_h\rangle:=& \int_{t_n}^{t_{n+1}}\langle \mathscr{F}(t), v_h\rangle \d t, \quad \tau\langle \mathscr{G}^{n+1}, v_h\rangle:= \int_{t_n}^{t_{n+1}}\langle \mathscr{G}(t), v_h\rangle \d t.\label{def-F-G-n+1}
\end{align}
\end{subequations}

The consistency errors $ Q_{0}^{n+1} $ to $ Q_{4}^{n+1} $ are rigorously established in the following lemma, employing Taylor's expansion as a key analytical tool.
\begin{lemma}\label{truncation-error-lemma}
	For sufficiently smooth solutions $ (\widetilde{u},\widetilde{p}) $, the following estimates hold.
	\begin{subequations}
	\begin{align}\label{Q-1-estimate}
		&|\langle Q^{n+1}_0, v_h\rangle|+|\langle Q^{n+1}_1, v_h\rangle|+|\langle Q^{n+1}_2, v_h\rangle|+|\langle Q^{n+1}_3, v_h\rangle|\leq C\tau^2 \|v_h\|_{\Omega_h^{n+1}},\\
		&|\langle Q^{n+1}_4, v_h\rangle|\leq C\tau \|v_h\|_{\Omega_h^{n+1}}.
	\end{align}
	\end{subequations}
Moreover, if $ v_{h} \in \mathring V_{h} $ satisfies $ (v_{h},q_{h})_{\Omega_{h}^{n+1/2}} = 0 $ for any $ q_{h} \in Q_{h} $, we have
\begin{align}\label{improve-Q4-estimate}
    |\langle Q^{n+1}_4, v_h\rangle|\leq C\tau^{2} \|v_h\|_{\Omega_h^{n+1}}.
\end{align} 
\end{lemma}
\begin{proof}
	We proceed with the proof by considering the inequality for $ Q_{2}^{n+1} $ as an example. We define an auxiliary function $g(t)=(\nabla\hat{u}_h,\nabla v_h)_{\Omega_h(t)}$. Utilizing Taylor's expansion, we obtain
\begin{align*}
	\Big|\int_{t_n}^{t_{n+1}} g(t)dt-\frac{1}{2}\tau (g(t_{n+1})+g(t_n))\Big|\leq C\tau^3\max_{t\in [t_n,t_{n+1}]}|g''(t)|.
\end{align*}
By computing the explicit formula for $g''(t)$ and introducing another function $h(t)$, which is obtained by replacing $\hat{u}_h$ with $\widetilde{u}$, $w_h$ with $w$, and $D_{t,h}$ with $D_t$ in the explicit formula of $g''(t)$, we can utilize Lemma \ref{error-Ritz-H1-optimal}, Lemma \eqref{error-stokes-ritz-H^1-3}, and assumption \eqref{assumption-w-phi} to establish
\begin{align}\label{g-second-h}
	|g''(t)-h(t)|\leq Ch\|\nabla v_h\|_{\Omega_h^{n+1}}\leq C\|v_h\|_{\Omega_h^{n+1}} \quad \mbox{(inverse estimates used)}.
\end{align}
Moreover, from the formula of $h(t)$, we can integrate by parts to move the gradient on $v_h$ to terms involving $\widetilde{u}$ and $w$. Thus, we have
\begin{align*}
	|h(t)|\leq C\|v_h\|_{\Omega_h^{n+1}}.
\end{align*}
Consequently, we have shown that
\begin{align*}
	\Big|\int_{t_n}^{t_{n+1}} g(t)\,\d t-\frac{1}{2}\tau (g(t_{n+1})+g(t_n))\Big|\leq C\tau^3\|v_h\|_{\Omega_h^{n+1}}.
\end{align*}
Next, we focus on 
\begin{align*}
	&\,\frac{1}{2}\tau \left(g(t_{n+1})+g(t_n)\right)-\tau(\nabla\hat{u}_h^{n+1/2}, \nabla v_h)_{\Omega_h^{n+1/2}}\\
	=&\,\frac{1}{2}\tau\Big[\frac{1}{2}(\nabla(\hat{u}_h^{n+1}-\hat{u}_h^n), \nabla v_h)_{\Omega_h^{n+1}}-\frac{1}{2}(\nabla(\hat{u}_h^{n+1}-\hat{u}_h^n), \nabla v_h)_{\Omega_h^{n}}\Big]\\
	=&\,\frac{\tau}{4}\int_{t_n}^{t_{n+1}}D_{t,h}(\nabla (\hat{u}_{h}^{n+1}-\hat{u}_{h}^{n}),\nabla v_{h})\;\d t.
\end{align*}
Similar to the proof in \eqref{g-second-h}, by replacing $D_{t,h}\hat{u}_h$ with $D_{t,h}\widetilde{u}$, $ \hat{u}_{h} $ with $ \widetilde{u} $, and $ w_{h} $ with $ w $, and using inverse estimates and Lemma \ref{error-Ritz-H1-optimal}, we obtain the desired result.

Regarding $ Q_{4}^{n+1} $, if $ v_{h} $ satisfies $ (v_{h},q_{h})_{\Omega_{h}^{n+1/2}} = 0 $ for all $ q_{h} \in Q_{h} $, then we can establish 
\begin{align*}
	\tau|\langle Q^{n+1}_4, v_h\rangle| = \Big|(v_{h},\hat{p}_{h}^{n+1/2})_{\Omega_{h}^{n+1/2}} - \int_{t_{n}}^{t_{n+1}}({\rm div} w_{h}^{n},\hat{p}_{h})_{\Omega_{h}(t)}\d t\Big| \le C \tau^{3} \| v_{h} \|_{\Omega_{h}^{n+1}},
\end{align*} 
employing the same method as the estimation of $ Q_{2}^{n+1} $.
\end{proof}
 
Considering the constraint $ ({\rm div}\hat{u}_{h}^{n+1},q_{h})_{\Omega_{h}^{n+1}} = 0 $ in equation \eqref{Ritz-consistency-equation-b}, we can readily obtain the estimate for $ \mathscr E_{2}^{n+1} $ as follows:
\begin{align}\label{E1-E3-estimate}
	|\langle \mathscr{E}_2^{n+1}, q_h\rangle|\leq C\beta\tau^2\|\nabla q_h\|_{\Omega_h^{n+1}}.
\end{align}
Furthermore, by utilizing \eqref{FEM-consistency-error-2} we can deduce estimates for $ \mathscr F^{n+1} $ and $ \mathscr G^{n+1} $, yielding
\begin{align}
    |\langle \mathscr{F}^{n+1}, v_h\rangle|+|\langle \mathscr{G}^{n+1}, v_h\rangle|\leq Ch^{r+1}\|\nabla v_h\|_{\Omega_h^{n+1}}.\label{F-G-estimate}
\end{align} 
As a corollary of Lemma \ref{truncation-error-lemma} and estimate \eqref{F-G-estimate}, we can derive the estimates for $ \mathscr E_{1}^{n+1} $:
\begin{subequations}
    \begin{align}
	    |\langle \mathscr E_{1}^{n+1},v_{h}\rangle| &\le C \tau \| v_{h} \|_{\Omega_{h}^{n+1}} + C h^{r+1}\| \nabla v_{h} \|_{\Omega_{h}^{n+1}},\label{E-1-estimate}\\
	    |\langle \mathscr E_{1}^{n+1},v_{h}\rangle| &\le C \tau^{2} \| v_{h} \|_{\Omega_{h}^{n+1}} + C  h^{r+1}\| \nabla v_{h} \|_{\Omega_{h}^{n+1}}, \;\text{if}\; (v_{h},q_{h})_{\Omega_{h}^{n+1/2}}=0 \; \forall q_{h} \in Q_{h}.\label{improve-E-1-estimate}
    \end{align}
\end{subequations}


\subsection{Error equations and the energy error estimates}
Let us define the errors between Stokes--Ritz projection and the numerical solutions as $e^{n}_h:=\hat{u}_h^n-u_h^n$ and $\eta_h^n:=\hat{p}_h^n-p_h^n$. By subtracting equation \eqref{fully-scheme} from \eqref{Ritz-consistency-equation}, we obtain the error equations for $ (e_{h}^{n+1},\eta_{h}^{n+1}) $ at step $n+1$ ($n\geq 0$) as follows:
\begin{subequations}\label{main-scheme-error-equation}
	\begin{align}
		\label{main-scheme-error-equation-1}
		& \left(d_t e^{n+1}_h, v_h\right)_{\Omega^{n+1/2}_h}+(\nabla e^{n+1/2}_h, \nabla v_h)_{\Omega^{n+1/2}_h}-(w^{n+1/2}_h \cdot \nabla e_h^{n+1/2}, v_h)_{\Omega^{n+1/2}_h}-\left({\rm div} v_h, \eta^{n}_h\right)_{\Omega^{n+1/2}_h}\nonumber\\
		&=\langle \mathscr{E}_1^{n+1}, v_h\rangle\quad \forall v_h\in \mathring V_{h},\\
		\label{main-scheme-error-equation-2}
		& ({\rm div}e^{n+1}_h, q_h)_{\Omega^{n+1}_h}+\beta\tau \left(\nabla(\eta^{n+1}_h-\eta^n_h), \nabla q_h\right)_{\Omega^{n+1}_h}=\langle \mathscr{E}_2^{n+1}, q_h\rangle \quad \forall q_h\in Q_{h}.
	\end{align}
\end{subequations}

Utilizing the error equations \eqref{main-scheme-error-equation}, we can employ an energy method to estimate the $ L^{2} $ norm of the error $ e_{h}^{n} $. Directly proving first-order convergence in time is demonstrated in the following lemma. However, to establish second-order convergence in time, further in-depth analysis is required, which will be presented in the subsequent subsections.
\begin{lemma}\label{energy-error-estimate-0}
Under the assumptions of Theorem \ref{main-theorem}, there exists a constant $\tau_0$, which is independent of the exact solutions and the mesh size $h$, such that for $\tau\le \tau_0$ the solution of \eqref{main-scheme-error-equation} satisfies the following energy estimate: 
	\begin{align}
		&\|e^{m+1}_h\|^2_{\Omega_h^{m+1}}+\sum_{n=0}^{m}\tau\|\nabla e_h^{n+1/2}\|^2_{\Omega_h^{n+1/2}}+\beta\tau^2\|\nabla \eta^{m+1}_h\|^2_{\Omega_h^{m+1}}\notag\\
		\leq & Ch^{2(r+1)}+C\tau^2 + C \| e_{h}^{0} \|_{\Omega_{h}^{0}}^{2}+C \tau^2\|\nabla \eta^{0}_h\|^2_{\Omega_h^0} .
	\end{align}
\end{lemma} 
\begin{proof}
	Testing \eqref{main-scheme-error-equation-1} with $v_h=e_h^{n+1/2}$ yields 
	\begin{align}
		&\frac{1}{2\tau}(\|e^{n+1}_h\|^2_{\Omega_h^{n+1/2}}-\|e^{n}_h\|^2_{\Omega_h^{n+1/2}})+\|\nabla e_h^{n+1/2}\|^2_{\Omega_h^{n+1/2}}-({\rm div}e_h^{n+1/2}, \eta_h^{n})_{\Omega_h^{n+1/2}}\nonumber\\
		=&(w_h^{n+1/2}\cdot \nabla e_h^{n+1/2}, e_h^{n+1/2})_{\Omega_h^{n+1/2}}+\langle \mathscr{E}_1^{n+1}, e_h^{n+1/2}\rangle.
	\end{align}
By utilizing the estimate \eqref{E-1-estimate}, we obtain the following inequality for $n\geq0$:
\begin{align}\label{energy-estimate-n-half}
	&\frac{1}{2\tau}\Big(\|e^{n+1}_h\|^2_{\Omega_h^{n+1/2}}-\|e^{n}_h\|^2_{\Omega_h^{n+1/2}}\Big)+\|\nabla e_h^{n+1/2}\|^2_{\Omega_h^{n+1/2}}-({\rm div}e_h^{n+1/2}, \eta_h^{n})_{\Omega_h^{n+1/2}}\nonumber\\
	\leq & C\|\nabla e_h^{n+1/2}\|_{\Omega_h^{n+1/2}}\|e_h^{n+1/2}\|_{\Omega_h^{n+1/2}}+Ch^{r+1}\|\nabla e_h^{n+1/2}\|_{\Omega_h^{n+1/2}}+C\tau \|e_h^{n+1/2}\|_{\Omega_h^{n+1/2}}.
\end{align}
Furthermore, from Lemma \ref{transport-theorem-lemma}, we have 
\begin{align}
	\Big|\|e_h^{n+1}\|^2_{\Omega_h^{n+1}}-\|e_h^{n+1}\|^2_{\Omega_h^{n+1/2}}\Big|\leq C\tau\|e_h^{n+1}\|^2_{\Omega_h^{n+1}},\quad 	\Big|\|e_h^{n}\|^2_{\Omega_h^{n}}-\|e_h^{n}\|^2_{\Omega_h^{n+1/2}}\Big|\leq C\tau\|e_h^{n}\|^2_{\Omega_h^{n}}.
\end{align}
To address the term $-({\rm div}e_h^{n+1/2}, \eta_h^{n})_{\Omega_h^{n+1/2}}$, we observe that
\begin{align}
	-({\rm div}e_h^{n+1/2}, \eta_h^{n})_{\Omega_h^{n+1/2}}
	=\frac{1}{2}\Big(-({\rm div}e_h^{n+1}, \eta_h^{n})_{\Omega_h^{n+1}}-({\rm div}e_h^{n}, \eta_h^{n})_{\Omega_h^n}\Big)+\mathscr{A}_1^{n},
\end{align}
where 
$$
2\mathscr{A}_1^{n}:=\left({\rm div}e_h^{n+1}, \eta_h^{n}\right)_{\Omega_h^{n+1}}-\left({\rm div}e_h^{n+1}, \eta_h^{n}\right)_{\Omega_h^{n+1/2}}+\left({\rm div}e_h^{n}, \eta_h^{n}\right)_{\Omega_h^{n}}-\left({\rm div}e_h^{n}, \eta_h^{n}\right)_{\Omega_h^{n+1/2}}.
$$ 
Applying Taylor's expansion, similar to the proof of Lemma \ref{truncation-error-lemma}, we obtain
\begin{align*}
	|\mathscr{A}_1^{n}|&\leq C\tau\left(\|e_h^{n+1}\|_{\Omega_h^{n+1}}+\|e_h^n\|_{\Omega_h^n}\right)\|\nabla\eta_h^n\|_{\Omega_h^n} \leq C \tau^2\|\nabla\eta_h^n\|^2_{\Omega_h^n}+C\left(\|e_h^{n+1}\|^2_{\Omega_h^{n+1}}+\|e_h^{n}\|^2_{\Omega_h^{n}}\right).
\end{align*}
For $n\geq 1$, by testing \eqref{main-scheme-error-equation-2} with $q_h=\eta_h^{n}$ for steps $n+1$ and $n$, respectively, we have
\begin{align}\label{divergence-estimate-e-h}
	&-\frac{1}{2}\left(\left({\rm div}e_h^{n+1}, \eta_h^{n}\right)_{\Omega_h^{n+1}}+\left({\rm div}e_h^{n}, \eta_h^{n}\right)_{\Omega_h^n}\right)\nonumber\\
	=&\frac{\beta \tau}{2}\left(\left(\nabla(\eta_h^{n+1}-\eta_h^n), \nabla\eta_h^{n}\right)_{\Omega_h^{n+1}}+\left(\nabla(\eta_h^{n}-\eta_h^{n-1}), \nabla\eta_h^{n}\right)_{\Omega_h^{n}}\right)-\frac{1}{2}\langle\mathscr{E}^{n+1}_2+\mathscr{E}^{n}_2, \eta_h^{n}\rangle,
\end{align}  
where the term $-\frac{1}{2}\langle\mathscr{E}^{n+1}_2+\mathscr{E}^{n}_2, \eta_h^{n}\rangle$ can be estimated using \eqref{E1-E3-estimate}, as follows: 
\begin{align}
	\Big|-\frac{1}{2}\langle\mathscr{E}^{n+1}_2+\mathscr{E}^{n}_2, \eta_h^{n}\rangle\Big| \leq C\beta\tau^2\|\nabla\eta_h^n\|_{\Omega_h^n}\leq C\tau^2\|\nabla \eta_h^n\|_{\Omega_h^n}^2+C\tau^2.
\end{align}
Additionally, we have the relation
\begin{align*}
	&\left(\nabla(\eta_h^{n+1}-\eta_h^n), \nabla\eta_h^{n}\right)_{\Omega_h^{n+1}}+\left(\nabla(\eta_h^{n}-\eta_h^{n-1}), \nabla\eta_h^{n}\right)_{\Omega_h^{n}}\nonumber\\
	=&\frac{1}{2}\left(\|\nabla \eta^{n+1}_h\|^2_{\Omega_h^{n+1}}-\|\nabla \eta^{n-1}_h\|^2_{\Omega_h^{n-1}}\right)-\frac{1}{2}\left(\|\nabla(\eta_h^{n+1}-\eta_h^n)\|^2_{\Omega_h^{n+1}}-\|\nabla(\eta_h^{n}-\eta_h^{n-1})\|^2_{\Omega_h^{n}}\right)+\mathscr{A}_2^{n},
\end{align*}
where $\mathscr{A}_2^{n}:=\frac{1}{2}(-\|\nabla \eta^{n}_h\|^2_{\Omega_h^{n+1}}+\|\nabla \eta^{n}_h\|^2_{\Omega_h^{n}}+\|\nabla \eta^{n-1}_h\|^2_{\Omega_h^{n-1}}-\|\nabla \eta^{n-1}_h\|^2_{\Omega_h^{n}})$ can be estimated using a similar approach to the proof of Lemma \ref{truncation-error-lemma}:
\begin{align}\label{lemma-4.5-temp2}
	|\mathscr{A}_2^{n}|\leq C\tau\left(\|\nabla\eta_h^n\|^2_{\Omega_h^n}+\|\nabla\eta_h^{n-1}\|^2_{\Omega_h^{n-1}}\right).
\end{align}
Combining estimates \eqref{energy-estimate-n-half}-\eqref{lemma-4.5-temp2}, we have
\begin{align}\label{lemma-4.5-temp3}
	&\frac{1}{2\tau}\left(\|e^{n+1}_h\|^2_{\Omega_h^{n+1}}-\|e^{n}_h\|^2_{\Omega_h^{n}}\right)+\frac{1}{2}\|\nabla e_h^{n+1/2}\|^2_{\Omega_h^{n+1/2}}+\frac{\beta\tau}{4}\left(\|\nabla \eta^{n+1}_h\|^2_{\Omega_h^{n+1}}-\|\nabla \eta^{n-1}_h\|^2_{\Omega_h^{n-1}}\right)\nonumber\\
	&\leq C\left(\|e_h^{n+1}\|^2_{\Omega_h^{n+1}}+\|e_h^{n}\|^2_{\Omega_h^{n}}\right)+\frac{\beta\tau}{4}\left(\|\nabla(\eta_h^{n+1}-\eta_h^n)\|^2_{\Omega_h^{n+1}}-\|\nabla(\eta_h^{n}-\eta_h^{n-1})\|^2_{\Omega_h^{n}}\right)\nonumber\\
	&\quad +Ch^{2(r+1)}+C\tau^2+C \tau^2\left(\|\nabla\eta_h^n\|^2_{\Omega_h^n}+\|\nabla\eta_h^{n-1}\|^2_{\Omega_h^{n-1}}\right).
\end{align}
By testing equation \eqref{main-scheme-error-equation-2} with $q_h=\eta_h^{m+1}-\eta_h^m$ at step $m+1$, we deduce the inequality:
\begin{align}\label{lemma-4.5-temp5}
	\beta\tau\|\nabla(\eta_h^{m+1}-\eta_h^m)\|_{\Omega_h^{m+1}}\leq \|e_h^{m+1}\|_{\Omega_h^{m+1}}+C\beta\tau^2.
\end{align}
By summing up \eqref{lemma-4.5-temp3} from $n=1$ to $n=m$ and incorporating inequality \eqref{lemma-4.5-temp5} and the fact $ \beta > 1 $, we obtain

\begin{align}\label{e-h-m-estimate-mid}
	&\frac{1}{4}\|e^{m+1}_h\|^2_{\Omega_h^{m+1}}+\sum_{n=1}^{m}\tau\|\nabla e_h^{n+1/2}\|^2_{\Omega_h^{n+1/2}}+\frac{\beta\tau^2}{2}\left(\|\nabla \eta^{m+1}_h\|^2_{\Omega_h^{m+1}}+\|\nabla \eta^{m}_h\|^2_{\Omega_h^{m}}\right)\nonumber\\
	\le & C\sum_{n=1}^m\beta\tau^3\left(\|\nabla\eta_h^n\|^2_{\Omega_h^n}+\|\nabla\eta_h^{n-1}\|^2_{\Omega_h^{n-1}}\right)+\|e^{1}_h\|^2_{\Omega_h^{1}}+\frac{\beta\tau^2}{2}\left(\|\nabla \eta^{1}_h\|^2_{\Omega_h^{1}}+\|\nabla \eta^{0}_h\|^2_{\Omega_h^0}\right)\notag\\
	&+ C\sum_{n=1}^m\tau\left(\|e_h^{n+1}\|^2_{\Omega_h^{n+1}}+\|e_h^{n}\|^2_{\Omega_h^{n}}\right)+Ch^{2(r+1)}+C\tau^2.
\end{align}
By applying the discrete Gronwall's inequality, it follows that for sufficiently small $ \tau $:
\begin{align}\label{lemma-4.5-temp6}
		&\|e^{m+1}_h\|^2_{\Omega_h^{m+1}}+\sum_{n=1}^{m}\tau\|\nabla e_h^{n+1/2}\|^2_{\Omega_h^{n+1/2}}+\beta\tau^2\left(\|\nabla \eta^{m+1}_h\|^2_{\Omega_h^{m+1}}+\|\nabla \eta^{m}_h\|^2_{\Omega_h^{m}}\right)\nonumber\\
	\leq & Ch^{2(r+1)}+C\tau^2+C\|e^{1}_h\|^2_{\Omega_h^{1}}+C\beta\tau^2\left(\|\nabla \eta^{1}_h\|^2_{\Omega_h^{1}}+\|\nabla \eta^{0}_h\|^2_{\Omega_h^0}\right).
\end{align}
To complete the proof, we need to estimate $\|e^{1}_h\|_{\Omega_h^{1}}$ and $\|\nabla \eta^{1}_h\|_{\Omega_h^{1}}$ using additional techniques. By setting $ n = 0 $ in \eqref{energy-estimate-n-half}, we obtain the following estimate:
\begin{align}\label{e-h-1-estimate}
&\frac{1}{2\tau}(\|e^{1}_h\|^2_{\Omega_h^{1/2}}-\|e^{0}_h\|^2_{\Omega_h^{1/2}})+\|\nabla e_h^{1/2}\|^2_{\Omega_h^{1/2}}-({\rm div}e_h^{1/2}, \eta_h^0)_{\Omega_h^{1/2}}\nonumber\\
	\leq & C\|\nabla e_h^{1/2}\|_{\Omega_h^{1/2}}\|e_h^{1/2}\|_{\Omega_h^{1/2}}+Ch^{r+1}\|\nabla e_h^{1/2}\|_{\Omega_h^{1/2}}+C\tau\| e_h^{1/2}\|_{\Omega_h^{1/2}}.
\end{align}
Next, by employing Young's inequality, we have
\begin{align*}
	C\tau \| e_h^{1/2}\|_{\Omega_h^{1/2}}\leq C\tau (\| e_h^{1}\|_{\Omega_h^{1}}+\| e_h^{0}\|_{\Omega_h^{0}})\leq \frac{1}{16\tau}(\| e_h^{1}\|^2_{\Omega_h^{1}}+\| e_h^{0}\|^2_{\Omega_h^{0}})+C\tau^3.
\end{align*}
Using a similar technique as in the estimation of \eqref{divergence-estimate-e-h}, we can prove that
\begin{align*}
    \Big|\frac{\beta \tau}{4}\| \nabla \eta_{h}^{1} \|_{\Omega_{h}^{1}}^{2}+(\text{div}e_{h}^{1/2},\eta_{h}^{0})_{\Omega_{h}^{1/2}}\Big| \le C \tau^{-1}\| e_{h}^{0} \|_{\Omega_{h}^{0}}^{2} + C \beta \tau \| \nabla \eta_{h}^{0} \|_{\Omega_{h}^{0}}^{2} +( C \tau +\frac{5}{16 \tau}) \| e_{h}^{1} \|_{\Omega_{h}^{1}}^{2} + C \tau^{3}.
\end{align*} 
By handling the remaining terms of \eqref{e-h-1-estimate} in the same manner as in the proof of \eqref{lemma-4.5-temp3} using H\"older's inequality, we obtain
\begin{align}\label{middle-estimate-eh-1}
	&\|e_h^1\|^2_{\Omega_h^1}+\tau\|\nabla e_h^{1/2}\|^2_{\Omega_h^{1/2}}+\beta\tau^2\|\nabla\eta_h^1\|^2_{\Omega_h^1}\notag\\
	\leq& C\|e_h^0\|^2_{\Omega_h^0}+C\tau	\|e_h^1\|^2_{\Omega_h^1}+Ch^{2(r+1)} +C\tau^4+C\tau^2\beta\|\nabla\eta_h^0\|^2_{\Omega_h^0}.
\end{align} 
When $\tau$ is sufficiently small, the term $C\tau\|e_h^1\|^2_{\Omega_h^1}$ can be absorbed by the left-hand side. In accordance with \eqref{lemma-4.5-temp6}, this completes the proof.
\end{proof}

Under the assumptions of Theorem \ref{main-theorem}, along with the energy estimate provided in Lemma \ref{energy-error-estimate-0} and the estimate \eqref{middle-estimate-eh-1}, we can immediately derive the following result.  
\begin{lemma}\label{energy-error-estimate-1}
Under the assumptions of Theorem \ref{main-theorem}, the following estimates hold:
	\begin{align}
		\|e^{m+1}_h\|_{\Omega_h^{m+1}}+\Big(\sum_{n=0}^{m}\tau\|\nabla e_h^{n+1/2}\|^2_{\Omega_h^{n+1/2}}\Big)^{\frac{1}{2}}+\tau\|\nabla \eta^{m+1}_h\|_{\Omega_h^{m+1}} &\leq C\tau  + Ch^{r+1},\\
		\|e_h^1\|_{\Omega_h^1}+\tau\|\nabla\eta_h^1\|_{\Omega_h^1} & \leq C\tau^2+Ch^{r+1}.
	\end{align}
\end{lemma} 

\subsection{Estimates for $\|e_h^{n+1}-e_h^n\|_{\Omega_h^{n+1}}$ and $\|\nabla(\eta_h^{n+1}-\eta_h^n)\|_{\Omega_h^{n+1}}$}
In this subsection, we aim to provide more refined estimates for $ \|e_h^{n+1}-e_h^n\|_{\Omega_h^{n+1}} $ and  $\|\nabla(\eta_h^{n+1}-\eta_h^n)\|_{\Omega_h^{n+1}}$, with the objective of enhancing the convergence order of our numerical scheme.

To simplify the notations, we introduce a difference operator $ \delta $ defined as follows:
\begin{align*}
	\delta v_h^{n+1}:=v_h^{n+1}-v_h^n,& \; \; \delta v_{h}^{n+1/2} := v_{h}^{n+1/2} - v_{h}^{n-1/2},\\
	(v_{h},w_{h})_{\delta \Omega_{h}^{n+1}}: = (v_{h},w_{h})_{\Omega_{h}^{n+1}} - (v_{h},w_{h})_{\Omega_{h}^{n}},&\; \; (v_{h},w_{h})_{\delta \Omega_{h}^{n+1/2}}: = (v_{h},w_{h})_{\Omega_{h}^{n+1/2}} - (v_{h},w_{h})_{\Omega_{h}^{n-1/2}},
\end{align*} 
where $ v_{h} $ and $ w_{h} $ can be any functions in $ \mathring V_{h} $ or $ Q_{h} $, respectively.

For each $n\geq 1$, by taking the difference of \eqref{main-scheme-error-equation} at step $n+1$ and step $n$, we obtain:
\begin{subequations}\label{main-scheme-error-equation-diff}
	\begin{align}
		\label{main-scheme-error-equation-diff-1}
		& (d_t \delta e^{n+1}_h, v_h)_{\Omega^{n+1/2}_h}+(\nabla \delta e^{n+1/2}_h, \nabla v_h)_{\Omega^{n+1/2}_h}-(w^{n+1/2}_h \cdot \nabla \delta e_h^{n+1/2}, v_h)_{\Omega^{n+1/2}_h}\nonumber\\
		=&\left({\rm div} v_h, \delta \eta^{n}_h\right)_{\Omega^{n+1/2}_h}+\langle \delta\mathscr{E}_1^{n+1}, v_h\rangle+\langle\mathscr{A}_4^{n+1}, v_h\rangle\\
		\label{main-scheme-error-equation-diff-2}
		& ({\rm div}\delta e^{n+1}_h, q_h)_{\Omega^{n+1}_h}+\beta\tau \left(\nabla(\delta\eta^{n+1}_h-\delta\eta^n_h), \nabla q_h\right)_{\Omega^{n+1}_h}=\langle \delta\mathscr{E}_2^{n+1}, q_h\rangle+\langle\mathscr{A}_5^{n+1}, q_h\rangle ,
	\end{align}
\end{subequations}
for all $ v_{h} \in \mathring V_{h} $ and $ q_{h} \in Q_{h} $. Here, $\mathscr{A}^{n+1}_4$ and $\mathscr{A}_5^{n+1}$ are defined as 
\begin{align*}
	\langle\mathscr{A}^{n+1}_4, v_h\rangle:=&-\left(d_t e^{n}_h, v_h\right)_{\delta\Omega^{n+1/2}_h}-(\nabla e^{n-1/2}_h, \nabla v_h)_{\delta\Omega^{n+1/2}_h}+(w^{n-1/2}_h \cdot \nabla e_h^{n-1/2}, v_h)_{\delta\Omega^{n+1/2}_h}\\
						&+\left({\rm div} v_h, \eta^{n-1}_h\right)_{\delta\Omega^{n+1/2}_h} + (\delta w_{h}^{n+1/2}\cdot \nabla e_{h}^{n-1/2}, v_{h})_{\Omega_{h}^{n+1/2}},\\
	\langle\mathscr{A}_5^{n+1}, q_h\rangle:=&-({\rm div}e^{n}_h, q_h)_{\delta\Omega^{n+1}_h}-\beta\tau \left(\nabla(\eta^{n}_h-\eta^{n-1}_h), \nabla q_h\right)_{\delta\Omega^{n+1}_h}.
\end{align*}

The consistency errors $\delta \mathscr E_{1}^{n+1}$ and $\delta\mathscr E_{2}^{n+1}$ are readily established through the application of Taylor's expansion, as demonstrated in Lemma \ref{truncation-error-lemma}. 
\begin{lemma}\label{delta-EFG-estimates}
	For sufficiently smooth solutions $u$ and $p$, the quantities $\delta \mathscr E_{1}^{n+1}$ and $\delta\mathscr E_{2}^{n+1}$ satisfy the following estimates for any $v_{h}\in \mathring V_{h}$ and $q_{h} \in Q_{h}$:
	\begin{align}		
		|\langle \delta \mathscr{E}_1^{n+1}, v_h\rangle|&\leq C\tau^2\| v_h\|_{\Omega_h^{n+1}} + Ch^{r+1}\|\nabla v_h\|_{\Omega_h^{n+1}}, \label{deltaE1-E2-estimate}\\
		|\langle \delta\mathscr{E}_2^{n+1}, q_h\rangle|&\leq C\tau^3\|\nabla q_h\|_{\Omega_h^{n+1}}.\label{deltaE3-estimate}
	\end{align}
\end{lemma}
Furthermore, by using Lemma \ref{transport-theorem-lemma}, we obtain the estimates for $\mathscr{A}^{n+1}_4$ and $\mathscr{A}_5^{n+1}$:
\begin{align}
	|\langle\mathscr{A}^{n+1}_4, v_h\rangle|\leq& C\|v_h\|_{\Omega_h^{n+1}}\|\delta e_h^n\|_{\Omega_h^n}+C\tau\|\nabla e_h^{n-1/2}\|_{\Omega_h^{n-1/2}}\|\nabla v_h\|_{\Omega_h^{n+1/2}}\notag\\
	&+C\tau\|\nabla v_h\|_{\Omega_h^{n+1/2}}\|\eta_h^{n-1}\|_{\Omega_h^{n-1}},\label{mscraA1^n+1-estimate}\\
	|\langle\mathscr{A}^{n+1}_5, q_h\rangle|\leq& C\tau\|e_h^n\|_{\Omega_h^n}\|\nabla q_h\|_{\Omega_h^n}+C\tau^2\|\nabla \delta \eta_h^n\|_{\Omega_h^n}\|\nabla q_h\|_{\Omega_h^n}.\label{mscraA2^n+1-estimate}
\end{align}
To estimate $\|\eta_h^n\|_{\Omega_h^n}$, we exploit the inf-sup condition and the following relation 
$$
|\left({\rm div}v_h, q_h\right)_{\Omega_h^n}-\left({\rm div}v_h, q_h\right)_{\Omega_h^{n+1/2}}|\leq C\tau\|\nabla v_h\|_{\Omega_h^n}\|q_h\|_{\Omega_h^n}. 
$$
This leads to the inequality: 
\begin{align}\label{inf-sup-eg1}
	\|q_h\|_{\Omega_h^n}\leq C\sup_{v_h\in \mathring V_h}\frac{\Big|({\rm div}v_h, q_h)_{\Omega_h^{n+1/2}}\Big|}{\|\nabla v_h\|_{\Omega_h^n}},\; \forall q_{h} \in Q_{h},
\end{align}
which holds when $\tau$ is sufficiently small. Consequently, using \eqref{main-scheme-error-equation-1}, we deduce that:
\begin{align}\label{eta-error-L2}
	\|\eta_h^n\|_{\Omega_h^n}\leq C\left[\|d_t e_h^{n+1}\|_{\Omega_h^{n+1}}+\|\nabla e_h^{n+1/2}\|_{\Omega_h^{n+1/2}}\right]+Ch^{r+1}+C\tau.
\end{align}
By substituting \eqref{eta-error-L2} into \eqref{mscraA1^n+1-estimate}, we obtain the following estimate:
\begin{align}\label{mscraA1^n+1-estimate1}
		|\langle\mathscr{A}^{n+1}_4, v_h\rangle|\leq C\left[\|\delta e_h^n\|_{\Omega_h^n}+\tau\|\nabla e_h^{n-1/2}\|_{\Omega_h^{n-1/2}}+\tau(\tau +h^{r+1})\right]\|\nabla v_h\|_{\Omega_h^{n+1/2}}.
\end{align}

Having established the necessary groundwork, we are now poised to present the energy estimates for $\delta e_h^{n+1}$ and $\delta\eta_h^{n+1}$, which demonstrate second-order convergence in time.
\begin{lemma}\label{energy-error-estimate-2}
Under the assumptions of Theorem \ref{main-theorem}, the following inequality holds:
	\begin{align}\label{delta-e-estimate}
		&\|\delta e^{m+1}_h\|_{\Omega_h^{m+1}}+\Big(\sum_{n=1}^{m}\tau\|\nabla \delta e_h^{n+1/2}\|^2_{\Omega_h^{n+1/2}}\Big)^{\frac{1}{2}}+\tau\|\nabla \delta \eta^{m+1}_h\|_{\Omega_h^{m+1}}\leq Ch^{r+1}+C\tau^2 .
	\end{align}
\end{lemma} 
\begin{proof}
By choosing $ v_{h} = \delta e_{h}^{n+1/2} $ in equation \eqref{mscraA1^n+1-estimate1} and applying Young's inequality, we obtain
\begin{align*}
	|\langle\mathscr{A}^{n+1}_4, \delta e_h^{n+1/2}\rangle|\leq C\|\delta e_h^n\|^2_{\Omega_h^n}+C\tau^2\Big(\|\nabla e_h^{n-1/2}\|^2_{\Omega_h^{n-1/2}}+\tau^2 +h^{2(r+1)}\Big)+\frac{1}{4}\|\nabla \delta e_h^{n+1/2}\|^2_{\Omega_h^{n+1/2}}.
\end{align*}
Next, by testing \eqref{main-scheme-error-equation-diff-1} with $v_h=\delta e_h^{n+1/2}$ and following a similar procedure as in the proof of Lemma \ref{energy-error-estimate-0}, utilizing Lemma \ref{delta-EFG-estimates}, we deduce the following inequality for each $n\geq 1$:
	\begin{align}\label{lemma-4.9-temp1}
		&\frac{1}{2\tau}\left(\|\delta e^{n+1}_h\|^2_{\Omega_h^{n+1}}-\|\delta e^{n}_h\|^2_{\Omega_h^{n}}\right)+\frac{1}{2}\|\nabla \delta e_h^{n+1/2}\|^2_{\Omega_h^{n+1/2}}-({\rm div}\delta e_h^{n+1/2}, \delta\eta_h^{n})_{\Omega_h^{n+1/2}}\nonumber\\
		\leq & C\left(\|\delta e_h^{n+1}\|^2_{\Omega_h^{n+1}}+\|\delta e_h^{n}\|^2_{\Omega_h^{n}}\right)+C\tau^2\|\nabla e_h^{n-1/2}\|^2_{\Omega_h^{n-1/2}}+Ch^{2(r+1)}+C\tau^4.
	\end{align}
Following from the proof in Lemma \ref{energy-error-estimate-0} to estimate $- ({\rm div}\delta e_{h}^{n+1/2},\delta \eta_{h}^{n})_{\Omega_{h}^{n+1/2}} $, the procedure is the same except for the inclusion of extra terms $-\frac{1}{2}\langle\mathscr{A}_5^{n+1}+\mathscr{A}_5^n, \delta\eta_h^n\rangle$ (for $n\geq 1$). By utilizing \eqref{mscraA2^n+1-estimate}, the estimate of $ \| e_{h}^{n} \|_{\Omega_{h}^{n}} $ proved in Lemma \ref{energy-error-estimate-1} and Young's inequality, these extra terms can be bounded as follows:
\begin{align*}
	|\frac{1}{2}\langle\mathscr{A}^{n+1}_5+\mathscr{A}_5^n, \delta \eta_h^n\rangle|&\leq C\tau^2\Big(\tau^2+h^{2(r+1)}\Big)+C\tau^2\Big(\|\nabla \delta \eta_h^n\|^2_{\Omega_h^n}+\|\nabla \delta \eta^{n-1}_h\|_{\Omega_h^{n-1}}\Big).
\end{align*}
Hence, in accordance with \eqref{deltaE3-estimate}, we obtain the following estimate for $ n \ge 2 $:
\begin{align}\label{lemma-4.9-temp2}
	&\frac{1}{2\tau}\Big(\|\delta e^{n+1}_h\|^2_{\Omega_h^{n+1}}-\|\delta e^{n}_h\|^2_{\Omega_h^{n}}\Big)+\frac{1}{2}\|\nabla \delta e_h^{n+1/2}\|^2_{\Omega_h^{n+1/2}}+\frac{\beta\tau}{4}\Big(\|\nabla \delta\eta^{n+1}_h\|^2_{\Omega_h^{n+1}}-\|\nabla \delta\eta^{n-1}_h\|^2_{\Omega_h^{n-1}}\Big)\nonumber\\
	\leq & C\Big(\|\delta e_h^{n+1}\|^2_{\Omega_h^{n+1}}+\|\delta e_h^{n}\|^2_{\Omega_h^{n}}\Big)+\frac{\beta\tau}{4}\Big(\|\nabla(\delta\eta_h^{n+1}-\delta\eta_h^n)\|^2_{\Omega_h^{n+1}}-\|\nabla(\delta\eta_h^{n}-\delta\eta_h^{n-1})\|^2_{\Omega_h^{n}}\Big)\nonumber\\
	&+Ch^{2(r+1)}+C\tau^4+C\tau^2\Big(\|\nabla\delta\eta_h^n\|^2_{\Omega_h^n}+\|\nabla\delta\eta_h^{n-1}\|^2_{\Omega_h^{n-1}} + \|\nabla e_h^{n-1/2}\|^2_{\Omega_h^{n-1/2}}\Big).
\end{align}

Another estimate for $n=1$ can be obtained by estimating $ \langle\mathscr{E}_2^1, \delta \eta_h^1\rangle $ as follows:
\begin{align}
	|\langle\mathscr{E}_2^1, \delta \eta_h^1\rangle|=|\left(\delta e_h^1, \delta \nabla \eta_h^1\right)_{\Omega_h^1}| + \beta \tau \| \nabla \delta \eta_{h}^{1} \|_{\Omega_{h}^{1}}^{2}\leq C\tau^3+Ch^{2(r+1)}\tau^{-1}.
\end{align}
Using the same method as in the derivation of \eqref{lemma-4.9-temp2}, we can obtain an estimate for $n=1$:
\begin{align}\label{deduction-n-1}
	&\frac{1}{2\tau}\left(\|\delta e^{2}_h\|^2_{\Omega_h^{2}}-\|\delta e^{1}_h\|^2_{\Omega_h^{1}}\right)+\frac{1}{2}\|\nabla \delta e_h^{3/2}\|^2_{\Omega_h^{3/2}}+\frac{\beta\tau}{4}\left(\|\nabla \delta\eta^{2}_h\|^2_{\Omega_h^{2}}-\|\nabla \delta\eta^{1}_h\|^2_{\Omega_h^{1}}\right)\nonumber\\
	\leq & C\left(\|\delta e_h^{2}\|^2_{\Omega_h^{2}}+\|\delta e_h^{1}\|^2_{\Omega_h^{1}}\right)+\frac{\beta\tau}{4}\|\nabla(\delta\eta_h^{2}-\delta\eta_h^1)\|^2_{\Omega_h^{2}} +Ch^{2(r+1)}\tau^{-1}+C\tau^3\notag\\
	     & +C\tau^2\left(\|\nabla\delta\eta_h^1\|^2_{\Omega_h^1}+\|\nabla e_h^{1/2}\|^2_{\Omega_h^{1/2}}\right).
\end{align}
By summing up the inequalities \eqref{lemma-4.9-temp2} and \eqref{deduction-n-1} from $n=1$ to $n=m$ and utilizing the results of Lemma \ref{energy-error-estimate-1}, we obtain the following estimate for $m\geq 1$:
\begin{align*}
	&\|\delta e^{m+1}_h\|^2_{\Omega_h^{m+1}}+\sum_{n=1}^{m}\tau\|\nabla \delta e_h^{n+1/2}\|^2_{\Omega_h^{n+1/2}}+\frac{\beta\tau^2}{2}\left(\|\nabla \delta\eta^{m+1}_h\|^2_{\Omega_h^{m+1}}+\|\nabla \delta \eta^{m}_h\|^2_{\Omega_h^{m}}\right)\nonumber\\
	\leq& C\sum_{n=1}^m\tau\left(\|\delta e_h^{n+1}\|^2_{\Omega_h^{n+1}}+\|\delta e_h^{n}\|^2_{\Omega_h^{n}}\right)+
	\frac{\beta\tau^2}{2}\|\nabla(\delta \eta_h^{m+1}-\delta \eta_h^m)\|^2_{\Omega_h^{m+1}}+Ch^{2(r+1)}+C\tau^4\nonumber\\
	&+C\sum_{n=1}^m\tau^3\left(\|\nabla\delta\eta_h^n\|^2_{\Omega_h^n}+\|\nabla\delta\eta_h^{n-1}\|^2_{\Omega_h^{n-1}}\right)+\|\delta e^{1}_h\|^2_{\Omega_h^{1}}+C\tau^2\|\nabla \delta\eta^{1}_h\|^2_{\Omega_h^{1}}.
\end{align*}
Finally, by following the same procedure as in the proof of Lemma \ref{energy-error-estimate-0} (the part of the proof after \eqref{e-h-m-estimate-mid}), the desired result is obtained.
\end{proof}

\subsection{The time-reversed parabolic dual problem}
To prove the optimal $\ell^\infty L^2$-norm error estimate, we consider the fully discrete time-reversed parabolic dual problem. In each step, we solve $(\psi_h^n,\varphi_h^n)\in \mathring V_h\times Q_h$ satisfying the equations:
\begin{subequations}\label{duality-form}
\begin{align}
	\label{duality-form-1}
	&\left(d_t \psi^{n+1}_h, v_h\right)_{\Omega_h^{n+1/2}}-(w^{n+1/2}_h\cdot \nabla\psi^n_h, v_h)_{\Omega^{n+1/2}_h}-\left(\nabla \psi_h^n,\nabla v_h\right)_{\Omega_h^{n+1/2}}-({\rm div}v_h, \varphi^n_h)_{\Omega^{n+1/2}_h} \nonumber\\
	=&(g^n_h, v_h)_{\Omega^{n+1/2}_h} \quad \forall v_h\in \mathring V_h,\\
	\label{duality-form-2}
	&\left({\rm div}\psi^n_h, q_h\right)_{\Omega^{n+1/2}_h}=0 \quad \forall q_h\in Q_h,
\end{align}
\end{subequations}
where $g_{h}^{n}$ is a given function, and the initial value is $\psi^{M}_h=0$. By testing equation \eqref{duality-form-1} with $v_h=e_h^{n+1/2}$, we define:
\begin{align}
	\mathcal{I}_1:=  (d_t\psi^{n+1}_h, e^{n+1/2}_h)_{\Omega^{n+1/2}_h}-(w^{n+1/2}_h\cdot\nabla\psi^n_h, e^{n+1/2}_h)_{\Omega^{n+1/2}_h},\\
	\mathcal{I}_2:= -(\nabla \psi^n_h,\nabla e^{n+1/2}_h)_{\Omega^{n+1/2}_h},\quad \mathcal{I}_3:=-({\rm div} e^{n+1/2}_h, \varphi^n_h)_{\Omega^{n+1/2}_h}.
\end{align}
The equation \eqref{duality-form-1} can be rewritten as
\begin{align}\label{formula-g-h-n}
(g^n_h, e^{n+1/2}_h)_{\Omega^{n+1/2}_h}=\mathcal{I}_1+\mathcal{I}_2+\mathcal{I}_3.
\end{align}
It is worth noting that
\begin{equation*}
	\mathcal{I}_1=\frac{1}{\tau}\left[(\psi^{n+1}_h, e^{n+1}_h)_{\Omega^{n+1}_h}-(\psi^{n}_h, e^{n}_h)_{\Omega^{n}_h}\right]-(\psi^n_h,d_t e^{n+1}_h)_{\Omega^{n+1/2}_h}+(w^{n+1/2}_h\cdot \nabla e^{n+1/2}_h,\psi^n_h)_{\Omega^{n+1/2}_h}+\mathscr{A}_6^n,
\end{equation*}
where $\mathscr{A}_6^n := A_{1}^{n} + A_{2}^{n} + A_{3}^{n}$ represents the remainder term, and $A_{1}^{n}$, $A_{2}^{n}$, and $A_{3}^{n}$ are defined by:
\begin{align*}
	A_1^{n}:=&\frac{1}{\tau}\left[(\psi^{n}_h, e^{n+1/2}_h)_{\Omega^{n}_h}-(\psi^{n}_h, e^{n+1/2}_h)_{\Omega^{n+1}_h}\right]-(w^{n+1/2}_h\cdot \nabla e^{n+1/2}_h,\psi^n_h)_{\Omega^{n+1/2}_h}\\
	     &-(w^{n+1/2}_h\cdot\nabla\psi^n_h, e^{n+1/2}_h)_{\Omega^{n+1/2}_h},\\
	A_2^{n}:=&\frac{1}{\tau}\left[\left(\psi_h^{n+1}-\psi_h^n, e_h^{n+1}\right)_{\Omega_h^{n+1/2}}-\left(\psi_h^{n+1}-\psi_h^n, e_h^{n+1}\right)_{\Omega_h^{n+1}}\right],\\
	A_3^{n}:=&-\frac12 \tau (d_t \psi_h^{n+1}, d_te^{n+1}_h)_{\Omega^{n+1/2}_h}.
\end{align*}
Similarly to the proof of Lemma \ref{truncation-error-lemma} with integration by parts, we obtain the estimate for $ \mathscr A_{6}^{n} $:
\begin{equation}\label{scrA1-estimate}
				|\mathscr{A}_6^n|\leq C\tau \|e^{n+1/2}_h\|_{\Omega_h^{n+1/2}}\|\psi^n_h\|_{\Omega^n_h}+ C\left(\tau\|e_h^{n+1}\|_{\Omega_h^{n+1}}+\|e^{n+1}_h-e^n_h\|_{\Omega_h^{n+1}}\right)\|d_t\psi^{n+1}_h\|_{\Omega^{n+1}_h}.
\end{equation}
Testing \eqref{main-scheme-error-equation-1} with $v_h=\psi_h^n$ and combining it with \eqref{duality-form-2}, we have 
\begin{align}\label{expression-I1-I2}
	\mathcal{I}_1+\mathcal{I}_2&=\frac{1}{\tau}\left[(\psi^{n+1}_h, e^{n+1}_h)_{\Omega^{n+1}_h}-(\psi^{n}_h, e^{n}_h)_{\Omega^{n}_h}\right]+\mathscr{A}_6^{n}-\langle\mathscr{E}^{n+1}_1, \psi_h^n\rangle.
\end{align}
Regarding $\mathcal{I}_3$, when $n\geq 1$, we can test \eqref{main-scheme-error-equation-2} with $q_h=\phi_h^n$ at steps $n+1$ and $n$ to obtain the following expression:
\begin{equation}
	\mathcal{I}_3
	=\frac{1}{2}\beta\tau\left[\left(\nabla(\eta_h^{n+1}-\eta_h^n), \nabla \varphi_h^n\right)_{\Omega_h^{n+1}}+\left(\nabla(\eta_h^{n}-\eta_h^{n-1}), \nabla \varphi_h^n\right)_{\Omega_h^{n}}\right]+\mathscr{A}_7^{n}-\frac{1}{2}\langle\mathscr{E}_2^{n+1}+\mathscr{E}_2^n, \varphi_h^n\rangle,\label{expression-I3}
\end{equation}
where $ \mathscr A_{7}^{n} $ is given by 
\begin{align}
	\mathscr{A}_7^{n}:= &\,\frac{1}{2}\left({\rm div} e_h^{n+1}, \varphi_h^n\right)_{\Omega_h^{n+1}}+\frac{1}{2}({\rm div} e_h^{n}, \varphi_h^n)_{\Omega_h^{n}}-({\rm div}e_h^{n+1/2}, \varphi_h^n)_{\Omega_h^{n+1/2}}\notag\\
	\leq &\, C\tau\|e_h^{n+1}-e_h^n\|_{\Omega_h^{n+1}}\|\nabla \varphi_h^n\|_{\Omega_h^n}.\label{scraA2-estimate}
\end{align}

Summing up \eqref{formula-g-h-n} from $n=1$ to $n=M$ and combining the expressions \eqref{expression-I1-I2} and \eqref{expression-I3}, taking into account $\psi_h^{M}=0$ and the estimates for $\mathscr{A}_6^{n}$ and $ \mathscr A_{7}^{n} $ (given by \eqref{scrA1-estimate}, \eqref{scraA2-estimate}), and considering the estimates \eqref{E1-E3-estimate}, \eqref{F-G-estimate} and \eqref{improve-E-1-estimate}, we obtain
\begin{align}
	&\Big|\tau\sum_{n=1}^{M-1}(g^n_h, e^{n+1/2}_h)_{\Omega^{n+1/2}_h}\Big|\notag\\
	\leq& |(\psi_h^1, e_h^1)_{\Omega_h^1}|+\Big(Ch^{r+1}+C\tau^2+C\max_{1\leq n\leq M}\|e_h^{n+1}-e_h^n\|_{\Omega_h^{n+1}}+C\tau\max_{1\leq n\leq M}\|e_h^n\|_{\Omega_h^n}\notag\\
	    &+C\tau\max_{0\leq n\leq M}\|\nabla(\eta_h^{n+1}-\eta_h^n)\|_{\Omega_h^{n+1}}\Big)\cdot \sum_{n=1}^{M-1}\tau\Big(\|\nabla\psi_h^n\|_{\Omega_h^n}+\|\nabla\varphi^n_h\|_{\Omega_h^n}+\|d_t\psi_h^{n+1}\|_{\Omega_h^{n+1}}\Big).\label{section4-temp}
\end{align}
We now establish the $\ell^1L^2$-estimate for $(\psi_h^n,\varphi_h^n)$ and utilize the results from Lemma \ref{energy-error-estimate-1} and Lemma \ref{energy-error-estimate-2} to obtain the main result. The $\ell^1L^2$-estimate is given by:
\begin{equation}
	\|\psi_h^1\|_{\Omega_h^1}+\sum_{n=1}^{M-1}\tau\left(\|\nabla\psi_h^n\|_{\Omega_h^n}+\|\nabla\varphi_h^{n}\|_{\Omega_h^n}+\|d_t\psi_h^{n+1}\|_{\Omega_h^{n+1}}\right)\leq C\min\{\ell_\tau,\ell_h\}\sum_{n=1}^{M-1}\tau\|g^n_h\|_{\Omega_h^{n+1/2}},\label{ell1L2-form1}
\end{equation}
where we have used the notations $\ell_\tau:=\ln\left(\frac{1}{\tau}+1\right)$ and $\ell_h:=\ln\left(\frac{1}{h}+1\right)$.
A detailed proof of \eqref{ell1L2-form1} is provided in Appendix \ref{ell1L2-estimate-section}. By substituting the results from Lemma \ref{energy-error-estimate-1} and Lemma \ref{energy-error-estimate-2} into \eqref{section4-temp}, we obtain:
\begin{align}
	\Big|\tau\sum_{n=1}^{M-1}(g^n_h, e^{n+1/2}_h)_{\Omega^{n+1/2}_h}\Big|\leq C\min\{\ell_\tau,\ell_h\}\left(h^{r+1}+\tau^2\right)\Big(\sum_{n=1}^{M-1}\tau \|g^n_h\|_{\Omega_h^{n+1/2}}\Big).
\end{align}
Since $g^n_h$ can be chosen arbitrarily, utilizing the duality between $\ell^\infty$ and $\ell^1$, we can deduce: 
\begin{align}
	\max_{1\leq n\leq M-1}\|e_h^{n+1/2}\|_{\Omega_h^{n+1/2}}\leq C\min\{\ell_\tau,\ell_h\}(h^{r+1}+\tau^2 ).
\end{align}
Considering that $\|e_h^{n+1}-e_h^n\|_{\Omega_h^{n+1}}\leq Ch^{r+1}+C\tau^2$ (as shown in Lemma \ref{energy-error-estimate-2}), we have established the main result in \eqref{main-theorem}. 
\qed

\section{Numerical experiments}%
\label{sec:Numerical experiments}
In this section, we provide numerical tests to support the theoretical result proved in Theorem \ref{main-theorem} and to illustrate the effectiveness of the projection method in simulating Navier--Stokes flow in a time-varying domain which contains a rotating propeller.
\begin{example}[The Stokes equations]\label{numerical-solution-stokes}
Consider the following Stokes equation in a time-varying domain $ \Omega (t)\in \R^{2} $.
\begin{align*}\left\{\begin{array}{ll}
	\partial_t u - \Delta u + \nabla p = f & \; \; \mbox{in}\; \;\bigcup\limits_{t\in(0,T]}\Omega(t)\times \{t\},\\
\nabla \cdot u = 0 & \; \; \mbox{in}\; \;\bigcup\limits_{t\in(0,T]} \Omega(t)\times\{t\},\\
u = 0 & \; \; \mbox{on} \; \;\bigcup\limits_{t\in(0,T]}\partial\Omega(t)\times \{t\},\\
u = u_0 & \; \; \mbox{on}\; \;\Omega(0),
\end{array}\right.
\end{align*}
where $ \Omega (t) $ is a dumbbell-shape domain given by:
\[
    \Omega (t) = \{(x,y):e^{-t/8}\phi (x,y) - 1 \le 0\}\;\text{for}\;\phi (x,y) = \frac{y^{2}}{(0.7x^{2}+0.3)^{2}} + x^{2}.
\]
The exact solutions to the problem are 
$$
u (x,y;t) = F (x,y;t)(c_{1}(x,y),c_{2}(x,y))^{\top} ,\;p (x,y;t) = x+y ,
$$ 
where
\begin{align*}
&c_{1}(x,y) = y (0.7x^{2}+0.3)^{-2},\quad c_{2}(x,y) = -x + 1.4xy^{2}(0.7x^{2}+0.3)^{-3},\\
&F (x,y;t) = e^{-t/8}\phi (x,y) - 1.
\end{align*}
To ensure consistency, the source term $ f $ is chosen adaptively to the exact solutions, specifically as $ f = \partial_{t} u - \Delta u + \nabla p $. Additionally, we select the velocity function $ w $ to be
$
    w (x,y;t) = -(\partial_{t}F \nabla F)/(\nabla F^{\top}\nabla F).
$

The initial and final discretized domains, denoted as $ \Omega_{h} ( 0 ) $ and $ \Omega_{h} ( 1 ) $ respectively, are illustrated in Figure \ref{fig-evolving-1}. These domains are obtained by employing the $ P_1 $ element and $ P_2 $ element, representing the piecewise linear and quadratic finite elements, respectively. 
\begin{figure}
	\centering
	\subfigure[$ P1 $ element at $ t = 0 $]{\includegraphics[trim = .1cm .1cm .1cm .1cm, clip=true,width=0.48\textwidth,height=0.24\textwidth]{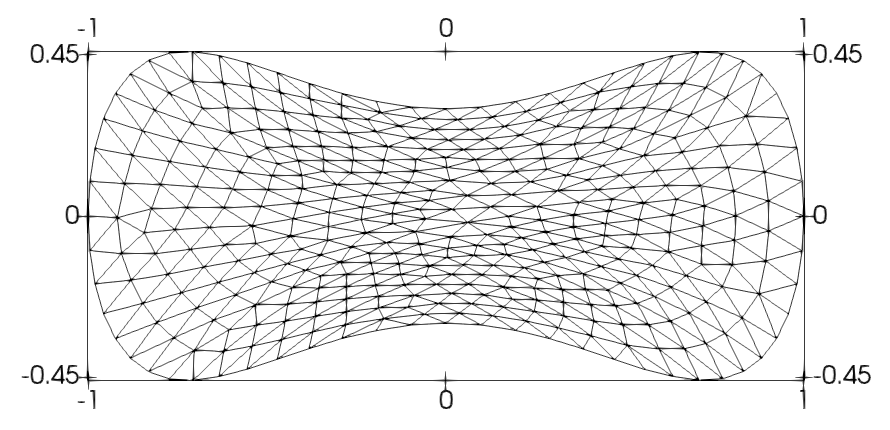}}
	\subfigure[$ P1 $ element at $ t = 1 $]{\includegraphics[trim = .1cm .1cm .1cm .1cm, clip=true,width=0.48\textwidth,height=0.24\textwidth]{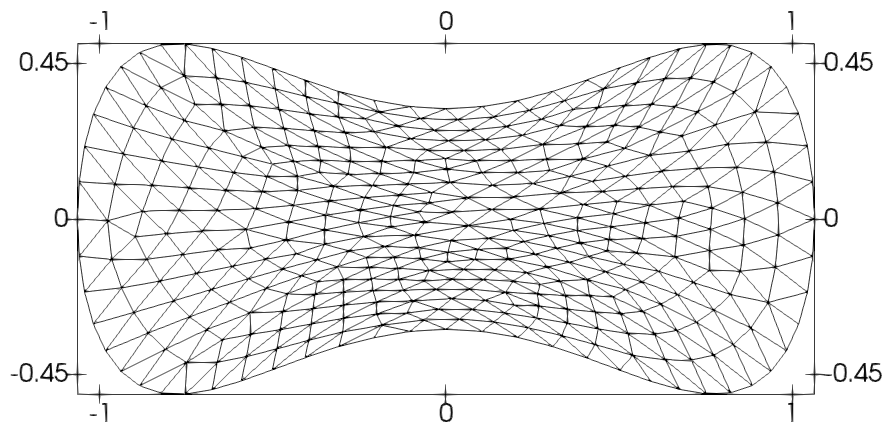}}
	\subfigure[$ P2 $ element at $ t = 0 $]{\includegraphics[trim = .1cm .1cm .1cm .1cm, clip=true,width=0.48\textwidth,height=0.24\textwidth]{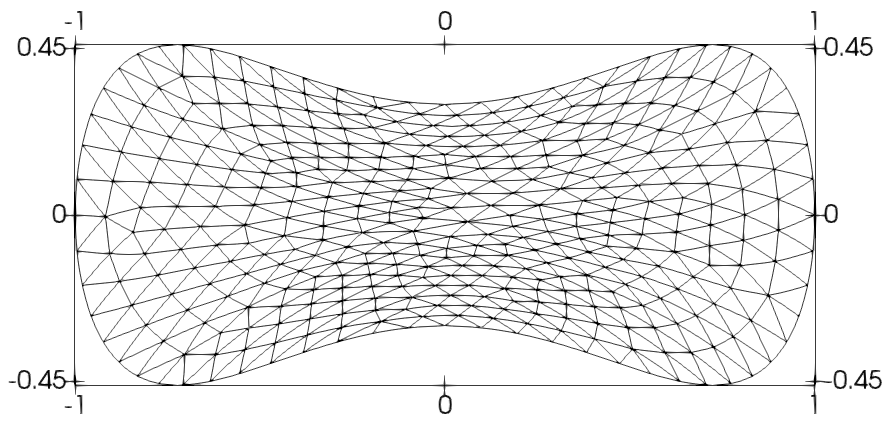}}
	\subfigure[$ P2 $ element at $ t = 1 $]{\includegraphics[trim = .1cm .1cm .1cm .1cm, clip=true,width=0.48\textwidth,height=0.24\textwidth]{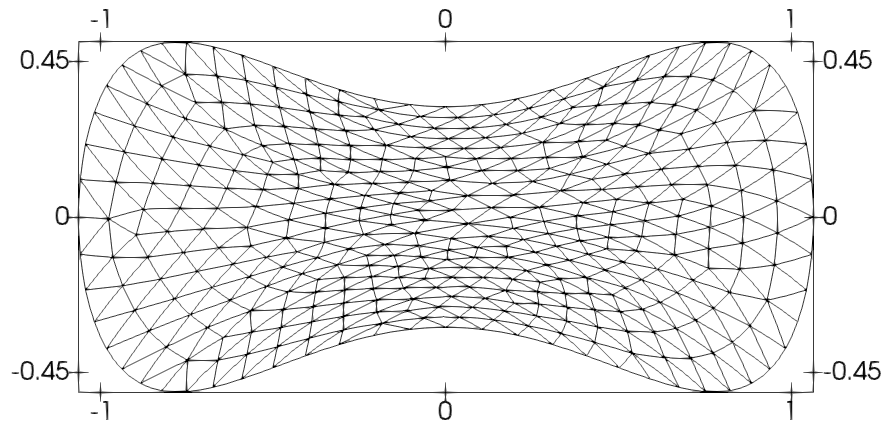}}
	\caption{Meshes of $ P_{1} $ and $ P_{2} $ elements at time $ T = 0 $ and $ T = 1 $.}
\label{fig-evolving-1}
\end{figure}

To assess the convergence properties of the numerical scheme, we conducted a convergence test at $ T = 1 $ using the $ P_2-P_1 $ element and a suitably small mesh size that ensures negligible errors from the space discretization. The resulting errors of the numerical solutions are depicted in Figure \ref{fig-1} (a) for different time step sizes: $ \tau = 1/32, 1/64, 1/128, 1/256 $. The observed errors demonstrate second-order convergence in time, which aligns with the theoretical findings established in Theorem \ref{main-theorem}.

In addition to investigating the convergence in time, we also conducted a convergence test to assess the spatial discretization. For this purpose, we employed three different sets of finite elements: $ P_{1b}-P_1 $, $ P_2-P_1 $, and $ P_3-P_2 $, while keeping the time step sizes sufficiently small to ensure minimal errors from the time discretization. The errors of the numerical solutions are presented in Figure \ref{fig-1} (b) for varying mesh sizes: $ h = 1/16,1/24,1/36,1/54 $. The results demonstrate that the numerical solutions exhibit $ r+1 $-th order convergence in space, where $ r $ corresponds to the order of the FEM. This finding aligns with the theoretical results established in Theorem \ref{main-theorem} for $ r = 2,3 $. Notably, for the $ P_{1b}-P_1 $ element, we verified that the inf-sup condition \eqref{inf-sup-condition} is satisfied. Therefore, we can attain second-order convergence using the same approach presented in this article.
\begin{figure}[htp!]
	\centering
	\subfigure[$ L^{2} $ error of $ u $ from time discretization]{\includegraphics[trim = .1cm .1cm .1cm .1cm, clip=true,width=0.48\textwidth,height=0.32\textwidth]{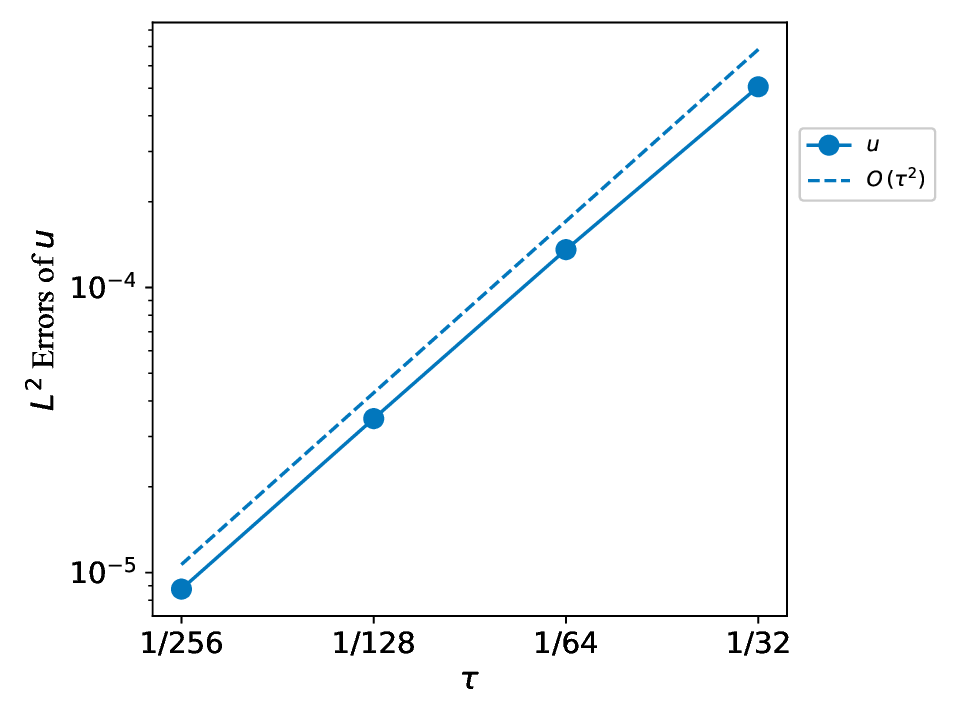}}
	\subfigure[$ L^{2} $ error of $ u $ from space discretization]{\includegraphics[trim = .1cm .1cm .1cm .1cm, clip=true,width=0.48\textwidth,height=0.32\textwidth]{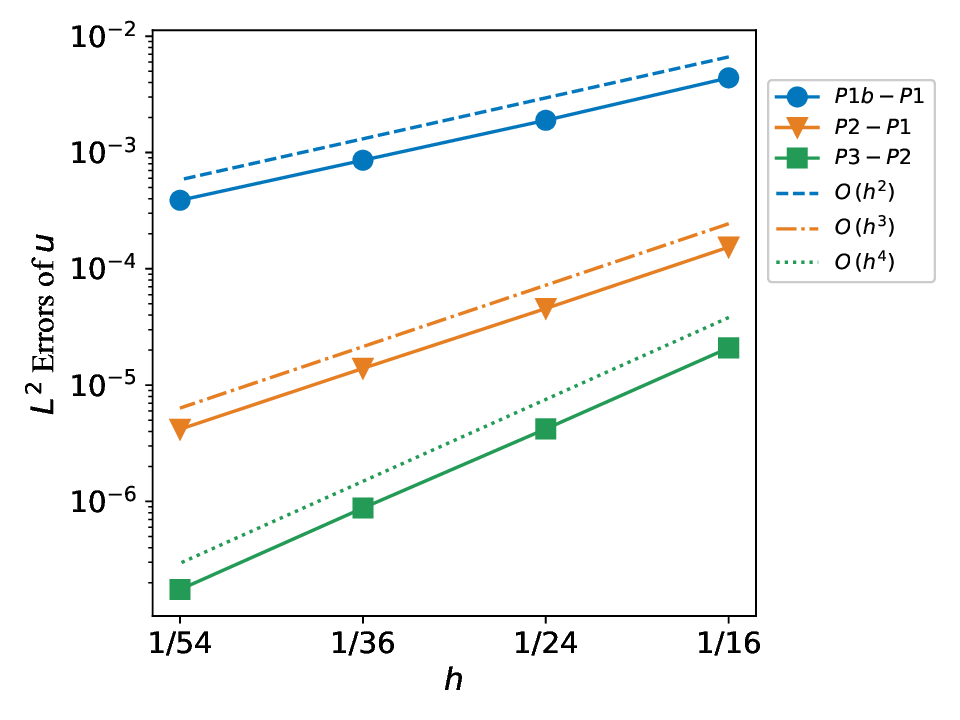}}
	\caption{Errors from time and space discretization at time $ T = 1 $.}
\label{fig-1}
\end{figure}

\end{example}

\begin{example}[The Navier--Stokes equations]
	In this example, we investigate the convergence order of numerical solutions for the Navier--Stokes equations. The governing equation for the evolving domain $ \Omega (t)\in \R^{2} $ is given by:
\begin{align}\label{NS-numerical-ex}
\left\{\begin{array}{ll}
\partial_t u + u \cdot \nabla u - \Delta u + \nabla p = f & \mbox{in}\,\,\bigcup\limits_{t\in(0,T]}\Omega(t)\times\{t\},\\
\nabla \cdot u = 0 & \mbox{in}\,\,\bigcup\limits_{t\in(0,T]}\Omega(t)\times\{t\},\\
u = 0 & \mbox{on} \,\,\bigcup\limits_{t\in(0,T]}\partial\Omega(t)\times \{t\},\\
u = u_0 & \mbox{on}\,\,\Omega(0),
\end{array}\right.
\end{align}
For this numerical test, we choose the evolving domain $ \Omega (t) $ and the exact solutions $ u $ and $ p $ to be the same as in Example \ref{numerical-solution-stokes}. To ensure accuracy, the source term $ f $ is adaptively chosen to match the first equation in \eqref{NS-numerical-ex}. It is important to note that \eqref{NS-numerical-ex} includes an additional nonlinear term $ u\cdot \nabla u $, which necessitates an adjustment in the numerical method. In the fully discrete scheme \eqref{main-scheme}, an extra term
$$
\Big((\frac{3}{2} u_{h}^{n} - \frac{1}{2} u_{h}^{n-1})\cdot \nabla u_{h}^{n+1/2},v_{h}\Big)_{\Omega_{h}^{n+1/2}} 
$$ 
should be added to the left-hand side of the equation.

Similarly to Example \ref{numerical-solution-stokes}, we assess the convergence behavior of the numerical solutions for the Navier--Stokes equations. Specifically, we examine the convergence in time at $ T = 1 $ using the $ P_2-P_1 $ element, with a sufficiently small mesh size that ensures the errors from space discretization are negligible. The errors of the numerical solutions are presented in Figure \ref{fig-2} (a) for various time step sizes: $ \tau = 1/32, 1/64, 1/128, 1/256 $. The numerical results demonstrate a second-order convergence in time, which aligns with the convergence order observed in the Stokes equation.

In addition, we investigate the convergence of spatial discretization using the $ P_{1b}-P_1 $, $ P_2-P_1 $, and $ P_3-P_2 $ elements, considering sufficiently small time step sizes that ensure the errors from time discretization are negligible. Figure \ref{fig-2} (b) illustrates the errors of the numerical solutions for different mesh sizes: $ h = 1/16,1/24,1/36,1/54 $. The results indicate that the numerical solutions exhibit $ r+1 $-th order convergence in space for $ r $-th order FEMs. This convergence behavior aligns with the observed convergence order in the Stokes equation.
\begin{figure}[htp!]
	\centering
	\subfigure[$ L^{2} $ error of $ u $ from time discretization]{\includegraphics[trim = .1cm .1cm .1cm .1cm, clip=true,width=0.48\textwidth,height=0.32\textwidth]{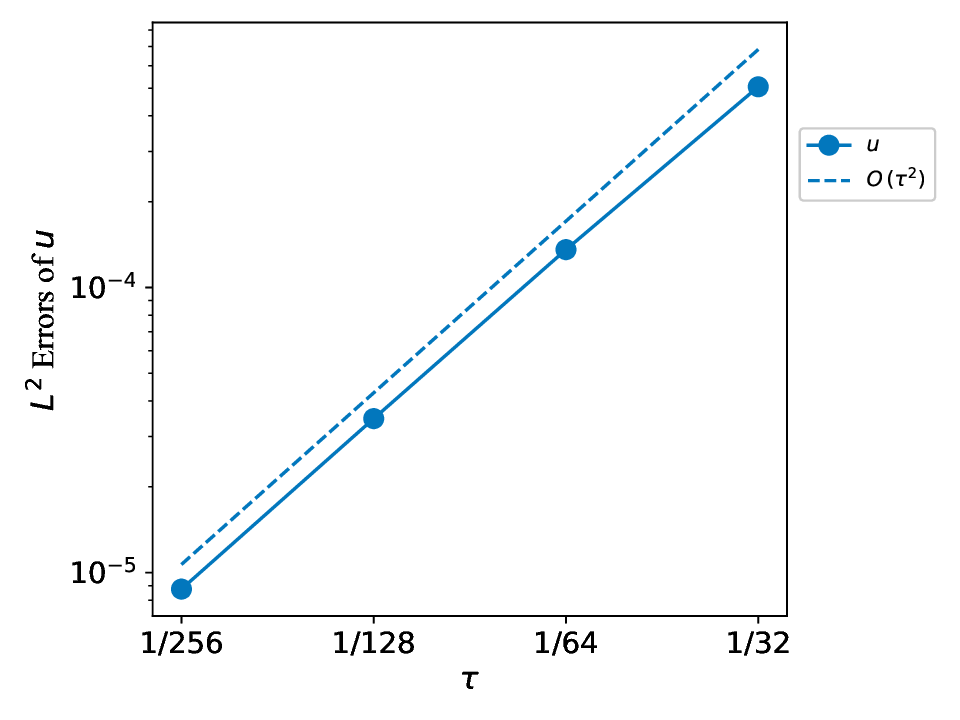}}
	\subfigure[$ L^{2} $ error of $ u $ from space discretization]{\includegraphics[trim = .1cm .1cm .1cm .1cm, clip=true,width=0.48\textwidth,height=0.32\textwidth]{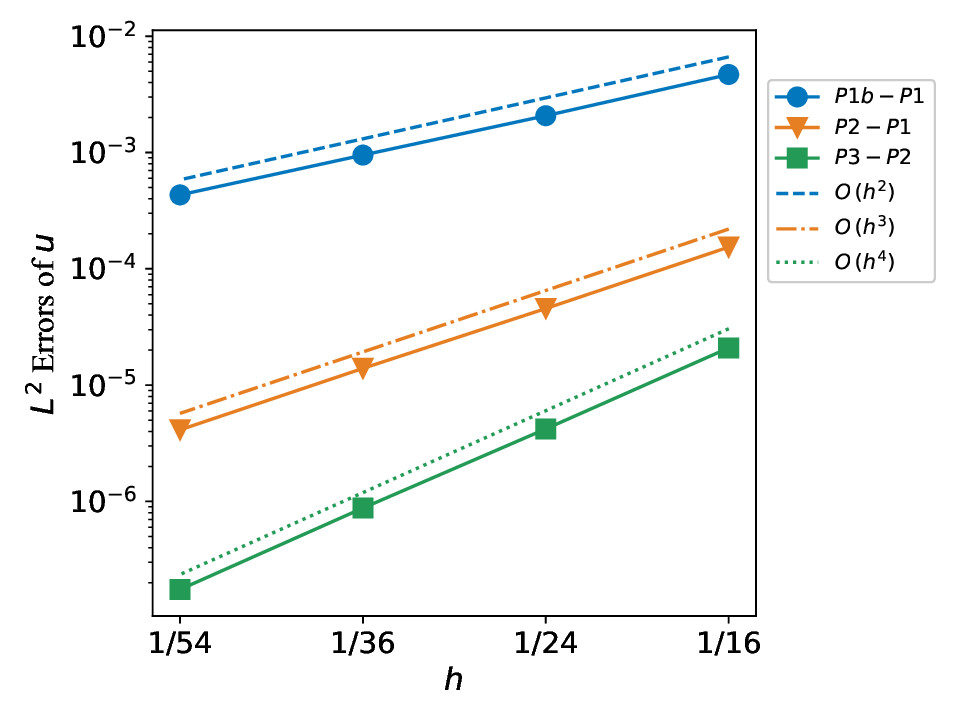}}
	\caption{Errors from time and space discretization at time $ T = 1 $.}
\label{fig-2}
\end{figure}
\end{example}

\begin{example}[Navier--Stokes flow in a domain with rotating propeller]
	In this example, we investigate the fluid motion surrounding a rotating propeller, governed by the Navier--Stokes equation with slip boundary conditions, i.e. 
	\begin{equation}\label{slip-NS-equation}
    \left\{\begin{array}{ll}
		    \partial_{t}u + u\cdot \nabla u - \nabla\cdot (2 \mu\mathbb{D}u - \frac{1}{\rho}p I) = 0 & \text{in}\;\bigcup\limits_{t \in (0,T]} \Omega (t)\times \{t\},\\
	\nabla \cdot u = 0,&\text{in}\;\bigcup\limits_{t \in (0,T]}\Omega (t)\times \{t\},\\
	u\cdot \textbf{n} = w\cdot \textbf{n} & \text{on}\;\bigcup\limits_{t \in (0,T]}\partial \Omega (t)\times \{t\},\\
	((2 \mu \mathbb{D} u - \frac{1}{\rho} p I) \cdot n)_{\rm tan} + k u_{\rm tan} = 0 & \text{on}\;\bigcup\limits_{t \in (0,T]} \partial \Omega (t)\times \{t\},\\
	u = u_{0} & \text{on}\; \Omega (0),
    \end{array}\right.
\end{equation}
where the subindex ${\rm tan}$ stands for the tangential component of a vector, $ w $ is the velocity of the propeller defined on the boundary and has a natural extension to the whole domain $ \Omega (t) $, $ \mu = 0.001 $ denotes the viscosity, $ \rho = 1000 $ is the fluid density, and $ \textbf{n} $ is the outward normal vector on the boundary. The initial domain $ \Omega ( 0 ) $ corresponds to a unit sphere with an ellipse removed, defined as:
\begin{align*}
    \Omega ( 0 ) = \{ ( x,y ):x^{2}+y^{2} \le 1, \;\text{and}\; 2 x^{2} + 4y^{2} \ge 1 \}.
\end{align*} 
The propeller, depicted as the middle ellipse in Figure \ref{Fig-rotation}, has a prescribed velocity profile. Specifically, when $0 \leq t \leq 2$, the propeller velocity is given by $w(x,y;t) := (-2ty, 2tx)$, and for $t > 2$, it is defined as $w(x,y;t) := (-2y, 2x)$. This velocity naturally extends to domain $\Omega$.

The numerical method studied in this paper can be extended to the Navier--Stokes equations with slip boundary conditions. The scheme needs to be modified to suit the slip boundary conditions, and optimal-order convergence in space and second-order convergence in time can be established similarly.

We perform convergence tests for the accuracy of the numerical scheme. To investigate the convergence in time, we select the $P_2-P_1$ element with a sufficiently small mesh size, ensuring that the errors from space discretization are negligible. The results, presented in Figure \ref{figure-slip-bc} (a), demonstrate the errors of the numerical solutions for various time step sizes: $\tau = 1/960, 1/1440, 1/2160, 1/3240$, and indicating that the numerical solutions exhibit second-order convergence in time. This is consistent with the theoretical results  in Theorem \ref{main-theorem}.

\begin{figure}[htp!]
	\centering
	\subfigure[$ L^{2} $ error of $ u $ from time discretization]{\includegraphics[trim = .1cm .1cm .1cm .1cm, clip=true,width=0.48\textwidth,height=0.32\textwidth]{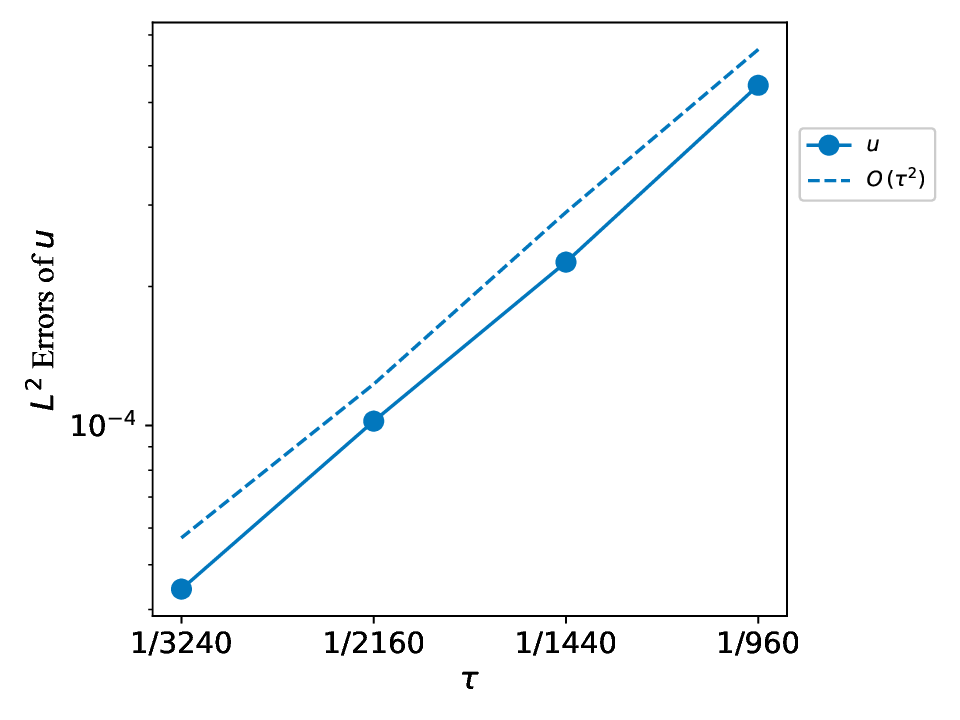}}
	\subfigure[$ L^{2} $ error of $ u $ from space discretization]{\includegraphics[trim = .1cm .1cm .1cm .1cm, clip=true,width=0.48\textwidth,height=0.32\textwidth]{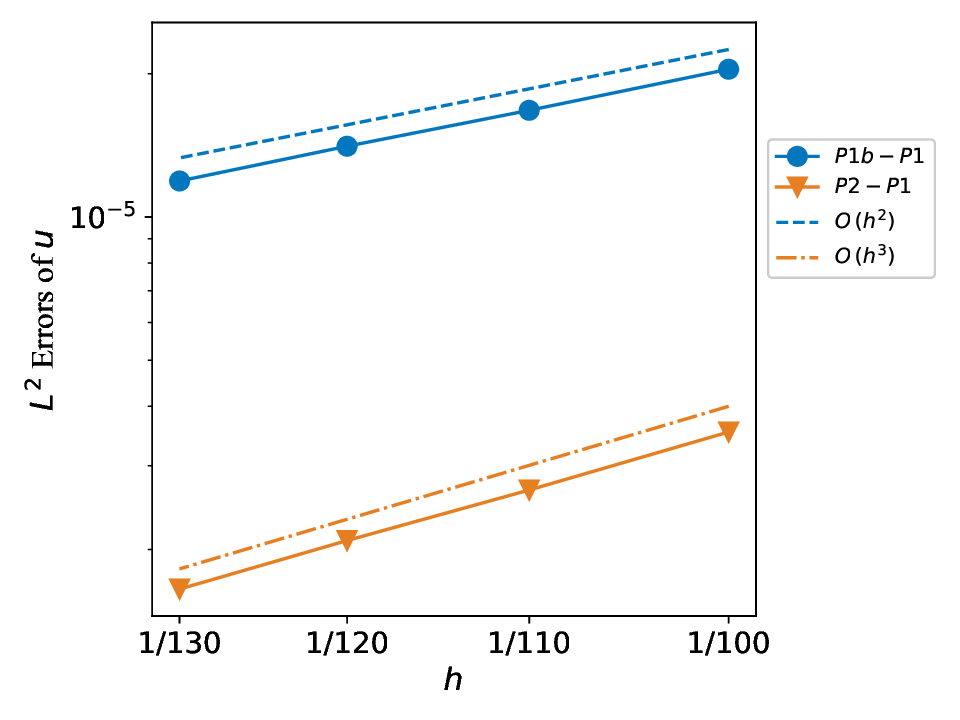}}
	\caption{Errors from time and space discretization at time $ T = 1 $.}
	\label{figure-slip-bc}
\end{figure}

In addition to time discretization, we illustrate the convergence of spatial discretization using both the $P_{1b}-P_1$ and $P_2-P_1$ elements, with sufficiently small time step sizes to ensure negligible errors from time discretization. Figure \ref{figure-slip-bc} (b) presents the errors of the numerical solutions for various mesh sizes: $h = 1/100, 1/110, 1/120, 1/130$. The results demonstrate that the numerical solutions exhibit $(r+1)$th-order convergence in space for finite elements of degree $r$. This aligns with the theoretical results presented in Theorem \ref{main-theorem} in the case $r = 2$. For the $P_{1b}-P_1$ element, since the inf-sup condition \eqref{inf-sup-condition} is satisfied, the same approach as described in this article yields second-order convergence.

To illustrate the propeller rotation, we conduct simulations with a mesh size of $h = 0.01$ and a time step size of $\tau = 0.001$. Figure \ref{Fig-rotation} depicts the process of the propeller rotation and displays the magnitude of the velocity field $|u|$. The figure portrays the flow of the fluid driven by the yellow elliptic propeller, offering insights into propeller-driven flows.

\begin{figure}[htp!]
	\centering
	\subfigure[$ t = 0 $]{\includegraphics[trim = .1cm .1cm .1cm .1cm, clip=true,width=0.24\textwidth,height=0.24\textwidth]{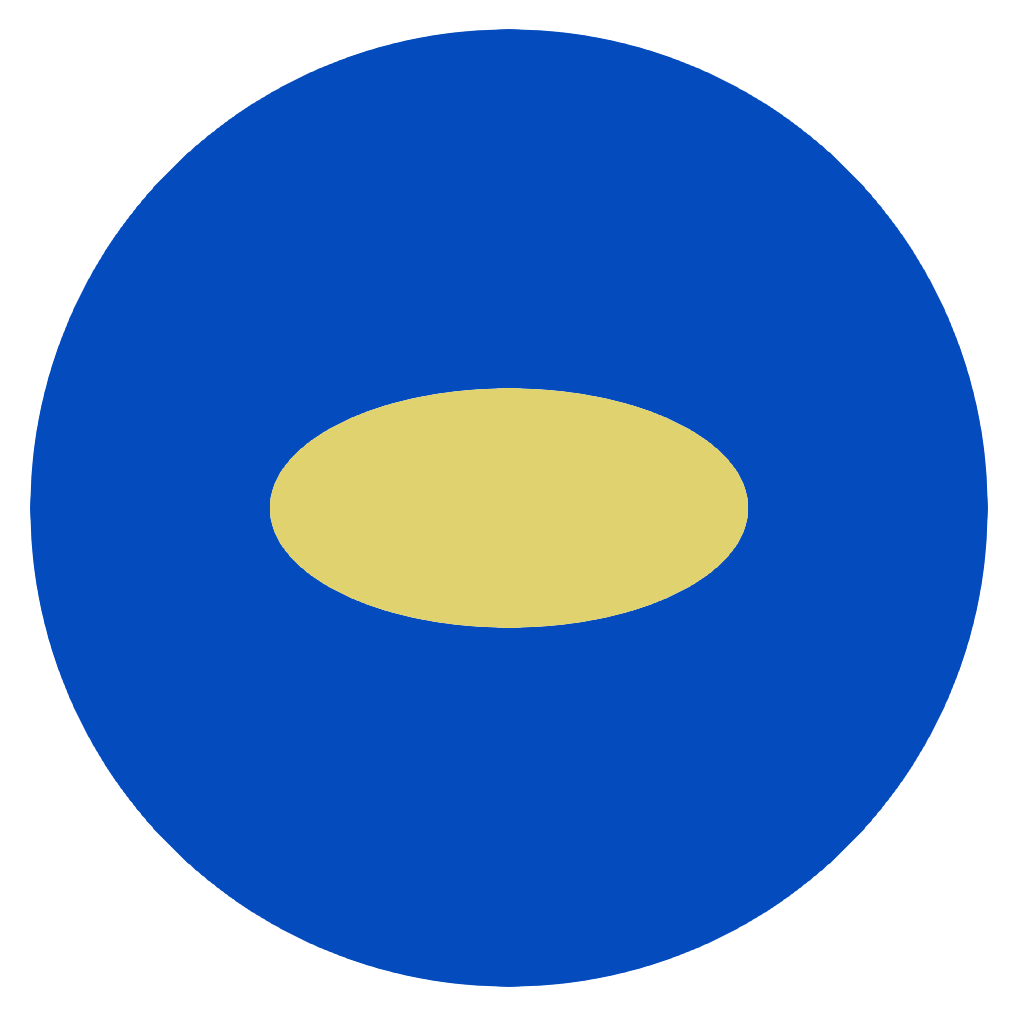}}
	\subfigure[$ t = 0.5 $]{\includegraphics[trim = .1cm .1cm .1cm .1cm, clip=true,width=0.24\textwidth,height=0.24\textwidth]{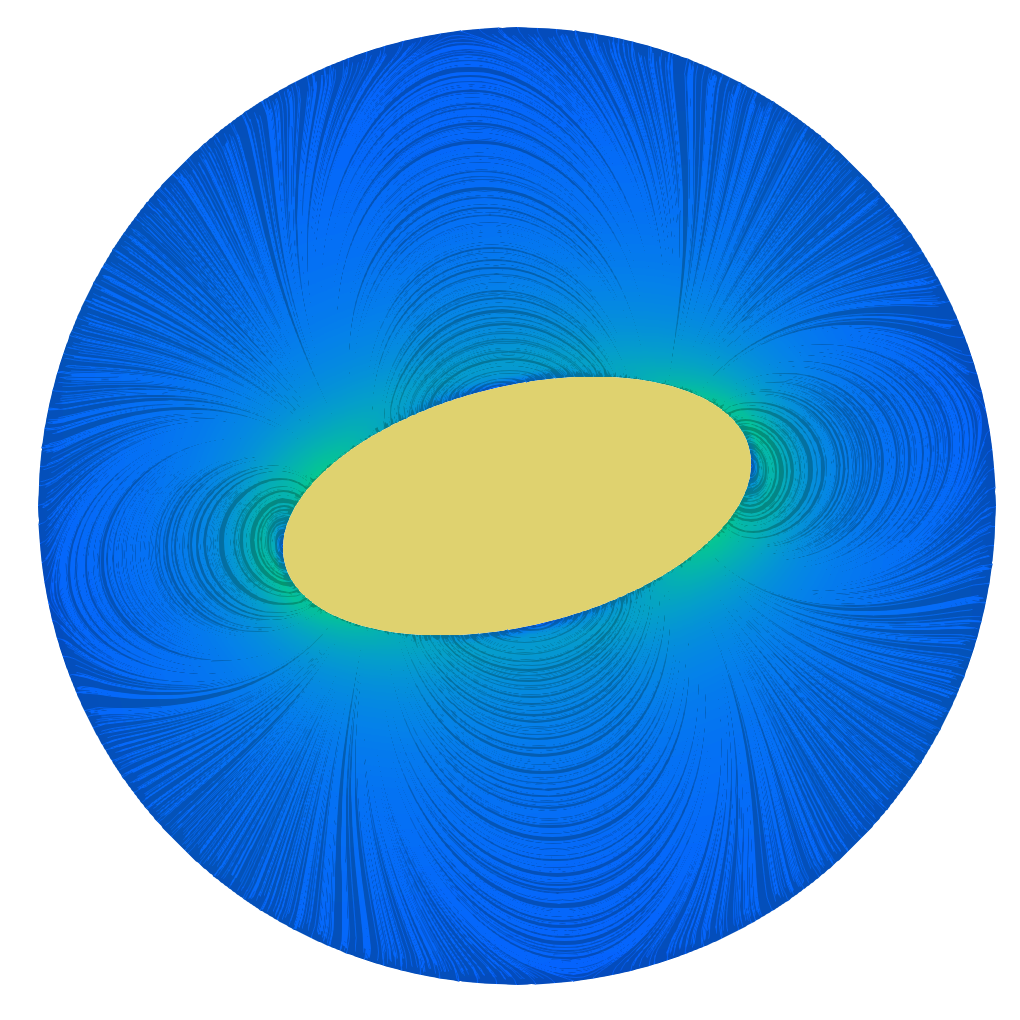}}
	\subfigure[$ t = 1 $]{\includegraphics[trim = .1cm .1cm .1cm .1cm, clip=true,width=0.24\textwidth,height=0.24\textwidth]{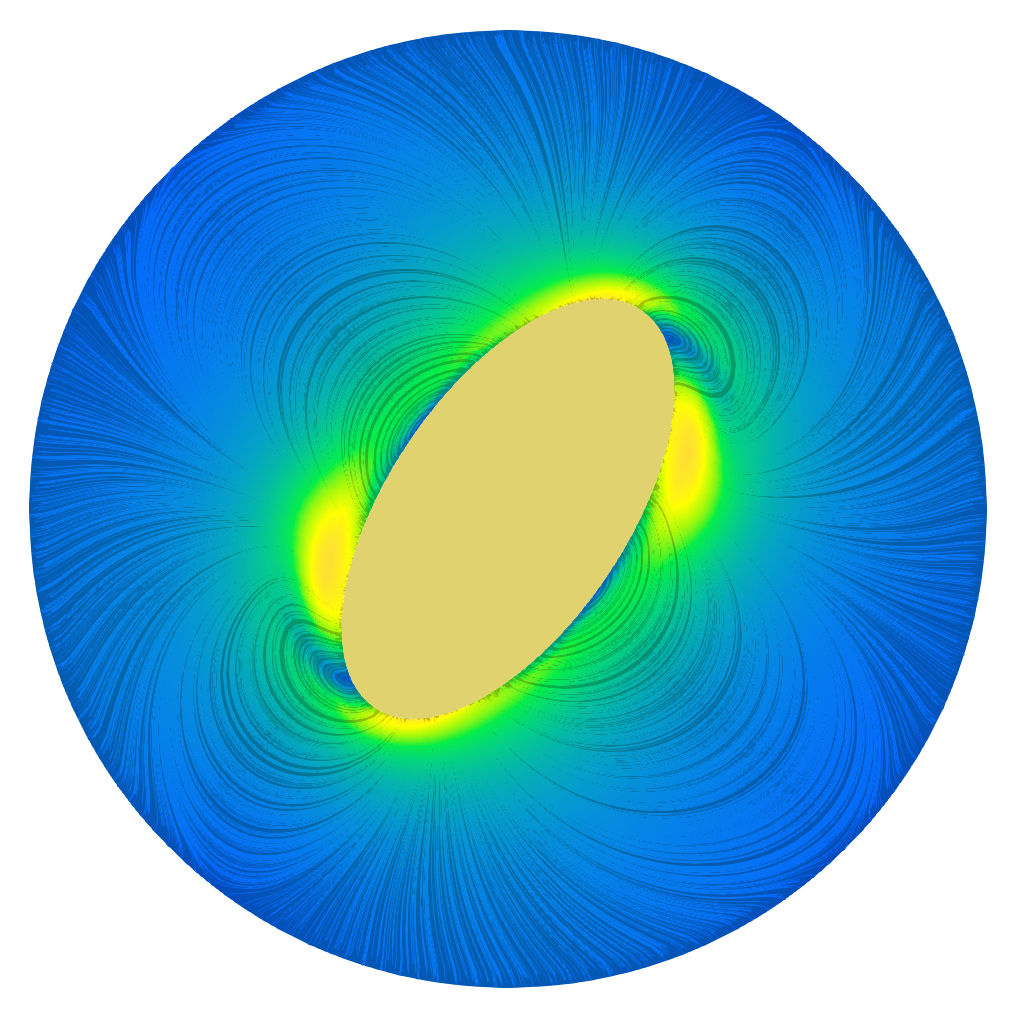}}
	\subfigure[$ t = 1.5 $]{\includegraphics[trim = .1cm .1cm .1cm .1cm, clip=true,width=0.24\textwidth,height=0.24\textwidth]{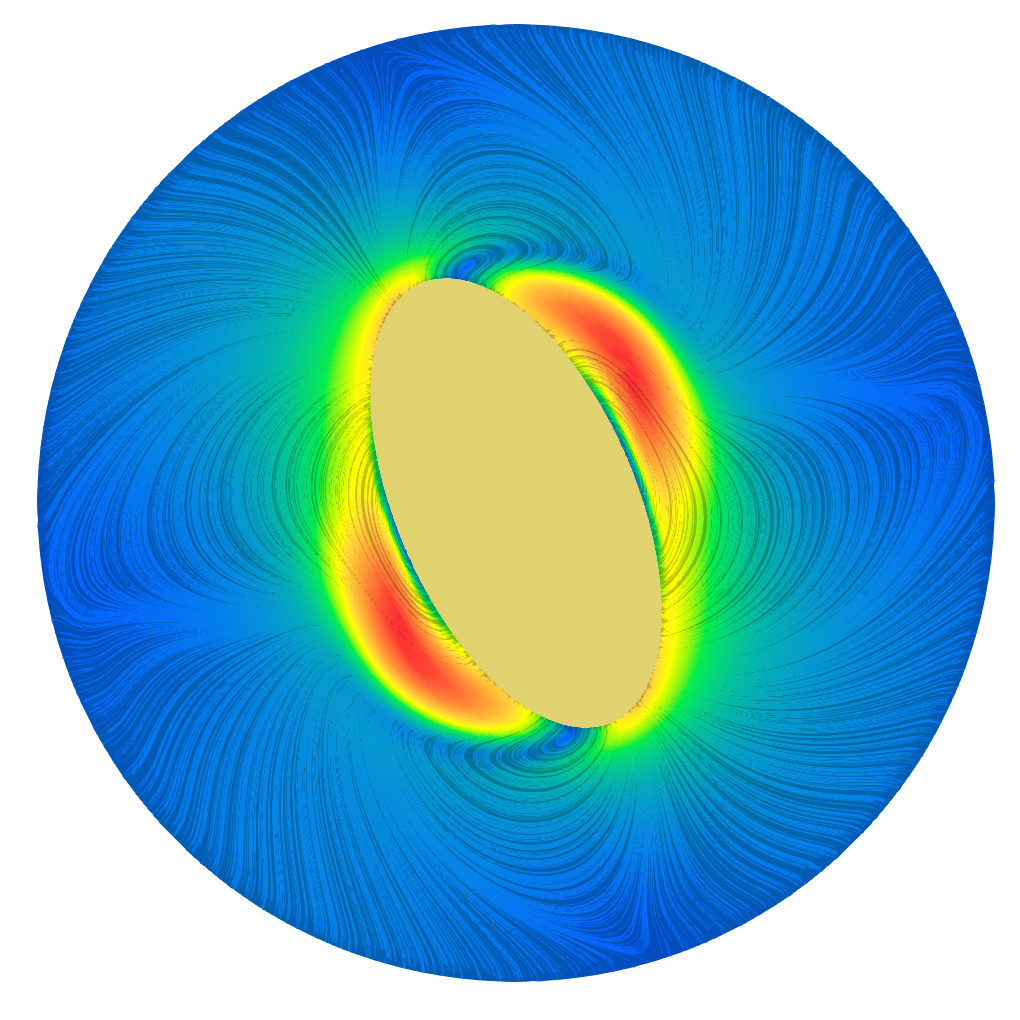}}
	\subfigure[$ t = 2 $]{\includegraphics[trim = .1cm .1cm .1cm .1cm, clip=true,width=0.24\textwidth,height=0.24\textwidth]{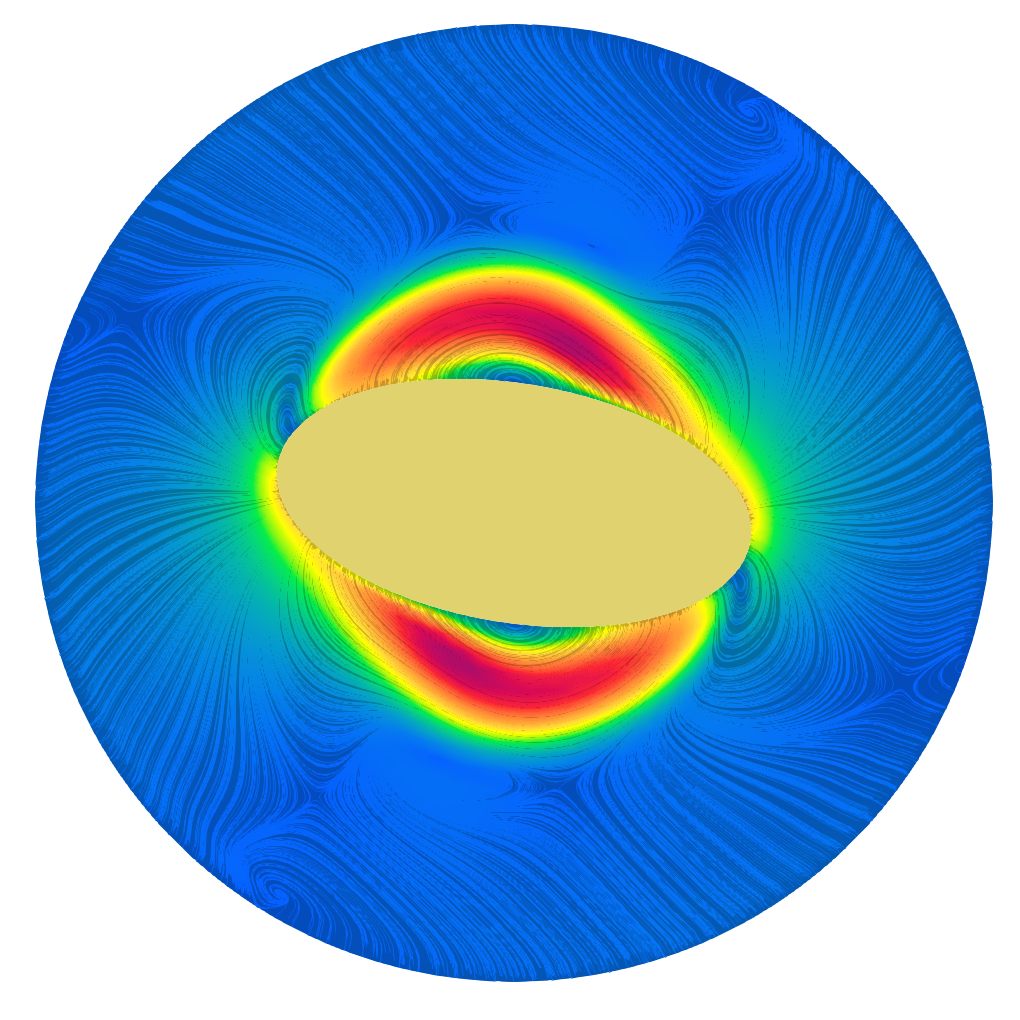}}
	\subfigure[$ t = 2.5 $]{\includegraphics[trim = .1cm .1cm .1cm .1cm, clip=true,width=0.24\textwidth,height=0.24\textwidth]{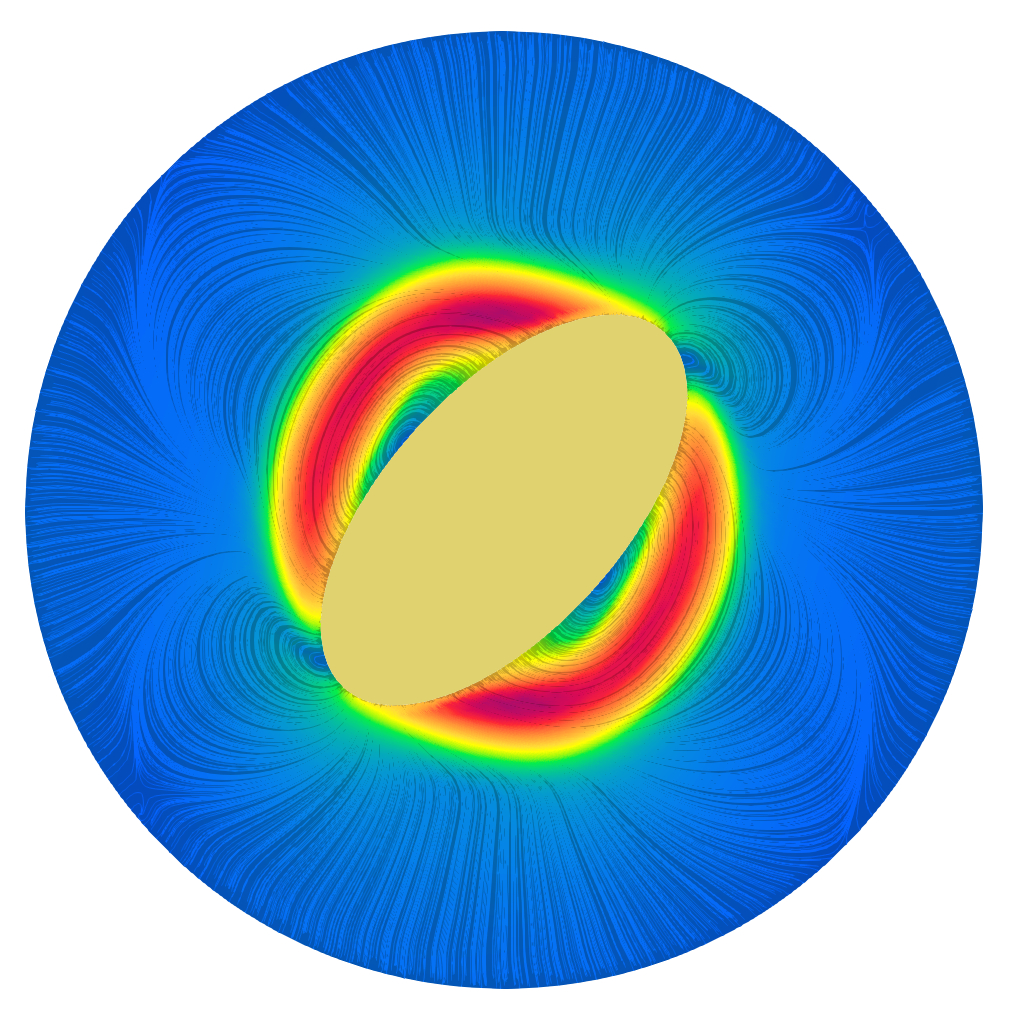}}
	\subfigure[$ t = 3 $]{\includegraphics[trim = .1cm .1cm .1cm .1cm, clip=true,width=0.24\textwidth,height=0.24\textwidth]{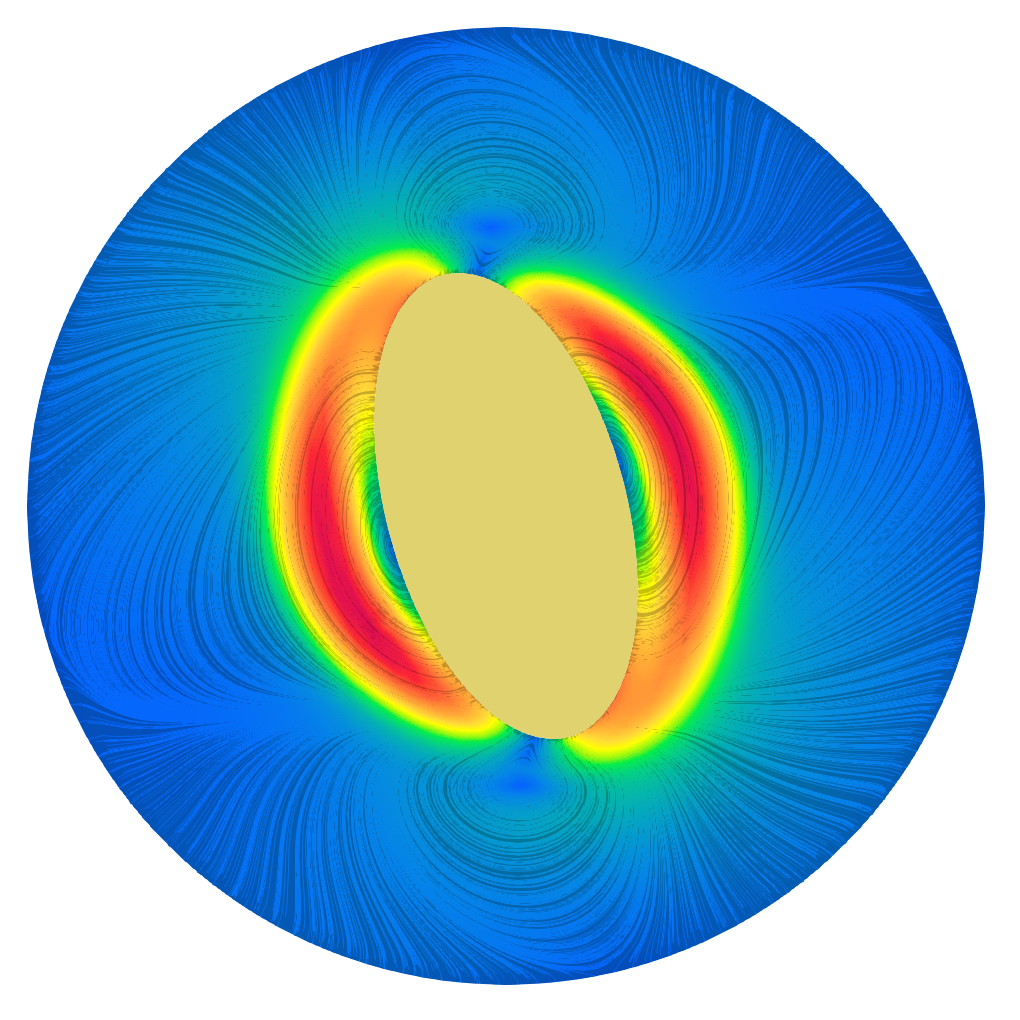}}
	\subfigure[$ t = 3.5 $]{\includegraphics[trim = .1cm .1cm .1cm .1cm, clip=true,width=0.24\textwidth,height=0.24\textwidth]{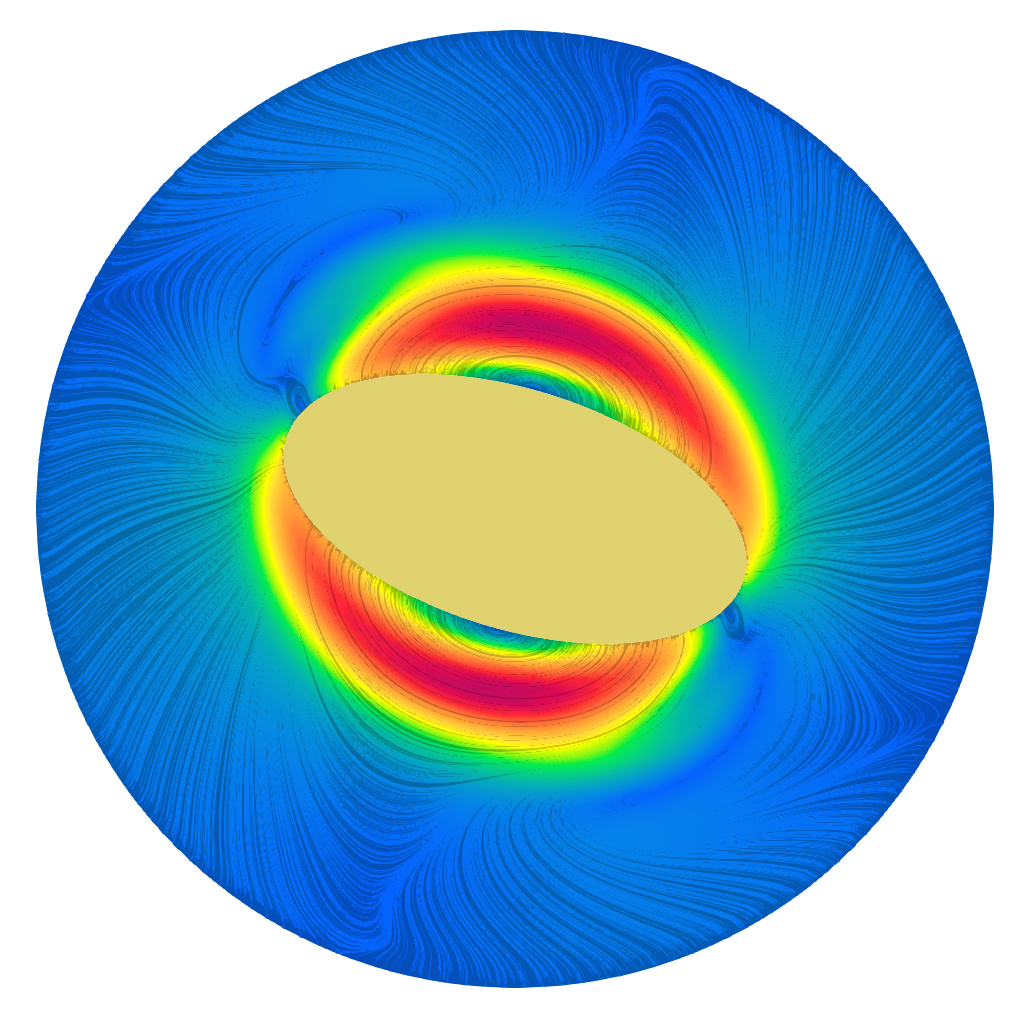}}
\end{figure}	
\begin{figure}[htp!]	
		\subfigure[$ t = 4 $]{\includegraphics[trim = .1cm .1cm .1cm .1cm, clip=true,width=0.24\textwidth,height=0.24\textwidth]{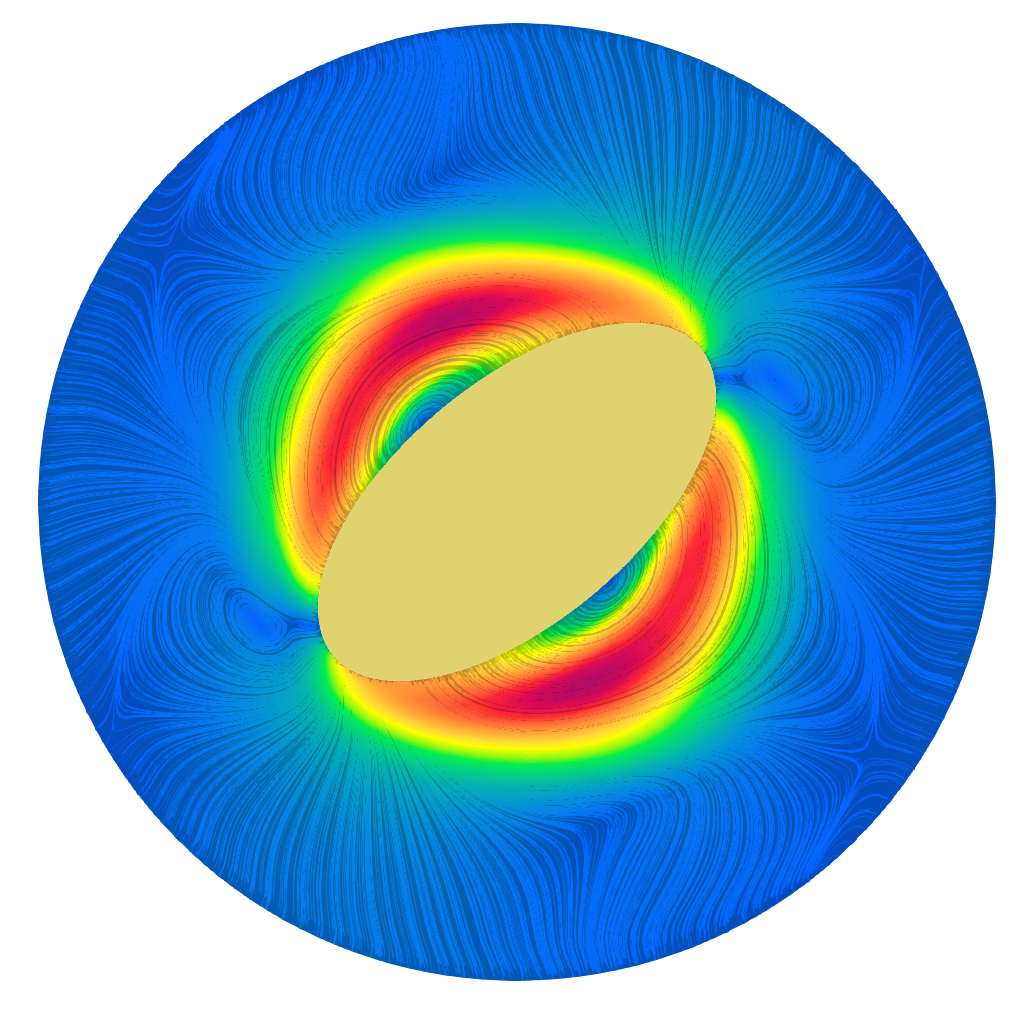}}
	\subfigure[$ t = 4.5 $]{\includegraphics[trim = .1cm .1cm .1cm .1cm, clip=true,width=0.24\textwidth,height=0.24\textwidth]{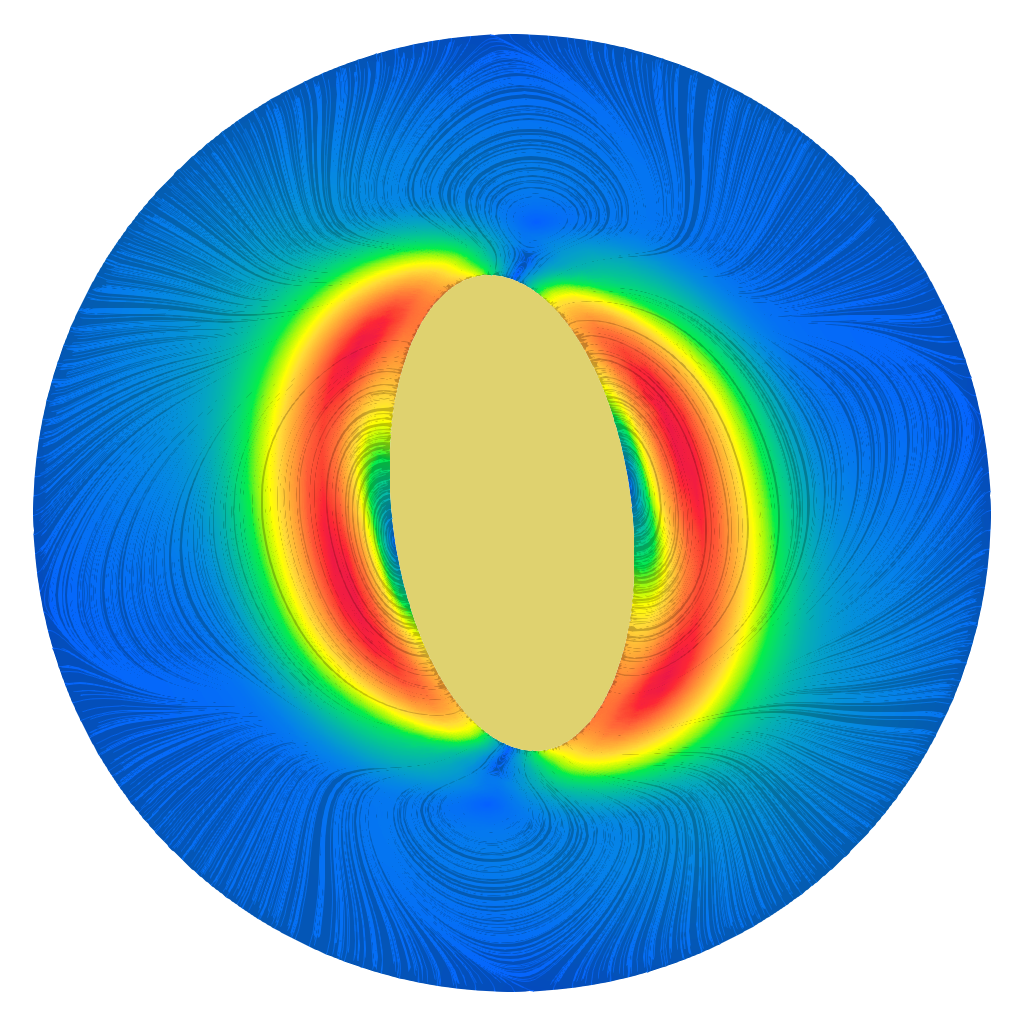}}
	\subfigure[$ t = 5 $]{\includegraphics[trim = .1cm .1cm .1cm .1cm, clip=true,width=0.24\textwidth,height=0.24\textwidth]{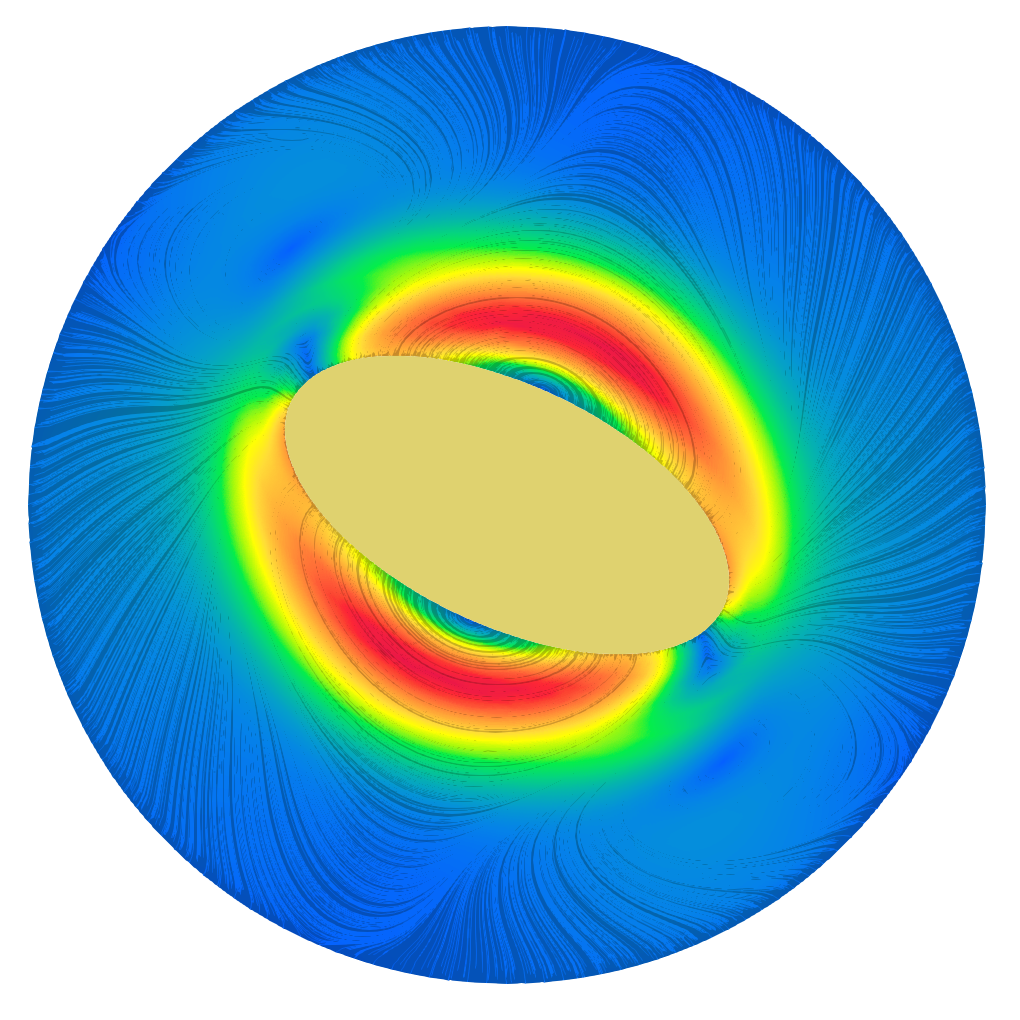}}
	\subfigure[$ t = 5.5 $]{\includegraphics[trim = .1cm .1cm .1cm .1cm, clip=true,width=0.24\textwidth,height=0.24\textwidth]{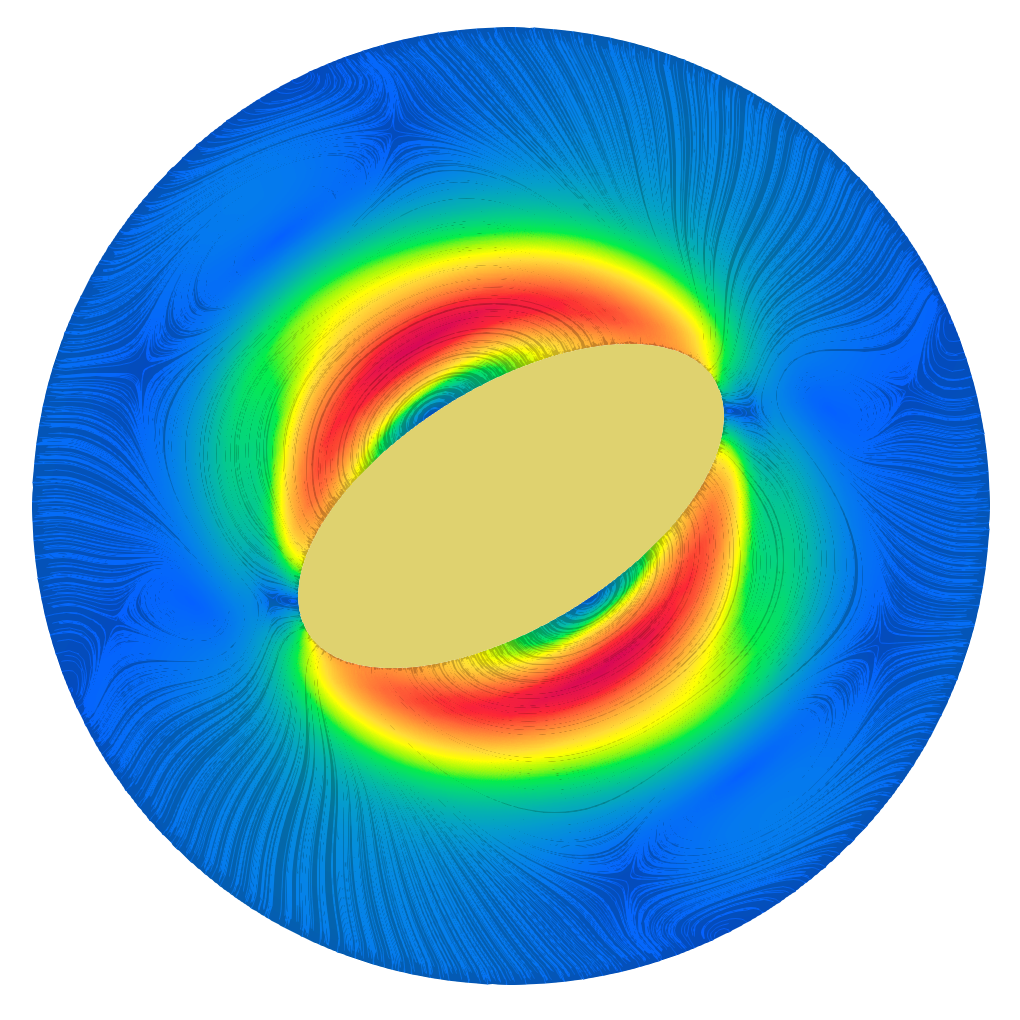}} 
	\subfigure[$ t = 6 $]{\includegraphics[trim = .1cm .1cm .1cm .1cm, clip=true,width=0.24\textwidth,height=0.24\textwidth]{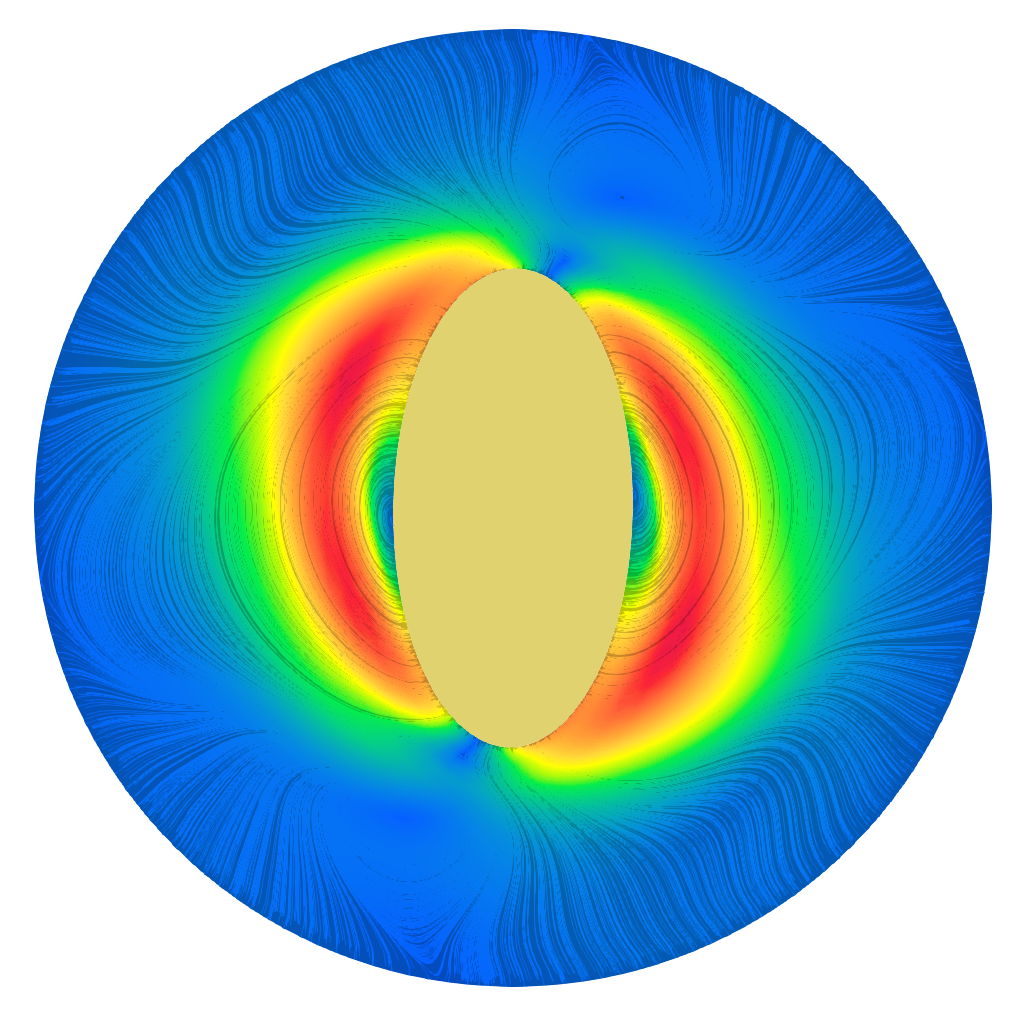}}
	\subfigure[$ t = 6.5 $]{\includegraphics[trim = .1cm .1cm .1cm .1cm, clip=true,width=0.24\textwidth,height=0.24\textwidth]{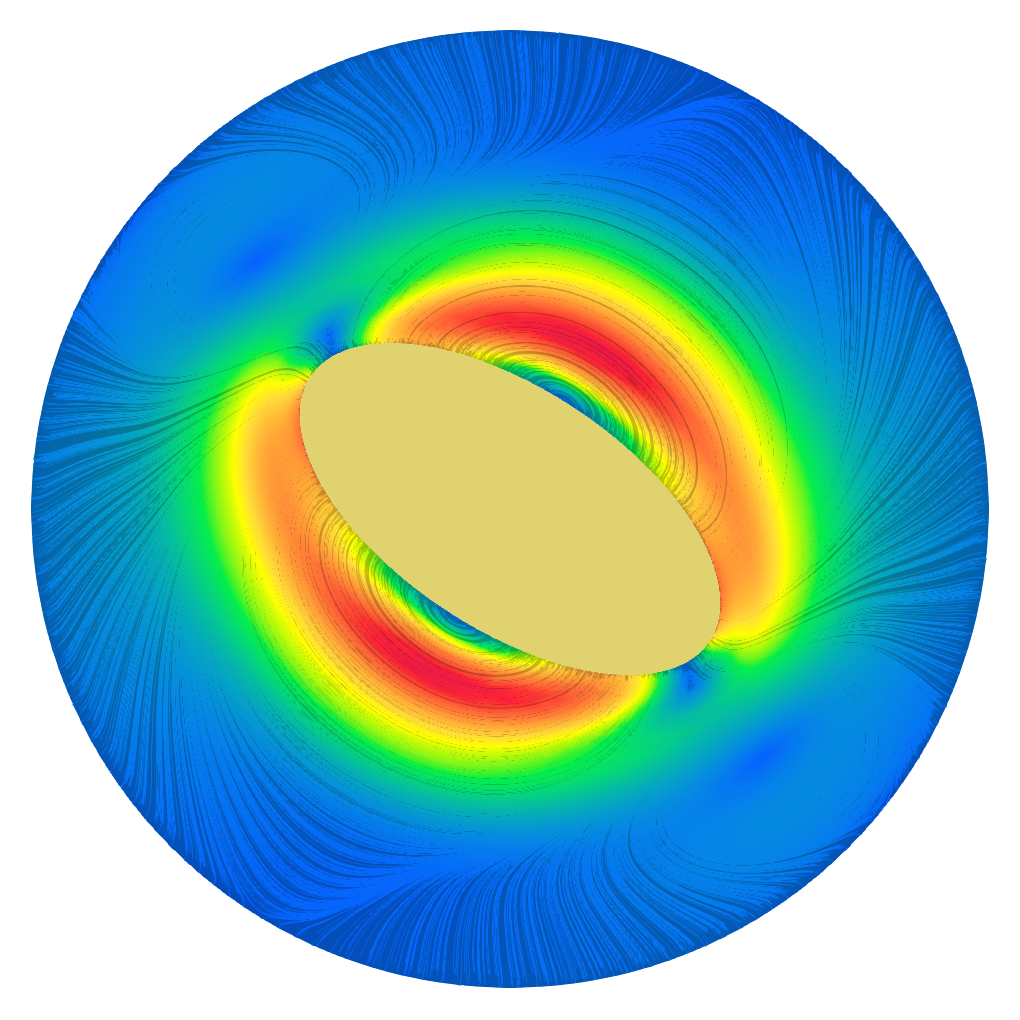}}
	\subfigure[$ t = 7 $]{\includegraphics[trim = .1cm .1cm .1cm .1cm, clip=true,width=0.24\textwidth,height=0.24\textwidth]{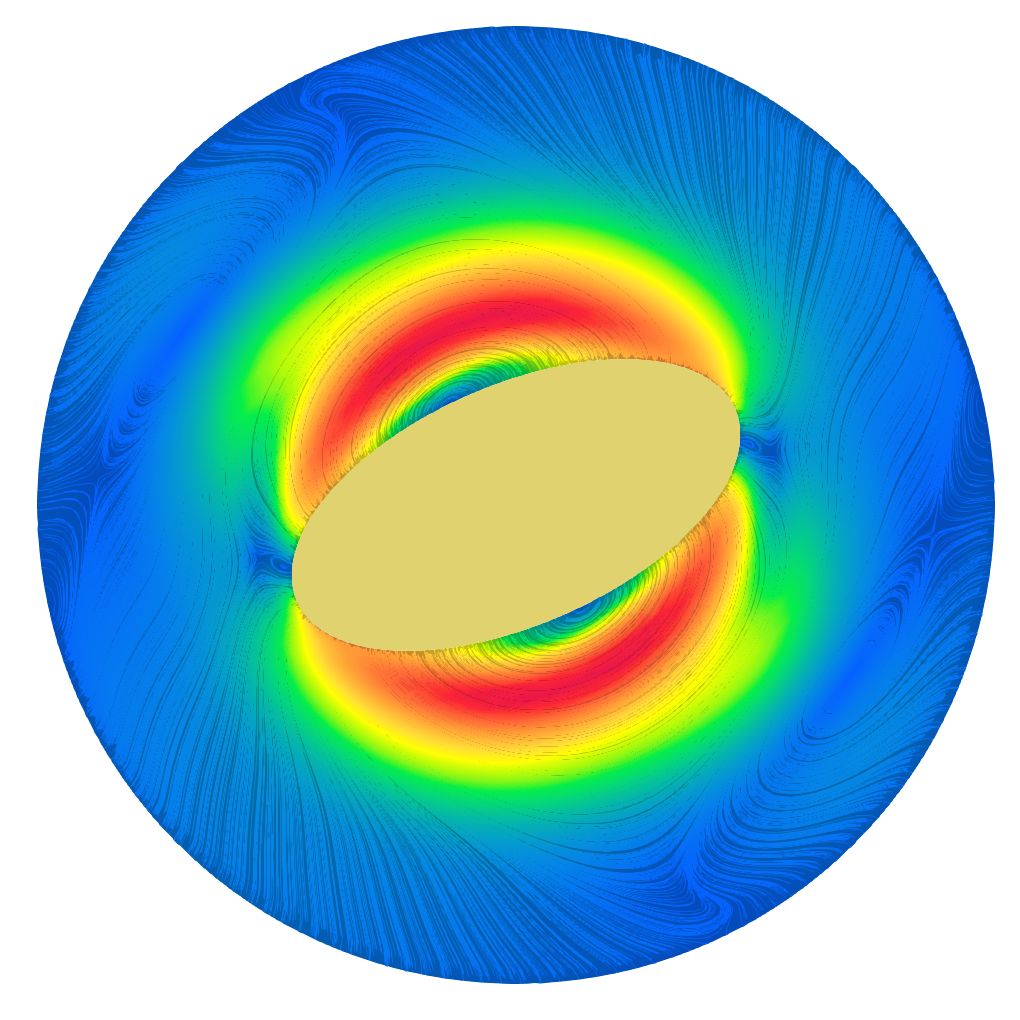}}
	\subfigure[$ t = 7.5 $]{\includegraphics[trim = .1cm .1cm .1cm .1cm, clip=true,width=0.24\textwidth,height=0.24\textwidth]{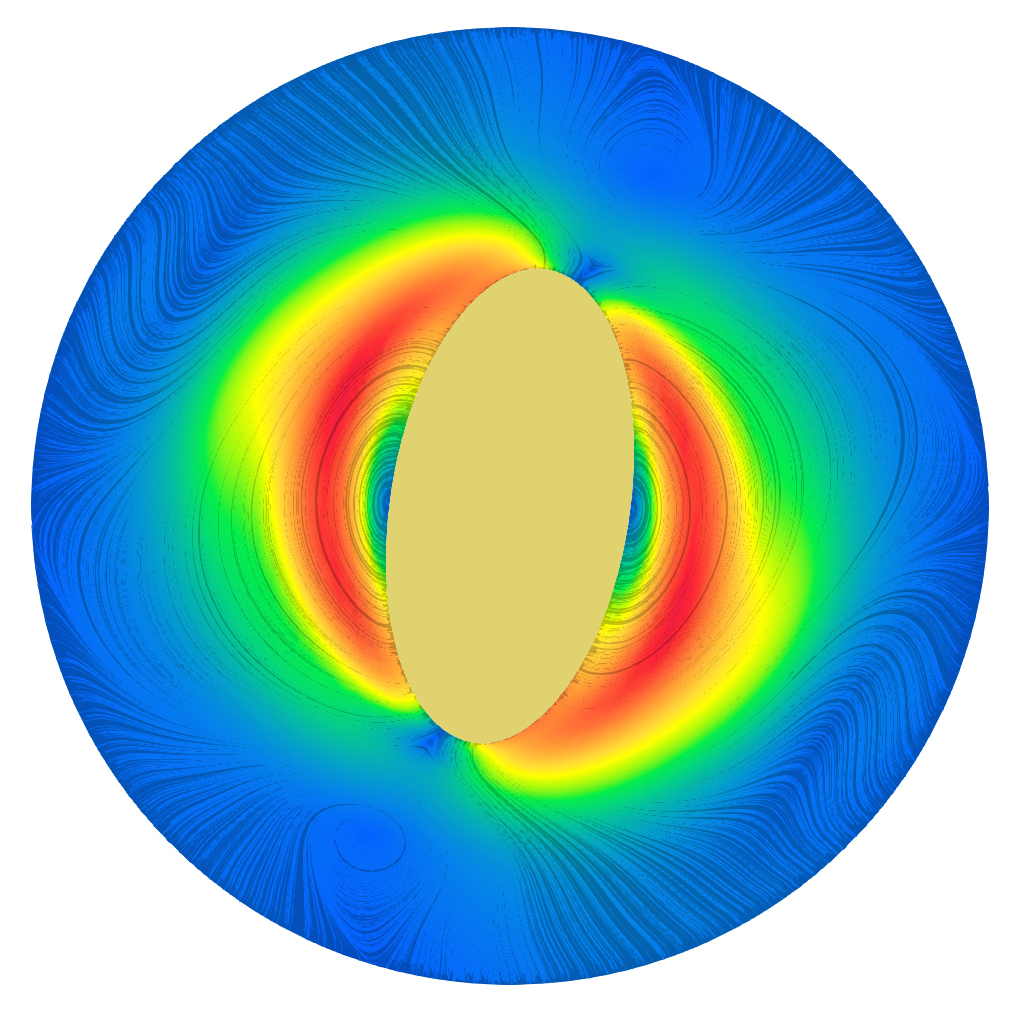}}
	\subfigure[The range of $ |u| $]{\includegraphics[trim = .1cm .1cm .1cm .1cm, clip=true,width=0.8\textwidth,height=0.05\textwidth]{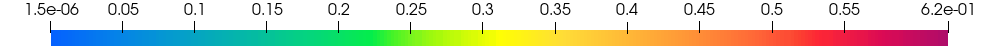}}
	\caption{Flow of the fluid driven by propeller rotation.}\label{Fig-rotation}
\end{figure}

\end{example}

\pagebreak

\newpage
\appendix
\section{$ L^{2} $ and $ H^{-1} $ Estimates for Ritz projection}\label{appendix-A}
\begin{lemma}\label{L^2-duality-1}
For arbitrary $\chi_h\in\mathring V_h$, if we let $g_h=R_hv-\chi_h$ in \eqref{Nitche_trick}, then \eqref{L^2-duality-1-inequality} holds.
\end{lemma}
\begin{proof}
    By employing \eqref{stokes-ritz} for $g_h=R_hv-\chi_h$, we obtain
\begin{align*}
	& \left|(\nabla\widetilde\psi_h, \nabla g_h )_{\Omega_h(t)}-({\rm div}g_h, \widetilde\varphi_h)_{\Omega_h(t)}\right|\\
	\leq & \left|(\nabla(v-\chi_h), \nabla \widetilde{\psi})_{\Omega_h(t)}+({\rm div}\widetilde{\psi}, R_hq-q)_{\Omega_h(t)}-({\rm div}(v-\chi_h), \widetilde{\varphi})_{\Omega_h(t)}\right|\\
	&+Ch\|g_h\|_{\Omega_h(t)}\left(\|\nabla(v-\chi_h)\|_{\Omega_h(t)}+\|R_hq-q\|^{'}_{\Omega_h(t)}\right).
\end{align*}	
It is important to note that
$$
|({\rm div}(v-\chi_h), \widetilde{\varphi})_{\Omega_h(t)}|\leq |(v-\chi_h, \nabla\widetilde{\varphi})_{\Omega_h(t)}|+|(v\cdot \mathbf{n}, \widetilde{\varphi})_{\partial\Omega_h(t)}|,
$$
where $\mathbf{n}$ represents the outward normal vector of the piece-wise smooth boundary $\partial\Omega_h(t)$. Utilizing Lemma \ref{boundary-integral}, we deduce
\begin{equation}\label{bounadry-integral-comp}
|(v\cdot \mathbf{n}, \widetilde{\varphi})_{\partial\Omega_h(t)}|\leq C\|v\|_{L^\infty(\partial\Omega_h(t))}\|\widetilde{\varphi}\|_{L^1(\partial\Omega_h(t))}\leq C\|v\|_{L^\infty(\partial\Omega_h(t))}\|\widetilde{\varphi}\|_{H^1(\mathbb{R}^d)}.
\end{equation}
Since ${\rm div}\widetilde{\psi}|_{\Omega(t)}=0$, we obtain
\begin{align}\label{div-boundary-skin}
|({\rm div}\widetilde{\psi}, R_hq-q)_{\Omega_h(t)}|&\leq \|{\rm div}\widetilde{\psi}\|_{L^2(\Omega_h(t)\setminus\Omega(t))}\|R_hq-q\|^{'}_{\Omega_h(t)}\nonumber\\
&\leq |\Omega_h(t)\backslash\Omega(t)|^{1/2-1/6}\|{\rm div}\widetilde{\psi}\|_{L^6(\mathbb{R}^d)}\|R_hq-q\|^{'}_{\Omega_h(t)}\nonumber\\
&\leq Ch\|\widetilde{\psi}\|_{H^2(\mathbb{R}^d)}\|R_hq-q\|^{'}_{\Omega_h(t)},
\end{align}
where Sobolev embedding and the condition $r\geq 2$ are employed to establish the last inequality. Finally,
\begin{align*}
	|(\nabla(v-\chi_h), \nabla \widetilde{\psi})_{\Omega_h(t)}|&\leq |(v-\chi_h, \Delta \widetilde{\psi})_{\Omega_h(t)}|+|(v,\mathbf{n}\cdot\nabla \widetilde{\psi})_{\partial\Omega_h(t)}|\\
	&\leq C\|v-\chi_h\|_{\Omega_h(t)}\|\widetilde{\psi}\|_{H^2(\mathbb{R}^d)}+C\|v\|_{L^\infty(\partial\Omega_h(t))}\|\widetilde{\psi}\|_{H^2(\mathbb{R}^d)}.
\end{align*}
Combining the above estimates completes the proof.
\end{proof} 

\begin{lemma}\label{negative-norm-pressure-error}
	Let $\widehat{R_hq}:=R_hq-q+\overline{q}$, where $ \overline{q} $ is the average of $ q $ over $ \Omega_{h}(t) $. For each $0\neq \lambda\in H^1(\mathbb{R}^d)$, \eqref{estimate-Rhq-hat} holds.
\end{lemma}
\begin{proof}
	For each $\lambda\in H^1(\mathbb{R}^d)$, let $\overline\lambda_* :=\frac{1}{|\Omega(t)|}\int_{\Omega(t)}\lambda\,\mathrm{d} x$ denote its average over $\Omega(t)$. By Lemma \ref{div-regular-inverse}, there exists $\chi\in H^2(\Omega(t))\cap H^1_0(\Omega(t))$ such that
$$
{\rm div} \chi=\lambda-\overline\lambda_* \quad\text{ in }\Omega(t), \quad \|\chi\|_{H^2(\Omega(t))}\leq C\|\lambda\|_{H^1(\mathbb{R}^d)}.
$$ 
We extend $\chi$ to $\widetilde{\chi}\in H^2(\mathbb{R}^d)$ as mentioned in \eqref{stein-extension}. By decomposing the integral, we have
\begin{align*}
	|(\widehat{R_hq}, \lambda)_{\Omega_h(t)}|&=|(\widehat{R_hq}, \lambda-\overline\lambda_*)_{\Omega_h(t)}|\\
	&\le|(\widehat{R_hq},\lambda-\overline\lambda_*)_{\Omega_h(t)\setminus\Omega(t)}|+|(\widehat{R_hq}, {\rm div}\widetilde{\chi})_{\Omega_h(t)\cap \Omega(t)}|\\
	&\leq|(\widehat{R_hq},\lambda-\overline\lambda_*)_{\Omega_h(t)\setminus\Omega(t)}|+|(\widehat{R_hq}, {\rm div}\widetilde{\chi})_{\Omega_h(t)}|+|(\widehat{R_hq}, {\rm div}\widetilde{\chi})_{\Omega_h(t)\setminus \Omega(t)}|.
\end{align*}
To estimate the boundary-skin integral, we use the following inequalities:
\begin{align}\label{lemma-3.8-tmp}
	|(\widehat{R_hq},\lambda-\overline\lambda_*)_{\Omega_h(t)\setminus\Omega(t)}|&\leq |\overline\lambda_*|\|\widehat{R_hq}\|_{\Omega_h(t)}|\Omega_h(t)\backslash\Omega(t)|^{1/2}+\|\widehat{R_hq}\|_{\Omega_h(t)}|\Omega_h(t)\backslash\Omega(t)|^{1/3}\|\lambda\|_{L^6(\mathbb{R}^d)}\notag\\
	&\leq Ch\|\widehat{R_hq}\|_{\Omega_h(t)}\|\lambda\|_{H^1(\mathbb{R}^d)},
\end{align}
where we have used Holder's inequality, Sobolev embedding $H^1(\mathbb{R}^d)\hookrightarrow L^6(\mathbb{R}^d)$, and the fact $r\geq 2$. Similarly, we can prove that
\begin{align*}
	|(\widehat{R_hq}, {\rm div}\widetilde{\chi})_{\Omega_h(t)\setminus \Omega(t)}|\leq Ch\|\widehat{R_hq}\|_{\Omega_h(t)}\|\widetilde{\chi}\|_{H^2(\mathbb{R}^d)}.
\end{align*}
Since $\widetilde{\chi}|_{\partial\Omega(t)}=0$, we can interpolate $\widetilde{\chi}$ to a function $\chi_h\in \mathring V_h$ using $I_h^{\mathring V}$, i.e., $\chi_h=I_h^{\mathring V}\widetilde{\chi}$. Then we have
\begin{align*}
	&|(\widehat{R_hq}, {\rm div}\widetilde{\chi})_{\Omega_h(t)}|\\
	\leq &|(\widehat{R_hq}, {\rm div}(\widetilde{\chi}-\chi_h))_{\Omega_h(t)}|+	|(\widehat{R_hq}, {\rm div}\chi_h)_{\Omega_h(t)}|\\
	=&|(\widehat{R_hq}, {\rm div}(\widetilde{\chi}-\chi_h))_{\Omega_h(t)}|+|(\nabla(R_hv-v),\nabla (\chi_h - \widetilde{\chi}))_{\Omega_h(t)}| + |(\nabla (R_{h}v - v), \nabla \widetilde{\chi})_{\Omega_{h}(t)}|\\
	\leq& Ch\|\widetilde{\chi}\|_{H^2(\mathbb{R}^d)}\left(\|\widehat{R_hq}\|_{\Omega_h(t)}+\|\nabla(R_hv-v)\|_{\Omega_h(t)}\right)+|(\nabla(R_hv-v),\nabla \widetilde{\chi})_{\Omega_h(t)}|.
\end{align*}
Integrating by parts and dealing the boundary integral term as in \eqref{bounadry-integral-comp}, we obtain
\begin{align*}
	|(\nabla(R_hv-v),\nabla \widetilde{\chi})_{\Omega_h(t)}|&\leq C\|R_hv-v\|_{\Omega_h(t)}\|\widetilde{\chi}\|_{H^2(\mathbb{R}^d)}+C\|v\|_{L^\infty(\partial\Omega_h(t))}\|\widetilde{\chi}\|_{H^2(\mathbb{R}^d)}.
\end{align*}
Combining the above estimates completes this proof.
\end{proof}
\begin{remark}
	In the derivation of \eqref{lemma-3.8-tmp}, the condition $r \geq 2$ is imposed, which, upon closer examination, is not strictly necessary. In fact, when considering the MINI-finite element space where $r = 1$, a more accurate estimation of $(\widehat{R_h q},\lambda)_{\Omega_h(t)\setminus\Omega(t)}$ can be achieved by utilizing the findings from Lemma \ref{boundary-skin} to estimate $\|\widehat{R_h q}\|_{\Omega_h(t)\setminus\Omega(t)}$.
\end{remark}

\section{$\ell^1 L^2$-estimates for discretized Stokes equations}\label{ell1L2-estimate-section}
In this section, we provide a comprehensive proof for \eqref{ell1L2-form1}, which is used in proving the main result of this article, i.e., Theorem \ref{main-theorem}. To prove \eqref{ell1L2-form1}, we begin with introducing an equivalent result, i.e., Lemma \ref{ell1L2-lemma2}, which serves as a foundation for our subsequent analysis. 
\begin{lemma}\label{ell1L2-lemma2}
Let $\left\{(\omega_h^n,\varphi_h^n)\in \mathring V_h\times Q_h\right\}_{0\leq n\leq M}$ be the solutions of 
\begin{subequations}\label{backward-euler}
	\begin{align}
		\label{backward-euler-1}
		&\left(d_t \omega^{n+1}_h, v_h\right)_{\Omega_h^{n+1/2}}+\left(\nabla \omega_h^{n+1},\nabla v_h\right)_{\Omega_h^{n+1/2}}-(c^{n+1}_h\cdot \nabla\omega^{n+1}_h, v_h)_{\Omega^{n+1/2}_h}-({\rm div}v_h, \varphi^{n+1}_h)_{\Omega^{n+1/2}_h} \nonumber\\
		=&(g^{n+1}_h, v_h)_{\Omega^{n+1/2}_h} \quad \forall v_h\in \mathring V_h,\\
		\label{backward-euler-2}
		&\left({\rm div}\omega^{n+1}_h, q_h\right)_{\Omega^{n+1/2}_h}=0 \quad \forall q_h\in Q_h,
	\end{align}
\end{subequations}	
with initial value $\omega_h^0=0$. Assume that
\begin{align}
	\|c_h^{n+1}\|_{W^{1,\infty}(\Omega_h^{n+1})}\leq C, \quad \|c_h^{n+1}-c_h^n\|_{W^{1,\infty}(\Omega_h^{n+1})}\leq C\tau.
\end{align}
Then we establish the following estimate
\begin{align}
	&\|\omega_h^{M}\|_{\Omega_h^{M}}+\sum_{n=0}^{M-1}\tau\left(\|\nabla\omega_h^{n+1}\|_{\Omega_h^{n+1}}+\|\nabla\varphi_h^{n+1}\|_{\Omega_h^{n+1}}+\|d_t\omega_h^{n+1}\|_{\Omega_h^{n+1}}\right)\notag\\
	\leq & C \min \left\{\ell_\tau,\ell_h\right\}\sum_{n=0}^{M-1}\tau\|g^{n+1}_h\|_{\Omega_h^{n+1/2}}, \quad \mbox{where $\ell_\tau:=\ln(1/\tau+1),\ell_h:=\ln(1/h+1)$}.
\end{align}
\end{lemma}

We now introduce some preliminary results before proving the lemma. We define the solution map $\mathcal{S}_{m,n}$ as follows:
\begin{definition}\label{definition-A-2}
	For $0\leq m < n\leq M$, let $u_h\in \mathring V_h$, and consider the unique solution $\{(u_h^k,p_h^k)\in \mathring V_h\times Q_h\}_{m+1\leq k\leq M}$ of the following equations, for any $ v_{h} \in \mathring V_{h} $, $ q_{h} \in Q_{h} $.
	\begin{subequations}\label{backward-euler-homo}
		\begin{align}
			\label{backward-euler-homo1}
			\left(d_t u^{k}_h, v_h\right)_{\Omega_h^{k-1/2}}+\left(\nabla u_h^{k},\nabla v_h\right)_{\Omega_h^{k-1/2}}-(c^{k}_h\cdot \nabla u^{k}_h, v_h)_{\Omega^{k-1/2}_h}-({\rm div}v_h, p^{k}_h)_{\Omega^{k-1/2}_h}&=0, \\
			\label{backward-euler-homo2}
			\left({\rm div}u_h^k, q_h\right)_{\Omega^{k-1/2}_h}&=0,
		\end{align}
	\end{subequations}	
	with initial value $u_h^m=u_h$. We define $\mathcal{S}_{m,n} u_h:= u_h^n$. Additionally, we adopt the convention that $\mathcal{S}_{m,m}u_h=u_h$.
\end{definition}
From the definition, it follows immediately that,
\begin{align}
	\mathcal{S}_{n,l}\circ \mathcal{S}_{m,n}=\mathcal{S}_{m,l} \quad \forall m\leq n\leq l.
\end{align}
Furthermore, we have the expression
\begin{align}
	\omega_h^{n+1}=\mathcal{S}_{n,n+1}(\omega_h^n+\tau g_h^{n+1}).
\end{align}
By applying induction, we obtain the following formulas:
\begin{align}
	\omega_h^{n+1}=&\sum_{m=1}^{n+1}\tau \mathcal{S}_{m-1,n+1} (g_h^m),\label{omega_h-formula}\\
	d_t\omega_h^{n+1}=&\sum_{m=1}^n \left(\mathcal{S}_{m-1,n+1}-\mathcal{S}_{m-1,n}\right)(g_h^m)+\mathcal{S}_{n,n+1}(g_h^{n+1}).\label{d_tomega_h-formula}
\end{align}

The following lemma provides estimates for $\mathcal{S}_{m,n}$ introduced in Definition \ref{definition-A-2}.
\begin{lemma}
	For $ u_{h} $ and $ \mathcal S_{m,n}u_{h} $ as defined in Definition \ref{definition-A-2}, the following estimates hold
	\begin{align}
		\|\mathcal{S}_{m,n}u_h\|_{\Omega_h^n}&\leq C\|u_h\|_{\Omega_h^m}, \label{solution-map-estimate1}\\
		\frac{1}{\tau}\|\mathcal{S}_{m,n+1}u_h-\mathcal{S}_{m,n}u_h\|_{\Omega_h^{n+1}}&\leq C \min \bigg\{\frac{1}{(n-m+1)\tau},\frac{1}{h^{2}}\bigg\}\|u_h\|_{\Omega_h^m}. \label{solution-map-estimate2}
\end{align}
\end{lemma}

By shifting the time and considering $\Omega_h^m$ as $\Omega_h^0$, we can simplify the problem. Utilizing the assumption that the mesh velocity $w_h$ is uniformly bounded with respect to time, it suffices to prove the lemma for the unique solution $\{(u_h^{n},p_h^n)\in \mathring V_h\times Q_h\}_{n=1}^{M}$ of the following system:
\begin{subequations}\label{backward-euler-homo-eg}
	\begin{align}
		\label{backward-euler-homo-eg1}
		\left(d_t u^{n}_h, v_h\right)_{\Omega_h^{n-1/2}}+\left(\nabla u_h^{n},\nabla v_h\right)_{\Omega_h^{n-1/2}}-(c^{n}_h\cdot \nabla u^{n}_h, v_h)_{\Omega^{n-1/2}_h}-({\rm div}v_h, p^{n}_h)_{\Omega^{n-1/2}_h} &=0,\\
		\label{backward-euler-homo-eg2}
		\left({\rm div}u^{n}_h, q_h\right)_{\Omega^{n-1/2}_h}&=0 ,
	\end{align}
\end{subequations}
where $v_{h} \in \mathring V_{h}$ and $q_{h} \in Q_{h}$, with the initial value $u_h^0$. The goal is to prove the following two lemmas, which establish the inequalities:
\begin{align}\label{A3-lemma-goal}
	\|u_h^n\|_{\Omega_h^n}\leq C\|u_h^0\|_{\Omega_h^0}, \quad \|d_tu_h^{n+1}\|\leq C\min\left\{\frac{1}{(n+1)\tau},\frac{1}{h^{2}}\right\}\|u_h^0\|_{\Omega_h^0}.
\end{align} 
To prove \eqref{A3-lemma-goal}, we will first establish the following foundational results: 
\begin{lemma}\label{appendix-energy-lemma1}
Let $\{(u_h^{n},p_h^n)\in \mathring V_h\times Q_h\}_{n=1}^{M}$ be the solution of \eqref{backward-euler-homo-eg}. Then the following estimates hold for $0\leq m\leq M$
\begin{align}
	\|u_h^{m+1}\|_{\Omega_h^{m+1}}+\bigg(\sum_{n=0}^{m}\tau\|\nabla u_h^{n+1}\|^2_{\Omega_h^{n+1}}\bigg)^{1/2}&\leq C\|u_h^0\|_{\Omega_h^0},\label{appendix-energy-1}\\
	\|\nabla u_h^{m+1}\|_{\Omega_h^{m+1}}+\bigg(\sum_{n=[m/2]+1}^m\tau\|d_t u_h^{n+1}\|^2_{\Omega_h^{n+1}}\bigg)^{1/2}&\leq Ct_{m+1}^{-1/2}\|u_h^0 \|_{\Omega_h^0},\label{appendix-energy-2}
\end{align}
where $[m/2]$ denotes the largest integer less than or equal to $m/2$.
\end{lemma}
\begin{proof}
	Using a similar approach to the proof of Lemma \ref{energy-error-estimate-0}, we can derive the first inequality by testing \eqref{backward-euler-homo-eg1} with $v_h=u_h^{n}$.
%
For the second inequality, we test \eqref{backward-euler-homo-eg1} with $v_h=d_t u_h^{n}$. Estimating the terms caused by domain variation by Lemma \ref{transport-theorem-lemma} and applying Young's inequality to $(c^{n}_h\cdot \nabla u^{n}_h, d_tu_h^{n})_{\Omega^{n-1/2}_h}$, we obtain:
\begin{align}\label{A3-lemma-temp2}
	& \|d_t u_h^{n}\|^2_{\Omega_h^{n-1/2}}+\frac{1}{2\tau}\left(\|\nabla u_h^{n}\|^2_{\Omega_h^{n}}-\|\nabla u_h^{n-1}\|^2_{\Omega_h^{n-1}}\right)\notag\\
	\leq & \frac{1}{4}\|d_t u_h^{n}\|^2_{\Omega_h^{n-1/2}}+ C\left(\|\nabla u_h^{n}\|^2_{\Omega_h^{n}}+\|\nabla u_h^{n-1}\|^2_{\Omega_h^{n-1}}\right)+\left({\rm div}d_t u_h^{n}, p_h^{n}\right)_{\Omega_h^{n-1/2}}.
\end{align}
Under the inf-sup condition, it follows from \eqref{backward-euler-homo-eg} that
\begin{align}\label{A3-lemma-temp3}
	\left({\rm div}d_t u_h^{n}, p_h^{n}\right)_{\Omega_h^{n-1/2}}=-\frac{1}{\tau}\left({\rm div} u_h^{n-1}, p_h^{n}\right)_{\Omega_h^{n-1/2}}\leq C\|\nabla u_h^{n-1}\|_{\Omega_h^{n-1}}\left(\|d_t u_h^{n}\|_{\Omega_h^{n}}+\|\nabla u_h^{n}\|_{\Omega_h^{n}}\right).
\end{align} 
Substituting \eqref{A3-lemma-temp3} into \eqref{A3-lemma-temp2}, we obtain
\begin{align}\label{A3-lemma-temp5}
	\|d_t u_h^{n}\|^2_{\Omega_h^{n}}+\frac{1}{\tau}\left(\|\nabla u_h^{n}\|^2_{\Omega_h^{n}}-\|\nabla u_h^{n-1}\|^2_{\Omega_h^{n-1}}\right)\leq C\left(\|\nabla u_h^{n}\|^2_{\Omega_h^{n}}+\|\nabla u_h^{n-1}\|^2_{\Omega_h^{n-1}}\right).
\end{align}
Using a technique similar to \cite[Appendix B]{LiMA-2022}, we construct a smooth cut-off function $0\leq \chi_m(t)\leq 1$ with the following properties: $\chi_{m}(t)=0$ for all $0\leq t\leq \frac{t_{m+1}}{2}$, $\chi_m(t)=1$ for all $t\geq t_{m+1}$, and $|\chi_m'(t)|\leq Ct_{m+1}^{-1}$ for all $t>0$. By multiplying \eqref{A3-lemma-temp5} with $\chi_m(t_{n})$ and summing from $n=1$ to $n=m+1$, taking into account of estimate \eqref{appendix-energy-1} we obtain for $m\geq 1$ that:
\begin{align}\label{A3-lemma-temp7}
	\|\nabla u_h^{m+1}\|^2_{\Omega_h^{m+1}}+\sum_{n=1}^{m+1}\tau\chi_m(t_{n})\|d_t u_h^{n}\|^2_{\Omega_h^{n}}\leq Ct_{m+1}^{-1}\|u_h^0 \|^2_{\Omega_h^0}.
\end{align}
For the special case $m=0$, \eqref{appendix-energy-2} trivially holds, thus \eqref{appendix-energy-2} follows directly from \eqref{A3-lemma-temp7}.
\end{proof}

\begin{lemma}
	Let $u_h^{n}$ be the solution of \eqref{backward-euler-homo}. 
	Then the following result holds for $0\leq m\leq M$:
	\begin{align}\label{appendix-energy-3}
		\|d_tu_h^{m+1}\|_{\Omega_h^{m+1}}\leq  C\min\bigg\{\frac{1}{(m+1)\tau},\frac{1}{h^{2}}\bigg\}\|u_h^0\|_{\Omega_h^0}.
	\end{align}
\end{lemma}
\begin{proof}
For $n\geq1$, we take the difference of \eqref{backward-euler-homo-eg1} between steps $n+1$ and $n$ and test the resulting equation with $ v_{h} = d_t u_{h}^{n+1} $. Using a similar method to the proof of Lemma \ref{energy-error-estimate-2}, we obtain: 
\begin{align}\label{A5-lemma-temp1}
	&\frac{1}{2\tau}\left(\|d_tu_h^{n+1}\|^2_{\Omega_h^{n+1}}-\|d_tu_h^{n}\|^2_{\Omega_h^{n}}\right)+\|\nabla d_t u_h^{n+1}\|^2_{\Omega_h^{n+1/2}}\nonumber\\
	\leq& \frac{1}{2}\|\nabla d_t u_h^{n+1}\|^2_{\Omega_h^{n+1/2}}+C\left(\|d_tu_h^{n+1}\|^2_{\Omega_h^{n+1}}+\|d_tu_h^{n}\|^2_{\Omega_h^{n}}\right)+\frac{1}{\tau}\left[\mathcal{I}^{n+1}(u_h^{n+1},  p_h^{n+1})-\mathcal{I}^{n}(u_h^n,  p_h^{n})\right]\nonumber\\
	& +C\left(\|p_h^n\|^2_{\Omega_h^n}+\|\nabla u_h^n\|^2_{\Omega_h^n}+\|p_h^{n+1}\|^2_{\Omega_h^{n+1}}\right), \quad\mbox{(for $n\geq 2$)}
\end{align}
where the bilinear form $\mathcal{I}^{n+1}$ is defined as follows:
\begin{align}
	\mathcal{I}^{n+1}(v_h, q_h):=\frac{-1}{2\tau}\int_{t_{n-1}}^{t_{n+1}}\left({\rm div} w_h(t){\rm div}v_h, q_h\right)_{\Omega_h(t)}-\left(\nabla v_h:\nabla w_h(t)^{\top}, q_h\right)_{\Omega_h(t)} \d t.
\end{align}

Similarly to the proof of Lemma \ref{appendix-energy-lemma1}, by combining discrete Gronwall's inequality and the following estimates
\begin{align}
	\|p_h^n\|_{\Omega_h^n}\leq C\left(\|d_t u_h^n\|_{\Omega_h^n}+\|\nabla u_h^n\|_{\Omega_h^n}\right) \quad\mbox{and}\quad |\mathcal{I}^n(u_h^n,p_h^n)|\leq C\|\nabla u_h^n\|_{\Omega_h^n}\|p_h^n\|_{\Omega_h^n},
\end{align}
after multiplying \eqref{A5-lemma-temp1} with $\chi_m(t_{n+1})$ and summing from $n=1$ to $m$, we obtain for $m\geq 3$
\begin{align}
	\|d_t u_h^{m+1}\|_{\Omega_h^{m+1}}\leq Ct_{m+1}^{-1}\|u_h^0\|_{\Omega_h^0}, \quad\mbox{(for $m\geq 3$)}.
\end{align}
While for $0\leq m\leq 2$ 
\begin{equation*}
	\|d_tu_h^{m+1}\|\leq \frac{C}{(m+1)\tau}\|u_h^0\|_{\Omega_h^0}.
\end{equation*}
trivially holds, thus we complete the first part of this proof.

As for the bound $ h^{-2} $, we test \eqref{backward-euler-homo-eg1} (at step $m+1$) with $ v_{h} = d_t u_{h}^{m+1} $ and obtain the following inequality
\begin{equation}\label{B5-lem-tmp0}
	\|d_t u_h^{m+1}\|^2_{\Omega_h^{m+1}}\leq C\left(\|\nabla u_h^{m+1}\|_{\Omega_h^{m+1}}\|\nabla d_t u_h^m\|_{\Omega_h^{m+1}}+({\rm div} d_tu_h^{m+1},p_h^{m+1})_{\Omega_h^{m+1/2}}\right).
\end{equation}
While by \eqref{backward-euler-homo-eg1} and transport theorem we have 
\begin{equation}\label{B5-lem-tmp1}
	({\rm div}d_t u_h^{m+1},p_h^{m+1})_{\Omega_h^{m+1/2}}=-\frac{1}{\tau}({\rm div}u_h^m,p_h^{m+1})_{\Omega_h^{m+1/2}}\leq C\|\nabla u_h^m\|_{\Omega_h^m}\|p_h^{m+1}\|_{\Omega_h^{m+1}},
\end{equation}
and from inf-sup condition it follows that
\begin{equation}\label{B5-lem-tmp2}
	\|p_h^{m+1}\|_{\Omega_h^{m+1}}\leq C\left(\|\nabla u_h^{m+1}\|_{\Omega_h^{m+1}}+\|d_tu_h^{m+1}\|_{\Omega_h^{m+1}}\right).
\end{equation}
Substituting \eqref{B5-lem-tmp2} and \eqref{B5-lem-tmp1} into \eqref{B5-lem-tmp0} and using inverse estimate of finite element functions and Young's inequality, we obtain
\begin{equation}
	\|d_t u_h^{m+1}\|^2_{\Omega_h^{m+1}}\leq Ch^{-2}(\|\nabla u_h^{m+1}\|^2_{\Omega_h^{m+1}}+\|\nabla u_h^{m}\|^2_{\Omega_h^{m+1}})\leq Ch^{-4}\|u_h^0\|^2_{\Omega_h^0},
\end{equation}
where in deducing the last inequality we have also used \eqref{appendix-energy-1}. Thus we complete the second part of this proof.

\end{proof}

\begin{proof}[Proof of Lemma \ref{ell1L2-lemma2}]
From \eqref{solution-map-estimate1} we can deduce that for every $0\leq \theta\leq 1$ there holds
\begin{align}
	\frac{1}{\tau}\|\mathcal{S}_{m,n+1}u_h-\mathcal{S}_{m,n}u_h\|_{\Omega_h^{n+1}}\leq C t_{n-m+1}^{\theta-1}h^{-2\theta}\|u_h\|_{\Omega_h^m}. \label{solution-map-estimate3}
\end{align}
By using \eqref{d_tomega_h-formula} and the estimate \eqref{solution-map-estimate1}-\eqref{solution-map-estimate3}, we obtain the inequality 
\begin{align}\label{dt-w-h-estimate}
		\|d_t\omega_h^{n+1}\|_{\Omega_h^{n+1}}\leq C\sum_{m=1}^{n+1} \tau t^{\theta-1}_{n-m+2}h^{-2\theta}\|g_h^m\|_{\Omega_h^{m-1/2}}.
\end{align}
Summing up \eqref{dt-w-h-estimate} from $ n = 0 $ to $ M - 1 $, we have 
\begin{align}
	\sum_{n=0}^{M-1} \tau\|d_t \omega_h^{n+1}\|_{\Omega_h^{n+1}}&\leq C\sum_{n=0}^{M-1}\tau\sum_{m=0}^{n}\tau t^{\theta-1}_{n-m+1}h^{-2\theta}\|g_h^{m+1}\|_{\Omega_h^{m+1/2}}\nonumber\\
	&\leq \sum_{m=0}^{M-1}\tau\|g_h^{m+1}\|_{\Omega_h^{m+1/2}}\sum_{n=m}^{M-1}\tau t^{\theta-1}_{n-m+1}h^{-2\theta}
\end{align}
If we take $\theta=0$, then we have 
\begin{equation*}
	\sum_{n=m}^{M-1}\tau t^{\theta-1}_{n-m+1}h^{-2\theta}=\sum_{n=m}^{M-1}\tau t^{-1}_{n-m+1}\leq C\ln(1/\tau+1)
\end{equation*}
If we take $\theta=1/\ln(1+1/h)$ then there holds
\begin{equation*}
	\sum_{n=m}^{M-1}\tau t^{\theta-1}_{n-m+1}h^{-2\theta}\leq C\theta^{-1}h^{-2\theta}\leq C\ln(1/h+1).
\end{equation*}
Therefore, combining the two different choices of $\theta$ we obtain
\begin{equation}
	\sum_{n=0}^{M-1} \tau\|d_t \omega_h^{n+1}\|_{\Omega_h^{n+1}}\leq C\min\{\ell_\tau,\ell_h\}\sum_{n=0}^{M-1}\tau\|g_h^{n+1}\|_{\Omega_h^{n+1/2}}
\end{equation}
Using \eqref{omega_h-formula} and the estimate \eqref{appendix-energy-1}, we have
\begin{align}
	\|\omega_h^{n+1}\|_{\Omega_h^{n+1}}\leq C\tau \sum_{m=0}^n \|g_h^{m+1}\|_{\Omega_h^{m+1/2}}.
\end{align}
By testing \eqref{backward-euler-1} with $v_h=\omega_h^{n+1}$ we deduce that
\begin{align}
	\|\nabla \omega_h^{n+1}\|_{\Omega_h^{n+1}}\leq C\left(\|g_h^{n+1}\|_{\Omega_h^{n+1/2}}+\|d_t \omega_h^{n+1}\|_{\Omega_h^{n+1}}+\|\omega_h^{n+1}\|_{\Omega_h^{n+1}}\right).
\end{align}
As a result, we obtain
\begin{equation}
	\|\omega_h^{M}\|_{\Omega_h^{M}}+\sum_{n=0}^{M-1}\tau\left(\|\nabla\omega_h^{n+1}\|_{\Omega_h^{n+1}}+\|d_t\omega_h^{n+1}\|_{\Omega_h^{n+1}}\right)\leq C\min\{\ell_\tau,\ell_h\}\sum_{n=0}^{M-1}\tau\|g^{n+1}_h\|_{\Omega_h^{n+1/2}}.
\end{equation}
Finally, the estimate for $\|\nabla \varphi_h^{n+1}\|_{\Omega_h^{n+1}}$ follows from
\begin{align}
	\|\nabla \varphi_h^{n+1}\|_{\Omega_h^{n+1}}\leq C\left(\|d_t\omega_h^{n+1}\|_{\Omega_h^{n+1}}+\|\nabla \omega_h^{n+1}\|_{\Omega_h^{n+1}}+\|g_h^{n+1}\|_{\Omega_h^{n+1/2}}\right),
\end{align}
which is a corollary of the Lemma \ref{FEM-pressure-gradient-estimate} below.
\end{proof}

\begin{lemma}\label{FEM-pressure-gradient-estimate}
	Let $u_h\in \mathring V_h$, $p_h\in Q_h$, $f_h\in V_h$ satisfy the following equations
	\begin{subequations}\label{Omega-n+1/2-fem}
		\begin{align}
			\left(\nabla u_h, \nabla v_h\right)_{\Omega_h^{n+1/2}}-\left({\rm div} v_h, p_h\right)_{\Omega_h^{n+1/2}}&=\left(f_h, v_h\right)_{\Omega_h^{n+1/2}}\quad \forall v_h\in \mathring V_h,\\
			\left({\rm div} u_h, q_h\right)_{\Omega_h^{n+1/2}}&=0 \quad \forall q_h\in Q_h.
		\end{align}
	\end{subequations}
	Then it holds that
	\begin{align}
		\|\nabla p_h\|_{\Omega_h^{n+1}}\leq C\|f_h\|_{\Omega_h^{n+1/2}}.
	\end{align}
\end{lemma}
\begin{proof}
	We define matrices functions $A_h^{n+1}$, $B_h^{n+1}$, and a scalar function $C_h^{n+1}$ on $\Omega_h^0$ as follows:
\begin{align*}
	&\left(A_h^{n+1}\nabla W_h, \nabla W'_h\right)_{\Omega_h^0}=\left(\nabla W_h, \nabla W_h'\right)_{\Omega_h^{n+1}} \quad \forall W_h, W_h'\in \mathring V_h,\\
	&\left(B_h^{n+1}:\nabla W_h, Y_h\right)_{\Omega_h^0}=\left({\rm div} W_h, Y_h\right)_{\Omega_h^{n+1}}\quad \forall W_h\in \mathring V_h, \;\forall Y_h\in Q_h, \\
	&\left(C^{n+1}_h W_h, W_h'\right)_{\Omega_h^0}=\left(W_h, W_h'\right)_{\Omega_h^{n+1}}\quad \forall W_h, W_h'\in \mathring V_h.
\end{align*}
Here, the matrices $ A_{h}^{n+1} $, $ B_{h}^{n+1} $ and the scaler $ C_{h}^{n+1} $ are defined as follows:
\begin{align*}
	A_h^{n+1}&:=\left({\rm det}\nabla\phi_h^{n+1}\right)\left(\nabla \phi_h^{n+1}\right)^{-1}\left(\nabla \phi_h^{n+1}\right)^{-\top}, \quad 	B_h^{n+1}:=\left({\rm det}\nabla\phi_h^{n+1}\right)\left(\nabla \phi_h^{n+1}\right)^{-\top},\\
	C_h^{n+1}&:=\left({\rm det}\nabla\phi_h^{n+1}\right).
\end{align*}
To facilitate our notation, we introduce the shorthand $A_h^{n+1/2}:= \frac{1}{2}\left(A_h^n+A_h^{n+1}\right)$. By denoting $U_h:= u_h\circ \phi_h^{n+1}\in V_h^0$, $P_h:= p_h\circ \phi_h^{n+1}\in Q_h^0$, $F_h:=f_h\circ \phi_h^{n+1}$, equation \eqref{Omega-n+1/2-fem} can be expressed equivalently as:
	\begin{subequations}\label{Omega-n+1/2-fem-1}
	\begin{align}
		(A_h^{n+1/2}\nabla U_h, \nabla W_h)_{\Omega_h^{0}}-(B_h^{n+1/2}:\nabla W_h, P_h)_{\Omega_h^{0}}&=(C_h^{n+1/2}F_h, W_h)_{\Omega_h^{0}},\quad \forall W_h\in \mathring V_h,\\
		(B^{n+1/2}_h:\nabla U_h, Y_h)_{\Omega_h^0}&=0, \quad \forall Y_h\in Q_h.
	\end{align}
\end{subequations}
We extend $F_h$ by zero outside of $\Omega_h^0$, and define $U\in H_1^0(\Omega^0)$, $P\in L^2_0(\Omega^0)$ to be solution of 
\begin{subequations}\label{Omega-n+1/2-fem-2}
	\begin{align}
		-\nabla \cdot \left(A^{n+1/2}\nabla U\right)+\nabla\cdot\left(B^{n+1/2} P\right)&=C^{n+1/2}F_h \quad \mbox{in $\Omega^0$},\\
		B^{n+1/2}:\nabla U&=0 \quad \mbox{in $\Omega^0$},
	\end{align}
\end{subequations}
with 
\begin{align*}
	& A^{n+1}=\left({\rm det}\nabla\phi^{n+1}\right)\left(\nabla \phi^{n+1}\right)^{-1}\left(\nabla \phi^{n+1}\right)^{-\top}, \quad 	B^{n+1}=\left({\rm det}\nabla\phi^{n+1}\right)\left(\nabla \phi^{n+1}\right)^{-\top},\\
	& C^{n+1}=\left({\rm det}\nabla\phi^{n+1}\right).
\end{align*}
Since $\Omega^0$ has smooth boundary and $A^{n+1/2}, B^{n+1/2}, C^{n+1/2}$ have smooth coefficients, we have elliptic regularity estimate
\begin{align*}
	\|U\|_{H^2(\Omega_h^0)}+\|P\|_{H^1(\Omega^0)}\leq C\|F_h\|_{\Omega_h^0}
\end{align*}
We extend $U,P$ to $\widetilde{U}, \widetilde{P}$ by Stein's extension operator and define an auxiliary function $\widetilde{\eta}$ as:
\begin{align*}
	\widetilde{\eta}:=-\nabla \cdot \left(A^{n+1/2}\nabla \widetilde{U}\right)+\nabla\cdot\left(B^{n+1/2} \widetilde{P}\right)-C^{n+1/2}F_h.
\end{align*}
Then $\widetilde{\eta}|_{\Omega^0}=0$ and $\|\widetilde{\eta}\|_{L^2(\mathbb{R}^d)}\leq C\|F_h\|_{\Omega_h^0}$. By integrating $\widetilde{\eta}$ against $W_h$ over $\Omega_h^0$, we obtain
\begin{subequations}
	\begin{align*}
		(A^{n+1/2}\nabla \widetilde{U}, \nabla W_h)_{\Omega_h^{0}}-(B^{n+1/2}:\nabla W_h, \widetilde{P})_{\Omega_h^{0}}&=(C_h^{n+1/2}F_h, W_h)_{\Omega_h^{0}}+(\widetilde{\eta}, W_h)_{\Omega_h^0\setminus\Omega^0}\; \forall W_h\in \mathring V_h,\\
		(B^{n+1/2}:\nabla \widetilde{U}, Y_h)_{\Omega_h^0}&=(B^{n+1/2}:\nabla \widetilde{U}, Y_h)_{\Omega_h^0\setminus \Omega^0}\; \forall Y_h\in Q_h.
	\end{align*}
\end{subequations}
Taking into account the fact that
\begin{align*}
	\|A_h^{n+1/2}-A^{n+1/2}\|_{L^p(\Omega_h^0)}+ \|B_h^{n+1/2}-B^{n+1/2}\|_{L^p(\Omega_h^0)}+\|C_h^{n+1/2}-C^{n+1/2}\|_{L^p(\Omega_h^0)}\leq Ch,
\end{align*}
by Lemma \ref{boundary-skin} we obtain 
\begin{subequations}
	\begin{align*}
		(A_h^{n+1/2}\nabla (\widetilde{U}-U_h), \nabla W_h)_{\Omega_h^{0}}-(B^{n+1/2}:\nabla W_h, \widetilde{P}-P_h)_{\Omega_h^{0}}&=\ell(W_h)\quad \forall W_h\in \mathring V_h,\\
		(B^{n+1/2}:\nabla (\widetilde{U}-U_h), Q_h)_{\Omega_h^0}&=\phi(Y_h)\quad \forall Y_h\in Q_h,
	\end{align*}
\end{subequations}
with $ \ell $ and $ \phi $ satisfying the following estimates
\begin{align*}
	|\ell(W_h)|\leq Ch\|F_h\|_{\Omega_h^0}\|\nabla W_h\|_{\Omega_h^0}, \quad |\phi(Y_h)|\leq Ch\|F_h\|_{\Omega_h^0}\|Y_h\|_{\Omega_h^0}.
\end{align*}
Analogous to Lemma \ref{H^1-error-1}, we can deduce that there is some constant $\lambda$ such that
\begin{align*}
	\|\widetilde{P}-P_h-\lambda\|_{\Omega_h^0}\leq Ch\|F_h\|_{\Omega_h^0}.
\end{align*}
Finally, through the application of inverse inequality, we can deduce
\begin{align*}
	\|\nabla P_h\|_{\Omega_h^0}\leq C\|F_h\|_{\Omega_h^0},
\end{align*}
which concludes the proof.
\end{proof}


\begin{thebibliography}{10}
\bibitem{Badia2006}
{S. Badia and R. Codina.} {Analysis of a stabilized finite element approximation of the transient convection-diffusion equation using an ALE framework.} {\em SIAM J. Numer. Anal.}, 44(5): 2159--2197, 2006.

\bibitem{boffi13}
{D. Boffi, F. Brezzi and M. Fortin.} {\em Mixed finite element methods and applications}. Vol. 44. Heidelberg: Springer, 2013.

\bibitem{Boffi2004}
{D. Boffi and L. Gastaldi.} {Stability and geometric conservation laws for ALE formulations.} {\em Comput. Methods Appl. Mech. Eng.}, 193(42-44): 4717--4739, 2004.

\bibitem{Bonito2013}
{A. Bonito, I. Kyza and R. H. Nochetto.} {Time-discrete higher-order ALE formulations: stability.} {\em SIAM J. Numer. Anal.}, 51(1): 577--604, 2013.

\bibitem{BonKyerror2013}
{A. Bonito, I. Kyza and R. H. Nochetto}. {Time-discrete higher order ALE formulations: A priori error analysis.} {\em Numer. Math.}, 125(2): 225--257, 2013.

\bibitem{Dziuk-Elliott-2007}
G. Dziuk and C. Elliott: Finite elements on evolving surfaces. {\em IMA J. Numer. Anal.} 27 (2007), pp. 262--292. 

\bibitem{Edelmann-2022}
D. Edelmann: 
Finite element analysis for a diffusion equation on a harmonically evolving domain. 
{\em IMA J. Numer. Anal.} 42 (2022), pp. 1866--1901.

\bibitem{ElliottUFEM2021}
{C. M. Elliott and T. Ranner.} {A unified theory for continuous-in-time evolving finite element space approximations to partial differential equations in evolving domains.} {\em IMA J. Numer. Anal.}, 41(3): 1696--1845, 2021.

\bibitem{EllioESFEM2015}
C. M. Elliott and C. Venkataraman. {Error analysis for an ALE evolving surface finite element method.} {\em Numer. Methods Partial Differ. Equ.}, 31(2): 459--499, 2015.

\bibitem{ForNo2004}
{L. Formaggia, F. Nobile.} {Stability analysis of second-order time accurate schemes for ALE–FEM.} {\em Comput. methods Appl. Mech. Eng.}, 193(39-41): 4097--4116, 2004.

\bibitem{farwig94}
{R. Farwig and H. Sohr}. {Generalized resolvent estimates for the Stokes system in bounded and unbounded domains}. {\em J. Math. Soc. Japan}, 46.4: 607--643, 1994.

\bibitem{Gas2001}
{L. Gastaldi.} {A priori error estimates for the arbitrary Lagrangian Eulerian formulation with finite elements.} 2001.

\bibitem{GawE2015}
{E. S. Gawlik and A. J. Lew.} {Unified analysis of finite element methods for problems with moving boundaries.} {\em SIAM J. Numer. Anal.}, 53(6): 2822--2846, 2015.

\bibitem{Girault2012}
V. Girault and P. A. Raviart. {\em Finite element methods for Navier--Stokes equations: theory and algorithms}. Springer Science and Business Media, 2012.

\bibitem{GongLiRao-2023}
W. Gong, B. Li, and Q. Rao. Convergent evolving finite element approximations of boundary evolution under shape gradient flow. {\em IMA J. Numer. Anal.}, 2023, DOI: 10.1093/imanum/drad080 

\bibitem{KLLP-2017}
{B. Kov\'acs, B. Li, C. Lubich and C. A. Power Guerra}. {Convergence of finite elements on an evolving surface driven by diffusion on the surface}. {\em Numer. Math.} 137 (2017), pp. 643--689.

%
\bibitem{Legendre2008}
{G. Legendre and T. Takahashi.} {Convergence of a Lagrange-Galerkin method for a fluid-rigid body system in ALE formulation.} {\em ESAIM Math. Model. Numer. Anal.}, 42(4): 609--644, 2008.

\bibitem{Lenoir1986}
{M.~Lenoir}. Optimal isoparametric finite elements and error estimates for domains 
involving curved boundaries. {\em SIAM J. Numer. Anal.}, 23(3):562--580, 1986.

\bibitem{LiMA-2022}
B. Li, S. Ma and Y. Ueda. {Analysis of Fully Discrete Finite Element Methods for $2D$ Navier--Stokes Equations with Critical Initial Data}. {\em ESAIM:M2AN}, 56: 2105--2139, 2022.

\bibitem{buyang22}
{B. Li, Y. Xia and Z. Yang}. Optimal convergence of arbitrary Lagrangian–Eulerian iso-parametric finite element methods for parabolic equations in an evolving domain. {\em IMA J. Numer. Anal.}, 43(1):501--534, 2023.

\bibitem{Liu2013}
{J. Liu}. {Simple and efficient ALE methods with provable temporal accuracy up to fifth order for the Stokes equations on time varying domains.} {\em SIAM J. Numer. Anal.}, 51(2): 743--772, 2013.

\bibitem{SanMart2009}
{J. S. Martín, L. Smaranda and T. Takahashi.} {Convergence of a finite element/ALE method for the Stokes equations in a domain depending on time.} {\em J. Comput. Appl. Math.}, 230(2): 521--545, 2009.

\bibitem{Nobie2001}
{F. Nobile.} {Numerical approximation of fluid-structure interaction problems with application to haemodynamics.} {\em Ph.D. thesis}, Department of Mathematics, Ecole Polytechnique Fdrale de Lausanne, Switzerland, 2001.

\bibitem{Formg1999}
{F. Nobile and L. Formaggia.} {A stability analysis for the arbitrary Lagrangian Eulerian formulation with finite elements.} {\em East-West J. Numer. Math.}, 7(ARTICLE): 105--132, 1999.

\bibitem{ShenJ1996}
{J. Shen.} {On error estimates of the projection methods for the Navier--Stokes equations: second-order schemes.} {\em Math. Comp.}, 65(215): 1039--1065, 1996.

\bibitem{Walker2015}
{S. W. Walker.} {The shapes of things: a practical guide to differential geometry and the shape derivative.} {\em SIAM}, 2015.


































\end{thebibliography}
\end{document}